\documentclass[reqno,a4paper,11pt]{amsart}
\usepackage[english]{babel}
\renewcommand{\baselinestretch}{1.1}
\usepackage{graphicx,a4wide,mathrsfs,tikz,mathtools,latexsym,ifthen,amssymb,stmaryrd}
\mathtoolsset{showonlyrefs,showmanualtags}
\usepackage{enumerate}
\usepackage{times}

\usepackage[utf8x]{inputenc}
\usepackage[toc,page]{appendix}

\usepackage{color}
\definecolor{MyDarkBlue}{rgb}{0.15,0.25,0.45}
\usepackage[pdftex]{hyperref}
\hypersetup{
colorlinks=true,
citecolor=red,
linkcolor=MyDarkBlue,
urlcolor=MyDarkBlue,
pdfauthor={Ugo Bruzzo and Francesco Sala},
pdftitle={Framed sheaves on projective stacks},
pdfsubject={}
breaklinks=true
}

\newcommand{\hilb}[2]{\mathrm{Hilb}^{#1}(#2)}

\newcommand{\triend}{\parbox{2mm}{\hfill} \hfill\mbox{\hspace{0.2mm}}\hfill$\triangle$}
\newcommand{\ocend}{\parbox{2mm}{\hfill} \hfill\mbox{\hspace{0.2mm}}\hfill$\oslash$}

\newtheorem{theorem}{Theorem}[section]
\newtheorem{proposition}[theorem]{Proposition}
\newtheorem{lemma}[theorem]{Lemma}
\newtheorem{corollary}[theorem]{Corollary}
\newtheorem*{theorem*}{Theorem}

\theoremstyle{remark}
\newtheorem{exa}[theorem]{Example}
\newenvironment{example}{\begin{exa}}{\triend\end{exa}}

\theoremstyle{remark}
\newtheorem{rem}[theorem]{Remark}
\newenvironment{remark}{\begin{rem}}{\triend\end{rem}}

\theoremstyle{definition}
\newtheorem{defin}[theorem]{Definition}
\newenvironment{definition}{\begin{defin}}{\ocend\end{defin}}

\newtheorem{condition}[theorem]{Condition}

\title[Framed sheaves on projective stacks]{\large Framed sheaves on projective stacks}
\begin{document}
\bigskip
\maketitle \thispagestyle{empty} \vspace{-5mm} 
\begin{center}
{\sc Ugo Bruzzo$^{\S\ddag}$ and Francesco Sala$^{\P\star}$} \\[3pt] 
with an appendix by {\sc Mattia Pedrini}$^{\S\ddag}$ \\[10pt] 
\small
$^\S$ Scuola Internazionale Superiore di Studi Avanzati {\sc (SISSA)},\\ Via Bonomea 265, 34136
Trieste, Italia \\[4pt] $^\ddag$ Istituto Nazionale di Fisica Nucleare, Sezione di Trieste \\[4pt]
$^\P$ Heriot-Watt University,  
School of Mathematical and Computer Sciences, \\ Department of Mathematics, 
Colin Maclaurin Building, Riccarton, \\ Edinburgh EH14 4AS, United Kingdom \\[4pt]  $^\star$ Maxwell Institute for Mathematical Sciences, Edinburgh, United Kingdom
 \end{center}
 \par
\vfill
 {\small {\sc Abstract.}
Given a normal projective irreducible stack $\mathscr X$ over an algebraically closed field of characteristic zero  we consider {\em framed sheaves} on $\mathscr X$, i.e., pairs $(\mathcal E,\phi_{\mathcal E})$, where $\mathcal E$ is a coherent sheaf on $\mathscr X$ and $\phi_{\mathcal E}$ is a morphism from $\mathcal E$ to a fixed coherent sheaf $\mathcal F$. 
After introducing a suitable notion of (semi)stability, we construct a projective scheme, which is a moduli space for semistable framed sheaves with fixed Hilbert polynomial, and an open subset of it, which is a fine moduli space for stable framed sheaves. If $\mathscr X$ is a projective irreducible orbifold of dimension two and $\mathcal F$ a locally free sheaf on a smooth divisor $\mathscr D\subset \mathscr X$ satisfying certain conditions, we consider {\em $(\mathscr{D}, \mathcal{F})$-framed sheaves}, i.e., framed sheaves $(\mathcal E,\phi_{\mathcal E})$ with $\mathcal E$ a torsion-free sheaf which is locally free in a neighborhood of $\mathscr D$, and ${\phi_{\mathcal{E}}}_{\vert \mathscr{D}}$   an isomorphism. These pairs are $\mu$-stable for a suitable choice of a parameter entering the (semi)stability condition, and of the polarization of $\mathscr X$. This implies the existence of a fine moduli space parameterizing isomorphism classes of $(\mathscr{D}, \mathcal{F})$-framed sheaves on $\mathscr{X}$ with fixed Hilbert polynomial, which is a quasi-projective scheme. In an appendix we develop the example of stacky Hirzebruch surfaces.
This is the first paper of a project  aimed to provide an algebro-geometric approach to the study of gauge theories on a wide class of 4-dimensional Riemannian manifolds by means of framed sheaves on ``stacky" compactifications of them. In particular, in a subsequent paper \cite{art:bruzzopedrinisalaszabo2013}  these results are used  to study gauge theories on ALE spaces of type $A_k$. 
} 

\vfill
\parbox{.94\textwidth}{\hrulefill}\par
\noindent \begin{minipage}[c]{\textwidth}\parindent=0pt \renewcommand{\baselinestretch}{1.2}
\small
{\em Date:}   \today  \par 
\par
{\em E-Mail:} {\tt ugo.bruzzo@sissa.it, mattia.pedro@gmail.com, salafra83@gmail.com}

The authors gratefully acknowledge financial support and hospitality during the respective visits to Heriot-Watt University and {\sc sissa}. Support for this work was partly provided by {\sc prin} ``Geometry of Algebraic Varieties''. The second author was supported in part by Grant RPG-404 from
the Leverhulme Trust. The first author is a member of the {\sc vbac} group.

Current address of the second author: Department of Mathematics, Middlesex College, The University of Western Ontario. London, Ontario, Canada, N6A 5B7.
\end{minipage}

\newpage

{\small\tableofcontents}

%%%%%%

\section{Introduction}
 
According to Donaldson \cite{art:donaldson1984}, the moduli space of gauge equivalence classes of framed $SU(r)$-instantons of charge $n$ on $\mathbb R^4$ is isomorphic to the moduli space $\mathrm{M}^{reg}(r,n)$ of isomorphism classes of {\em framed vector bundles} of rank $r$ and second Chern class $n$ on the complex projective plane $\mathbb{P}^2=\mathbb{C}^2\cup l_\infty$. A framed vector bundle is a pair $(E,\phi_E)$, where $E$ is a vector bundle on $\mathbb{P}^2$ of rank $r$ and $c_2(E)=n$ and $\phi_E\colon E_{\vert l_{\infty}}\xrightarrow{\sim}\mathcal{O}_{l_{\infty}}^{\oplus r}$ a trivialization along the line $l_{\infty}$. The moduli space $\mathrm{M}^{reg}(r,n)$ is an open subset in the moduli space $\mathrm{M}(r,n)$ of \textit{framed sheaves} on $\mathbb{P}^2$, i.e., the moduli space of pairs $(E,\phi_E)$ modulo isomorphism, with $E$ a torsion-free sheaf on $\mathbb{P}^2$ of rank $r$ and $c_2(E)=n$, locally free in a neighbourhood of $l_{\infty}$, and $\phi_E\colon E_{\vert l_{\infty}}\xrightarrow{\sim}\mathcal{O}_{l_{\infty}}^{\oplus r}$ a \textit{framing at infini\-ty}. $\mathrm{M}(r,n)$ is a nonsingular quasi-projective variety of dimension $2rn$, which has a description in terms of {\em monads} and linear data, the so-called {\em Atiyah-Drinfeld-Hitchin-Manin (ADHM) data} \cite{book:nakajima1999}, Ch.~2. It is a resolution of singularities of the moduli space 
\begin{equation}
\mathrm{M}_0(r,n):=\bigsqcup_{i=0}^n \mathrm{M}^{reg}(r,i)\times S^{n-i}(\mathbb{C}^2)
\end{equation}
of {\em ideal framed instantons} on $\mathbb R^4$, i.e., instantons whose square curvature is allowed to degenerate to a Dirac delta at a number $m$ of points, with $ 1 \le m \le n$. Here we denote by $S^{n-i}(\mathbb{C}^2)$ the $n-i$ symmetric product of $\mathbb{C}^2$. The moduli spaces $\mathrm{M}(1,n)$ coincide with the Hilbert schemes $(\mathbb C^2)^{[n]}$ of $n$ points in the complex affine plane. 

In view of its relation with framed instantons, the moduli space $\mathrm{M}(r,n)$ has been studied quite intensively (see, e.g., \cite{art:bruzzofucitomoralestanzini2003, book:nakajima1999, art:nakajimayoshioka2005-I, art:nakajimayoshioka2005-II}) and its geometry is quite well known. This is in particular important in physics, where the moduli spaces of ideal  instantons play the role of  parameter spaces  for the classical vacua of a (topological, supersymmetric) Yang-Mills theory; unfortunately this space is singular,  and therefore, in order to make computations, it is conveniently replaced by the the moduli space of framed sheaves, with the added bonus of being able to use powerful algebro-geometric techniques. The so-called Nekrasov partition function  introduced in \cite{art:nekrasov2003} (see also \cite{art:flumepoghossian2003,art:bruzzofucitomoralestanzini2003}) plays an important role in this theory, also in view  of its connection with the  {\em Seiberg-Witten prepotential} (see, e.g., \cite{art:nakajimayoshioka2005-I, art:nakajimayoshioka2005-II}) and  with  {\em Donaldson invariants} (see, e.g., \cite{art:gottschenakajimayoshioka2008}).

It is quite natural to wonder if Donaldon's correspondence can be generalized to other noncompact 4-dimensional Riemannian manifolds. King in his PhD thesis \cite{phd:king1989} studied the correspondence between framed instantons on the blowup $\widetilde{\mathbb{C}^2}$ of $\mathbb{C}^2$ at the origin and framed vector bundles on the blowup $\widetilde{\mathbb{P}^2}=\widetilde{\mathbb{C}^2}\cup l_\infty$ of $\mathbb P^2=\mathbb{C}^2\cup l_\infty$ at the same point; this was generalized by Buchdhal \cite{art:buchdahl1993} by considering framed instantons on the blowup of $\mathbb{C}^2$ at $n$ points and framed vector bundles on the blowup of $\mathbb P^2$ at the same points.

All the examples so far described share one pattern: in order to study instantons on a noncompact 4-dimensional Riemannian manifold, one endows the manifold of a K\"ahler structure, compactifies it by adding a projective line, and considers framed vector bundles on the corresponding smooth projective surface. A natural question to ask is if this algebro-geometric approach to the study of Yang-Mills theories holds for other  4-dimensional Riemannian manifolds. First, one should check if there is a consistent theory of moduli spaces parameterizing framed vector bundles (sheaves) with fixed Chern classes on any smooth projective surface. This is the case: indeed Bruzzo and Markushevich provide in \cite{art:bruzzomarkushevich2011} a general construction of a fine moduli space for framed sheaves on smooth projective surfaces, building on work by Huybrechts and Lehn \cite{art:huybrechtslehn1995-I,art:huybrechtslehn1995-II}. Moduli spaces of framed sheaves on multiple blowups of $\mathbb P^2$ were considered by Henni \cite{art:henni2009}, while the case of Hirzebruch surfaces were constructed by Bartocci, Bruzzo and Rava in \cite{art:bartoccibruzzorava2013}, in both cases using {\em monads}. Secondly, one needs a generalization of Donaldson's correspondence to (at least a wide class of) noncompact 4-dimensional Riemannian manifolds. 

In this connection, one should keep in mind that the ways a  K\"ahler four-dimensional manifold $M$ can be suitably compactified is  constrained by a result of Bando \cite{art:bando1993}.  According to Bando, if $\bar M$ is a K\"ahler compactification of $M$ by a smooth divisor $D$ whose normal line bundle is positive, holomorphic vector bundles on $\bar M$, which are isomorphic along the compactifying divisor $D$ to a fixed vector bundle endowed with a flat connection $\nabla$, correspond to holomorphic vector bundles on $M$ with anti-selfdual square integrable connections (i.e. instantons), with holonomy at infinity induced by $\nabla$.
 If $D$ is a projective line, only instantons on $M$ with trivial holonomy at infinity can be
described in terms of framed locally free sheaves on $\bar M$.  For example, this case occurs when $M$ is an ALE space.

To circumvent this restriction one can change the ambient space from a compact K\"ahler surface to a 2-dimensional compact K\"ahler Deligne-Mumford stack. Let $\mathscr{X}$ be a 2-dimensional connected compact K\"ahler Deligne-Mumford stack and $\mathscr{D}$ a smooth 1-dimensional integral closed substack of $\mathscr{X}$ such that the line bundle $\mathcal{O}_{\mathscr{X}}(\mathscr{D})$ is positive and ample on $\mathscr{D}$ and $X_0=\mathscr{X}\setminus \mathscr{D}$ is a K\"ahler surface with cone-like singularities (for example, an ALE space).   Eyssidieux and the second author proved in \cite{art:eyssidieuxsala2013} that there is a correspondence between holomorphic vector bundles on $\mathscr{X}$, which are isomorphic along $\mathscr{D}$ to a fixed vector bundle $\mathcal F$, and holomorphic vector bundles on $X_0$ endowed with Hermite-Einstein metrics with holonomy at infinity given by a fixed flat connection on $\mathcal F$.

Evidence for this kind of generalization of Donaldson's correspondence can be found in \cite{art:bruzzopoghossiantanzini2011}. There Bruzzo, Poghossian and Tanzini computed the partition function of $\mathcal{N}=4$ supersymmetric Yang-Mills theories on the total spaces $\mathrm{Tot}(\mathcal O_{\mathbb P^1}(-p))$ of the line bundles $\mathcal O_{\mathbb P^1}(-p)$ by means of framed sheaves on the Hirzebruch surfaces $\mathbb F_p$, regarded as projective compactifications of $\mathrm{Tot}(\mathcal O_{\mathbb P^1}(-p))$. That analysis showed that the computation  {of the partition function} made sense also for framed sheaves $(E,\phi_E)$ on $\mathbb F_p$ with first Chern class $c_1(\mathcal E) = \frac{k}{p}C$, with $k$ any integer, and $C$ the class of the section of $\mathbb F_p\to\mathbb P^1$ squaring to $-p$. This of course makes little sense, and indeed in \cite{art:bruzzopoghossiantanzini2011} the authors conjectured that the computation actually was taking place on a ``stacky'' compactification of $\mathrm{Tot}(\mathcal O_{\mathbb P^1}(-p))$. Indeed, as we shall see in  Appendix \ref{Mattia}, a fractional first Chern class  only appears  when one considers instantons with nontrivial holonomy at infinity.

In this paper we construct a mathematically rigorous theory of moduli spaces of framed sheaves on projective stacks. Together with the work of Eyssidieux and the second author \cite{art:eyssidieuxsala2013}, this paper provides a completely algebro-geometric approach to the study of instantons on noncompact 4-dimensional Riemannian manifolds which can be compactified to projective orbifolds by adding one-dimensional smooth projective stacks. This is applied in \cite{art:bruzzopedrinisalaszabo2013} to study instantons on {\em ALE spaces of type $A_k$} \cite{art:kronheimernakajima1990}, with $k\geq 1$, by means of framed sheaves on 2-dimensional projective toric orbifolds, which are ``stacky" compactifications of the ALE spaces (in the complex analytic setting, the idea of compactifying the ALE spaces to complex V-manifolds was already suggested in \cite{art:nakajima2007}). In particular,   \cite{art:bruzzopedrinisalaszabo2013}  provides  a rigorously mathematical derivation of the partition functions for gauge theories on ALE spaces of type $A_k$, conjecturally described in \cite{art:bonellimaruyoshitanzini2011-I,art:bonellimaruyoshitanzini2011-II, art:bonellimaruyoshitanziniyagi2012}. Moreover, the study of the partition functions in \cite{art:bruzzopedrinisalaszabo2013} allows a comparison with the partition functions computed in \cite{art:fucitomoralespoghossian2004},  clarifying some ambiguities noticed in \cite{art:itomaruyoshiokuda2013} between these two different approaches to compute partition functions for gauge theories on ALE spaces of type $A_k$.

We state now  the main results of this paper. Let $\mathscr X$ be a normal projective irreducible stack of dimension $d$ defined over an algebraically closed field $k$ of characteristic zero, with a coarse moduli scheme $\pi\colon\mathscr X\to X$. Fix a polarization $(\mathcal G,\mathcal O_X(1))$ on $\mathscr X$. A \textit{framed sheaf} is a pair $(\mathcal E,\phi_{\mathcal E})$ where $\mathcal E$ is a coherent sheaf on $\mathscr X$ and $\phi_{\mathcal E}$ is a morphism from $\mathcal E$ to a fixed coherent sheaf $\mathcal F$. We  call $\phi_{\mathcal{E}}$ a {\em framing} of $\mathcal{E}$. This notion is more general than the one  discussed before;  framed sheaves for which $\mathcal{F}$ is a locally free sheaf over a divisor $\mathscr{D}$ and the framing is an isomorphism will be called $(\mathscr{D},\mathcal{F})$-framed sheaves.

We consider a generalization of Gieseker (semi)stability for framed sheaves that depends on the polarization and on a rational polynomial $\delta$ of degree $d-1$ with positive leading coefficient. We call it $\delta$-(semi)stability condition to emphasize the dependence on $\delta$.

Fix a numerical polynomial $P$ of degree $d$. Let 
\begin{equation}
\underline{\mathcal{M}}^{(s)s}\colon (Sch/k)^{\circ}\to (Sets)
\end{equation}
be the contravariant functor of $\delta$-(semi)stable framed sheaves on $\mathscr{X}$ with Hilbert polynomial $P$, which associates with any scheme $S$ of finite type over $k$ the set of isomorphism classes of flat families of $\delta$-(semi)stable framed sheaves on $\mathscr{X}$ with Hilbert polynomial $P$ parameterized by $S$. The first problem we have addressed in this paper is the study of the (co)representability of this functor (which is equivalent to ask if there exists a (fine) moduli space of $\delta$-(semi)stable framed sheaves on $\mathscr{X}$ with Hilbert polynomial $P$). We obtain the following result.
\begin{theorem}\label{thm:uno}
Let $\mathscr X$ be a normal projective irreducible stack of dimension $d$ defined over an algebraically closed field $k$ of characteristic zero, with a coarse moduli scheme $\pi\colon\mathscr X\to X$, and $(\mathcal G,\mathcal O_X(1))$ a polarization on $\mathscr X$. Fix a coherent sheaf $\mathcal{F}$ on $\mathscr{X}$ and a rational polynomial $\delta$ of degree $d-1$ with positive leading coefficient $\delta$. Then for any numerical polynomial $P$ of degree $d$, there exists an algebraic stack $\mathfrak{M}^{(s)s}$ of finite type over $k$ such that
\begin{itemize}\setlength{\itemsep}{2pt}
\item $\mathfrak{M}^{(s)s}$ admits a good moduli space $\pi\colon \mathfrak{M}^{(s)s}\to \mathrm{M}^{(s)s}$ (in the sense of Alper);
\item $\mathrm{M}^{ss}$ is a projective scheme and $\mathrm{M}^{s}$ is an open subscheme of $\mathrm{M}^{ss}$;
\item the contravariant functor $[\mathfrak{M}^{(s)s}]$ which associates with any scheme $S$ of finite type over $k$ the set of isomorphism classes of objects in $\mathfrak{M}^{(s)s}(S)$ is isomorphic to the moduli functor $\underline{\mathcal{M}}^{(s)s}$ of $\delta$-(semi)stable framed sheaves on $\mathscr{X}$ with Hilbert polynomial $P$;
\item $\mathrm{M}^{ss}$ corepresents the contravariant functor $[\mathfrak{M}^{ss}]$, while $\mathrm{M}^{s}$ represents the contravariant functor $[\mathfrak{M}^{s}]$.
\end{itemize}
\end{theorem}
Roughly speaking, the previous Theorem states that one can construct a projective scheme $\mathrm{M}^{ss}$, which is a moduli space for $\delta$-semistable framed sheaves on $\mathscr{X}$ with Hilbert polynomial $P$, and a quasi-projective scheme $\mathrm{M}^{s}$, which is a fine moduli space for $\delta$-stable framed sheaves on $\mathscr{X}$ with Hilbert polynomial $P$, i.e., it comes with a {\em universal} family of $\delta$-stable framed sheaves on $\mathscr{X}$ with Hilbert polynomial $P$.
 
As in the case of framed sheaves on smooth projective varieties, we have the following characterization of the tangent space of $\mathrm{M}^s$.
\begin{theorem}\label{thm:uno-bis}
The tangent space to the fine moduli space $\mathrm{M}^s$ at a point $[(\mathcal E,\phi_{\mathcal E})]$ can be identified with the hyper-Ext group $\mathbb E \mathrm{xt}^1(\mathcal E, \mathcal E \xrightarrow{\phi_{\mathcal E}}\mathcal F)$, while the hyper-Ext group $\mathbb E \mathrm{xt}^2(\mathcal E, \mathcal E \xrightarrow{\phi_{\mathcal E}}\mathcal F)$ contains the obstruction to the smoothness of $\mathrm{M}^s$ at the point $[(\mathcal E,\phi_{\mathcal E})]$.
\end{theorem}
 
Let $\mathscr X$ be a projective irreducible orbifold of dimension two. Extending the original definition of framed sheaves on $\mathbb{P}^2$, we shall consider framings along a fixed 1-dimensional smooth integral closed substack $\mathscr D$, whose coarse moduli space $D$ is a $\mathbb Q$-Cartier big and nef smooth curve. Let $\mathcal F$ be a coherent sheaf on $\mathscr X$, supported on $\mathscr D$, such that $\mathcal F$ is a locally free $\mathcal{O}_{\mathscr D}$-module satisfying a suitable semistability condition.  We shall call $D$ a {\em good framing divisor} and $\mathcal F$ a {\em good framing sheaf}. A $(\mathscr D,\mathcal F)$\textit{-framed sheaf} is a framed torsion-free sheaf $(\mathcal E,\phi_{\mathcal E}\colon \mathcal E\to  \mathcal F)$, with $\mathcal E$ locally free in a neighbourhood of $\mathscr D$ and ${\phi_{\mathcal E}}_{\vert \mathscr D}$ an isomorphism.

One can introduce a $\mu$-(semi)stability condition for framed sheaves on $\mathscr{X}$ depending on a positive rational number $\delta$. With a suitable choice of the parameter $\delta_1$ and of the polarization of $\mathscr X$, one obtains that all $(\mathscr D,\mathcal F)$-framed sheaves are $\mu$-stable, hence one has the following result.
\begin{theorem} \label{thm:due} 
There exists a fine moduli space parameterizing isomorphism classes of $(\mathscr{D}, \mathcal{F}_{\mathscr{D}})$-framed sheaves $(\mathcal E,\phi_{\mathcal E})$ on $\mathscr{X}$ with given Hilbert polynomial $P$, which is a quasi-projective scheme. Its tangent space at a point $[(\mathcal E,\phi_{\mathcal E})]$ is $\operatorname{Ext}^1_{\mathscr{X}}(\mathcal E,\mathcal E(-\mathscr{D}))$. If $\operatorname{Ext}^2_{\mathscr{X}}(\mathcal E,\mathcal E(-\mathscr{D}))=0$ for all points $[(\mathcal E,\phi_{\mathcal E})]$, the moduli space is a smooth quasi-projective variety.
\end{theorem}

In the last part of the paper (Section \ref{FramedsheavesToricorbifolds}) we apply the theory  to two-dimensional projective irreducible orbifolds that are toric root stacks. Let $X$ be a normal projective toric surface and $D$ a torus-invariant rational curve which contains the singular locus $\mathrm{sing}(X)$ of $X$ and   is a good framing divisor. Let $\pi^{can}\colon\mathscr{X}^{can}\to X$ be the {\em canonical toric orbifold} of $X$. It is the unique (up to isomorphism) two-dimensional toric orbifold for which the locus where $\pi^{can}$ is not an isomorphism has a nonpositive dimension. Denote by $\tilde{\mathscr{D}}$ the smooth effective Cartier divisor $(\pi^{can})^{-1}(D)_{red}$. Let $\mathscr{X}:=\sqrt[k]{\tilde{\mathscr{D}}/\mathscr{X}^{can}}$ be the toric orbifold over $\mathscr{X}^{can}$ obtained by performing a {\em $k$-th root construction along $\tilde{\mathscr{D}}$}. The stack $\mathscr{X}$ is the fibred product of $\mathscr{X}^{can}\times_{[\mathbb{A}^1/\mathbb{G}_m]} [\mathbb{A}^1/\mathbb{G}_m]$, where the morphism $\mathscr{X}^{can}\to [\mathbb{A}^1/\mathbb{G}_m]$ is induced by $\tilde{\mathscr{D}}$ and the morphism $[\mathbb{A}^1/\mathbb{G}_m]\to [\mathbb{A}^1/\mathbb{G}_m]$ is induced by the map sending a complex number to its $k$-th power. The induced natural morphism $\mathscr{X}\to [\mathbb{A}^1/\mathbb{G}_m]$ corresponds to a smooth effective Cartier divisor $\mathscr{D}$. Away from $\tilde{\mathscr{D}}$ the stacks $\mathscr{X}$ and $\mathscr{X}^{can}$ are isomorphic and  $\mathscr{D}$ is an \'etale $\mu_k$-gerbe over $\tilde{\mathscr{D}}$; so --- roughly speaking --- if we endow locally the stack $\mathscr{X}^{can}$ of a $\mu_k$-action along $\tilde{\mathscr{D}}$, the divisor $\tilde{\mathscr{D}}$ is globally replaced by a $\mu_k$-gerbe over itself, and we obtain $\mathscr{X}$. The next theorem states that if $D$ is in addition a good framing divisor and the line bundle $\mathcal{O}_{\mathscr{X}^{can}}(\tilde{\mathscr{D}})$ is $\pi^{can}$-ample, there exist fine moduli spaces for $(\mathscr{D}, \mathcal{F})$-framed sheaves on $\mathscr{X}$ with fixed Hilbert polynomial for any choice of the good framing sheaf $\mathcal{F}$.
\begin{theorem}\label{thm:tre}
Let $X$ be a normal projective toric surface and $D$ a torus-invariant rational curve which contains the singular locus $\mathrm{sing}(X)$ of $X$ and is a good framing divisor. Let $\pi^{can}\colon\mathscr{X}^{can}\to X$ be the canonical toric orbifold of $X$ and $\tilde{\mathscr{D}}$ the smooth effective Cartier divisor $(\pi^{can})^{-1}(D)_{red}$. Assume that the line bundle $\mathcal{O}_{\mathscr{X}^{can}}(\tilde{\mathscr{D}})$ is $\pi^{can}$-ample. Let $\mathscr{X}:=\sqrt[k]{\tilde{\mathscr{D}}/\mathscr{X}^{can}}$, for some positive integer $k$, and $\mathscr{D}\subset\mathscr{X}$ the effective Cartier divisor corresponding to the morphism $\mathscr{X}\to [\mathbb{A}^1/\mathbb{G}_m]$. Then for any good framing sheaf $\mathcal{F}$ on $\mathscr{D}$ and any numerical polynomial $P\in \mathbb{Q}[n]$ of degree two, there exists a fine moduli space parameterizing isomorphism classes of $(\mathscr{D}, \mathcal{F})$-framed sheaves on $\mathscr{X}$ with Hilbert polynomial $P$, which is a quasi-projective scheme over $\mathbb{C}$.
\end{theorem}
  
This paper is organized as follows. In Section \ref{Projective stacks} we give the definitions of projective stacks and of polarizations on them, and describe the notion of support, purity and  Hilbert polynomial of a coherent sheaf on them. In Section \ref{Framed sheaves on projective stacks}, by generalizing \cite{art:huybrechtslehn1995-I,art:huybrechtslehn1995-II}, we define the notion of framed sheaf and the related (semi)stability conditions. Moreover, we give a notion of flat family of framed sheaves, following \cite{art:bruzzomarkushevichtikhomirov2010}, Sect.~2, and prove a boundedness theorem for flat families of $\delta$-semistable framed sheaves. In Section \ref{Moduli spaces of framed sheaves on projective stacks}, by using the GIT machinery, we prove Theorems \ref{thm:uno} and \ref{thm:uno-bis}. The proofs use some arguments from \cite{art:nironi2008,art:huybrechtslehn1995-I,art:huybrechtslehn1995-II}. In Section \ref{sec:framedsheaves} we introduce the notion of $(\mathscr D,\mathcal F)$-framed sheaf on projective orbifolds of dimension two, give a boundedness result, and prove Theorem \ref{thm:due}. As a byproduct, we obtain a boundedness result for $(D,F)$-framed sheaves on a normal irreducible projective surface $X$ with rational singularities, where $D$ is $\mathbb{Q}$-Cartier big and nef divisor containing the singularities of $X$ and $F$ a locally free sheaf on $D$ (Theorem \ref{thm:boundedness-normal}). In Section \ref{FramedsheavesToricorbifolds}, after giving a brief introduction to the theories of root and toric stacks, we prove Theorem \ref{thm:tre}.

Finally, three appendixes are devoted to prove some results about coherent sheaves on (smooth) projective stacks: a semicontinuity theorem for Hom groups of framed sheaves, a Serre duality theorem, and a characterization of the dual of a coherent sheaf on a smooth projective stack. 

A last appendix, due to Mattia Pedrini, is devoted to the study of framed sheaves on stacky Hirzebruch surfaces. 

\subsection*{Conventions}
 
Our standard reference for the theory of stacks is \cite{book:laumonmoretbailly2000}. We denote by $k$ an algebraically closed field of characteristic zero. All schemes are defined over $k$ and are Noetherian, unless otherwise stated. A variety is a reduced separated scheme of finite type over $k$. 

Let $S$ be a generic base scheme of finite type over $k$. By \emph{Deligne-Mumford $S$-stack} we mean a separated Noetherian Deligne-Mumford stack $\mathscr{X}$ of finite type over $S$. We  denote by $p\colon \mathscr{X}\to S$ the {\em structure morphism} of $\mathscr{X}$. When $S=\mathrm{Spec}(k)$, we  omit the letter $S$. An \emph{orbifold} is a smooth Deligne-Mumford stack with generically trivial stabilizer.

The \emph{inertia stack} $\mathcal{I}(\mathscr{X})$ of a Deligne-Mumford $S$-stack $\mathscr X$ is by definition the fibred product $\mathscr{X}\times_{\mathscr{X}\times\mathscr{X}}\mathscr{X}$ with respect to the diagonal morphisms $\Delta\colon\mathscr{X}\to \mathscr{X}\times \mathscr{X}$. For a scheme $T$, an object in $\mathcal{I}(\mathscr{X})(T)$ consists of pairs $(x,g)$ where $x$ is an object of $\mathscr{X}(T)$ and $g\colon x \xrightarrow{\sim} x$ is an automorphism. A morphism $(x,g)\to (x',g')$ is a morphism $f\colon x\to x'$ in $\mathscr{X}(T)$ such that $f\circ g=g'\circ f$. Let $\sigma\colon \mathcal{I}(\mathscr{X})\to \mathscr{X}$ be the {\em forgetful morphism} which for any scheme $T$ sends a pair $(x,g)$ to $x$.

Let $\mathscr{X}$ be a Deligne-Mumford $S$-stack. An \emph{\'etale presentation} of $\mathscr{X}$ is a pair $(U, u)$, where $U$ is a $S$-scheme and $u\colon U \to \mathscr{X}$ is a \emph{representable \'etale surjective} morphism (cf.\ \cite{book:laumonmoretbailly2000}, Def.~4.1). A morphism between two \'etale presentations $(U,u)$ and $(V,v)$ of $\mathscr{X}$ is a pair $(\varphi, \alpha)$, where $\varphi\colon U\to V$ is a $S$-morphism and $\alpha\colon u\xrightarrow{\sim} v\circ\varphi$ is a 2-isomorphism. We call \emph{\'etale groupoid} associated with the \'etale presentation $u\colon U\to \mathscr{X}$ the \'etale groupoid
\begin{equation}
  \begin{tikzpicture}[xscale=2.8,yscale=-.7, ->, bend angle=25]
\node (A0_1) at (0,1) {$V:=U\times_{\mathscr{X}} U$};
\node (A1_1) at (1.15,1) {$U$}; 
\node (B1) at (0.5,0.8) {$ $};
\node (B2) at (1,0.8) {$ $};
\node (C1) at (0.5,1.2) {$ $};
\node (C2) at (1,1.2) {$ $};
\node (Comma) at (1.4,1.25) {.};
    \path (B1) edge [->]node [auto] {$\scriptstyle{}$} (B2);
        \path (C1) edge [->]node [below] {$\scriptstyle{}$} (C2);
  \end{tikzpicture}
\end{equation}
If $\mathbf{P}$ is a property of schemes which is local in the \'etale topology (for example regular, normal, reduced, Cohen-Macaulay, etc), $\mathscr{X}$ has the property $\mathbf{P}$ if for one (and hence every) \'etale presentation $u\colon U\to \mathscr{X}$, the scheme $U$ has the property $\mathbf{P}.$

A (quasi)-coherent sheaf $\mathcal{E}$ on $\mathscr{X}$ is a collection of pairs $(\mathcal{E}_{U,u}, \theta_{\varphi,\alpha})$, where for any \'etale presentation $u\colon U\to \mathscr{X}$, $\mathcal{E}_{U,u}$ is a (quasi)-coherent sheaf on $U$, and for any morphism $(\varphi,\alpha)\colon (U,u)\to (V,v)$ between two \'etale presentations of $\mathscr{X}$, $\theta_{\varphi,\alpha}\colon \mathcal{E}_U \xrightarrow{\sim} \varphi^\ast \mathcal{E}_V$ is an isomorphism which satisfies a \emph{cocycle condition} with respect to three \'etale presentations (cf.\ \cite{book:laumonmoretbailly2000}, Lemma 12.2.1; \cite{art:vistoli1989}, Def.~7.18). A vector bundle on $\mathscr{X}$ is a coherent sheaf $\mathcal{E}$ such that all $\mathcal{E}_U$ are locally free.

If $(\mathscr{X}, p)$ is a Deligne-Mumford $S$-stack, by \cite{art:keelmori1997}, Cor.~1.3-(1), there exist a separated algebraic space $X$ and a morphism $\pi\colon \mathscr{X}\to X$ such that
\begin{itemize}\setlength{\itemsep}{2pt}
\item $\pi\colon \mathscr{X}\to X$ is proper and quasi-finite;
\item if $F$ is an algebraically closed field,   $\mathscr{X}(\mathrm{Spec}(F))/\mathrm{Isom}\to X(\mathrm{Spec}(F))$ is a bijection;
\item whenever $Y\to S$ is an algebraic space and $\mathscr X\to Y$ is a morphism, the morphism factors uniquely as $\mathscr X\to X \to Y$; more generally
\item whenever $S'\to S$ is a flat morphism of schemes, and whenever $Y\to S'$ is an algebraic space and $\mathscr X\times_S S'\to Y$ is a morphism, the morphism factors uniquely as $\mathscr X\times_S S'\to X\times_S S'\to Y$; in particular
\item the natural morphism $\mathcal{O}_X\to \pi_\ast\mathcal{O}_{\mathscr{X}}$ is an isomorphism.
\end{itemize}
We call the pair $(X,\pi)$ a {\em coarse moduli space} of $\mathscr X$. If the coarse moduli space of $\mathscr{X}$ is   a scheme $X$, we call it a \emph{coarse moduli scheme}. We recall some properties of Deligne-Mumford $S$-stacks that we shall use in this paper:
\begin{itemize}\setlength{\itemsep}{2pt}
\item the functor $\pi_\ast\colon \mathrm{QCoh}(\mathscr{X})\to \mathrm{QCoh}(X)$ is exact and maps coherent sheaves to coherent sheaves (cf.\ \cite{art:abramovichvistoli2002}, Lemma 2.3.4);
\item $H^\bullet(\mathscr{X}, \mathcal{E})\simeq H^\bullet (X, \pi_\ast\mathcal{E})$
 for any quasi-coherent sheaf $\mathcal{E}$ on $\mathscr{X}$ (cf.\ \cite{art:nironi2008}, Lemma 1.10);
\item $\pi_\ast\mathcal{E}$ is an $S$-flat coherent sheaf on $X$ whenever
$\mathcal{E}$ is an $S$-flat coherent sheaf on $\mathscr{X}$ (cf.\ \cite{art:nironi2008}, Cor.~1.3-(3)).
\end{itemize} 

The projectivity of a scheme morphism is understood in the sense of Grothendieck, i.e., $f\colon X\to Y$  is  projective  if there exists a coherent sheaf $E$ on $Y$ such that $f$ factorizes as a closed immersion of $X$ into $\mathbb{P}(E)$ followed by the structural morphism $\mathbb{P}(E)\to Y.$

We   use the letters $\mathcal{E}$,  $\mathcal{G}$, $ \mathcal{F}$, ...,  for sheaves on a Deligne-Mumford $S$-stack, and the letters $E$, $F$, $G$, ...,  for sheaves on a scheme.  For any coherent sheaf $\mathcal{F}$ on a Deligne-Mumford $S$-stack $\mathscr{X}$ we  denote by $\mathcal{F}^\vee$ its {\em dual} $\mathcal{H}om(\mathcal{F},\mathcal{O}_{\mathscr{X}})$. We denote in the same way the dual of a coherent sheaf on a scheme. The projection morphism  $T\times Y \to Y$  is written as  $p_Y$ or $p_{T\times Y,Y}$.

\subsection*{Acknowledgements.} We thank Dimitri Markushevich, Philippe Eyssidieux and Richard J. Szabo for useful suggestions and interesting discussions. We thank Niels Borne for explaining us his paper \cite{art:borne2007}, and Fabio Perroni for reading and commenting on a draft of this paper.

\bigskip\section{Projective stacks}\label{Projective stacks}
 
In this section we introduce projective stacks and collect some elements of the theory of coherent sheaves on them. Our main references are \cite{art:kresch2009,art:nironi2008}. To define projective stacks one needs the notion of {\em tameness} (cf.\ \cite{art:nironi2008}, Def.~1.1), but as in characteristic zero separatedness implies tameness (cf.\ \cite{art:abramovicholssonvistoli2008}) and our Deligne-Mumford stacks are separated, we do not need to introduce that notion.

\subsection{Preliminaries on projective stacks}

The projectivity of a scheme is related to the existence of a \emph{very ample} line bundle on it. In the stacky case, one can give an equivalent notion of projectivity only for a particular class of stacks. It was proven in \cite{art:olssonstarr2003} that, under certain hypotheses, there exist locally free sheaves, called \emph{generating sheaves}, which behave like ``very ample line bundles''. In \cite{art:edidinhassettkreschvistoli2001}, another class of locally free sheaves which resemble (very) ample line bundles were introduced. It was proved in \cite{art:olssonstarr2003} that these two classes of locally free sheaves coincide. We shall use one or the other definition according to convenience.

Let $\mathscr{X}$ be a Deligne-Mumford $S$-stack with coarse moduli space $\pi\colon \mathscr{X}\to X.$
\begin{definition}
Let $\mathcal{G}$ be a locally free sheaf on $\mathscr{X}.$ We define
\begin{align}
F_{\mathcal{G}}\colon \mathrm{QCoh}(\mathscr{X})&\to \mathrm{QCoh}(X), \; \mathcal{E}\longmapsto \pi_\ast(\mathcal{E}\otimes \mathcal{G}^\vee)\ ;\\
G_{\mathcal{G}}\colon \mathrm{QCoh}(X) &\to  \mathrm{QCoh}(\mathscr{X}), \; E\mapsto \pi^\ast E\otimes \mathcal{G}\ .
\end{align}
\end{definition}
\begin{remark}\label{rem:exactness}
The functor $F_{\mathcal{G}}$ is exact since $\mathcal{G}^\vee$ is locally free and the direct image functor $\pi_\ast$ is exact. The functor $G_{\mathcal{G}}$ is exact when the morphism $\pi$ is flat. This happens for instance if the stack is a gerbe over a scheme i.e., a stack over a scheme $Y$ which \'etale locally admits a section and such that any two local sections are locally 2-isomorphic, or in the case of root stacks over schemes.
\end{remark}
\begin{definition}
A locally free sheaf $\mathcal{G}$ is said to be a \emph{generator} for the quasi-coherent sheaf $\mathcal{E}$ if the adjunction morphism (left adjoint to the identity $\mathrm{id}\colon \pi_\ast(\mathcal{E}\otimes \mathcal{G}^\vee)\to \pi_\ast(\mathcal{E}\otimes \mathcal{G}^\vee)$)
\begin{equation}\label{eq:generating}
\theta_{\mathcal{G}}(\mathcal{E})\colon \pi^\ast\pi_\ast(\mathcal{E}\otimes \mathcal{G}^\vee)\otimes \mathcal{G}\to \mathcal{E}
\end{equation}
is surjective. It is a \emph{generating sheaf}  for $\mathscr{X}$ if it is a generator for every quasi-coherent sheaf on $\mathscr{X}.$
\end{definition}
A generating sheaf can be considered as a very ample sheaf relatively to the morphism $\pi\colon\mathscr{X}\to X.$ Indeed, the property expressed by \eqref{eq:generating} resembles a similar property for very ample line bundles (\cite{art:Grothendieck1961-III}, Thm.~2.1.1 Chap.~III): if $f\colon Y\to Z$ is a proper morphism, $\mathcal{O}_Y(1)$ is a very ample line bundle on $Y$ relative to $f$,  and $E$  is coherent sheaf on $Y$, there is a positive integer $N$ such that the adjunction morphism $f^\ast f_\ast\mathcal{H}om(\mathcal{O}_Y(-n),E)\otimes\mathcal{O}_Y(-n)\to E$ is surjective for any integer $n\geq N.$

Let $E$ be a quasi-coherent sheaf on $X.$ Since $\mathcal{G}$ is  locally free,
\begin{equation}
\operatorname{Hom}(\pi^\ast E\otimes \mathcal{G}, \pi^\ast E\otimes \mathcal{G})\simeq \operatorname{Hom}(\pi^\ast E,\mathcal{H}om(\mathcal{G},\pi^\ast E\otimes \mathcal{G}))\ .
\end{equation}
Define the morphism $\varphi_{\mathcal{G}}(E)$ as the right adjoint to the identity $\mathrm{id}\colon \pi^\ast E\otimes \mathcal{G}\to \pi^\ast E\otimes \mathcal{G}$:
\begin{equation}
\varphi_{\mathcal{G}}(E)\colon E\to\pi_\ast\left(\mathcal{H}om(\mathcal{G},\pi^\ast E\otimes \mathcal{G})\right)=F_{\mathcal{G}}(G_{\mathcal{G}}(E))\ .
\end{equation}

\begin{lemma}[{\cite[Cor.~5.4]{art:olssonstarr2003}}]\label{lem:projectionformula}
Let $\mathcal{F}$ be a quasi-coherent sheaf on $\mathscr{X}$ and $E$ a quasi-coherent sheaf on $X.$ A projection formula holds:
\begin{equation}
\pi_\ast(\pi^\ast(E)\otimes\mathcal{F}) \simeq E\otimes \pi_\ast\mathcal{F}\ .
\end{equation}
Moreover, this is functorial in the sense that if $f\colon \mathcal{F}\to \mathcal{F}'$ is a morphism of quasi-coherent sheaves on $\mathscr{X}$ and $g\colon E\to E'$ is a morphism of quasi-coherent sheaves on $X$, one has
\begin{equation}
\pi_\ast(\pi^\ast(g)\otimes f)=g\otimes \pi_\ast f \ .
\end{equation}
\end{lemma}
\proof The projection formula is proved at the beginning of the proof of Corollary 5.4 in \cite{art:olssonstarr2003}.
\endproof
According to this Lemma, $\varphi_{\mathcal{G}}(E)$ can be rewritten as
\begin{equation}
\varphi_{\mathcal{G}}(E)\colon E\to E\otimes \pi_\ast\left(\mathcal{E}nd(\mathcal{G})\right)\ ,
\end{equation}
and is the morphism given by tensoring a section by the identity endomorphism; in particular it is injective.
\begin{lemma}[{\cite[Lemma 2.9]{art:nironi2008}}]
Let $\mathcal{F}$ be a quasi-coherent sheaf on $\mathscr{X}$ and $E$ a coherent sheaf on $X$. The compositions
\begin{gather}
  \begin{tikzpicture}[xscale=4,yscale=-1.2]
    \node (A0_0) at (0, 0) {$F_{\mathcal{G}}(\mathcal{F})$};
    \node (A0_1) at (1, 0) {$F_{\mathcal{G}}\circ G_{\mathcal{G}}\circ F_{\mathcal{G}}(\mathcal{F})$};
    \node (A0_2) at (2, 0) {$F_{\mathcal{G}}(\mathcal{F})$};
    \path (A0_0) edge [->]node [auto] {$\scriptstyle{\varphi_{\mathcal{G}}(F_{\mathcal{G}}(\mathcal{F}))}$} (A0_1);
    \path (A0_1) edge [->]node [auto] {$\scriptstyle{F_{\mathcal{G}}(\theta_{\mathcal{G}}(\mathcal{F}))}$} (A0_2);
  \end{tikzpicture}
  \\
  \begin{tikzpicture}[xscale=4,yscale=-1.2]
    \node (A0_0) at (0, 0) {$G_{\mathcal{G}}(E)$};
    \node (A0_1) at (1, 0) {$G_{\mathcal{G}}\circ F_{\mathcal{G}}\circ G_{\mathcal{G}}(E)$};
    \node (A0_2) at (2, 0) {$G_{\mathcal{G}}(E)$\ .};
    \path (A0_0) edge [->]node [auto] {$\scriptstyle{G_{\mathcal{G}}(\varphi_{\mathcal{G}}(E))}$} (A0_1);
    \path (A0_1) edge [->]node [auto] {$\scriptstyle{\theta_{\mathcal{G}}(G_{\mathcal{G}}(E))}$} (A0_2);
  \end{tikzpicture}
\end{gather} are the identity endomorphisms.
\end{lemma}

Following \cite{art:edidinhassettkreschvistoli2001} we introduce another definition of ``ampleness'' for sheaves on stacks.
\begin{definition} \label{def:generating-sheaf}
A locally free sheaf $\mathcal{V}$ on $\mathscr{X}$ is $\pi$\emph{-ample} if for every geometric point of $\mathscr{X}$ the natural representation of the stabilizer group at that point on the fibre of $\mathcal{V}$ is faithful. A locally free sheaf $\mathcal{G}$ on $\mathscr{X}$ is $\pi$\emph{-very ample} if for every geometric point of $\mathscr{X}$ the natural representation of the stabilizer group at that point on the fibre of $\mathcal{G}$ contains every irreducible representation.
\end{definition}
The relation between these two notions is explained in \cite{art:kresch2009}, Sect.~5.2. In particular, we have the following result.
\begin{proposition}\label{prop:veryampleness}
Let $\mathcal{V}$ be a $\pi$-ample sheaf on $\mathscr{X}$ and $N$ the maximum between the numbers of conjugacy classes of any geometric stabilizer group of $\mathscr{X}$. Then, for any $r\ge N$, the locally free sheaf $\bigoplus_{i=1}^r\mathcal{V}^{\otimes i}$ is $\pi$-very ample. 
\end{proposition}
As shown in \cite{art:olssonstarr2003}, Thm.~5.2, a locally free sheaf $\mathcal{V}$ on $\mathscr{X}$ is $\pi$-very ample if and only if it is a   generating sheaf.

\begin{remark}\label{rem:generatingsheaf}
Let $\varphi\colon\mathscr{Y}\to\mathscr{X}$ be a representable morphism of Deligne-Mumford $S$-stacks. By the universal property of the coarse moduli spaces, $\varphi$ induces a morphism $\bar{\varphi}\colon Y\to X$ between the corresponding coarse moduli spaces together with a commutative diagram
\begin{equation}
  \begin{tikzpicture}[xscale=1,yscale=-1]
    \node (A0_0) at (0, 0) {$\mathscr{Y}$};
    \node (A0_2) at (2, 0) {$\mathscr{X}$};
    \node (A1_3) at (3, 1) {$.$};
    \node (A2_0) at (0, 2) {$Y$};
    \node (A2_2) at (2, 2) {$X$};
    \path (A0_0) edge [->]node [left] {$\scriptstyle{\pi_{\mathscr{Y}}}$} (A2_0);
    \path (A0_0) edge [->]node [auto] {$\scriptstyle{\varphi}$} (A0_2);
    \path (A0_2) edge [->]node [auto] {$\scriptstyle{\pi_{\mathscr{X}}}$} (A2_2);
    \path (A2_0) edge [->]node [auto] {$\scriptstyle{\bar{\varphi}}$} (A2_2);    
  \end{tikzpicture} 
\end{equation}
By \cite{book:laumonmoretbailly2000}, Prop.~2.4.1.3, for any geometric point of $\mathscr{Y}$ the morphism $\varphi$ induces an injective map between the stabilizer groups at that point and at the corresponding image point. So if $\mathcal{V}$ is a $\pi_{\mathscr{X}}$-ample sheaf on $\mathscr{X}$, then $\varphi^\ast \mathcal{V}$ is a $\pi_{\mathscr{Y}}$-ample sheaf on $\mathscr{Y}$. Denote by $N_{\mathscr{X}}$ (resp.\ $N_{\mathscr{Y}}$) the maximum of the numbers of conjugacy classes of any geometric stabilizer group of $\mathscr{X}$ (resp.\ $\mathscr{Y}$). If  $N_{\mathscr{X}}\geq N_{\mathscr{Y}}$ by Proposition \ref{prop:veryampleness} we get that $\oplus_{i=1}^r \varphi^\ast \mathcal{V}^{\otimes i}$ is $\pi_{\mathscr{Y}}$-very ample for any $r\geq N_{\mathscr{X}}$.
\end{remark}

\begin{definition}[{\cite[Def.~2.9]{art:edidinhassettkreschvistoli2001}}]
Let $\mathscr{X}$ be a stack of finite type over a base scheme $S.$ We say $\mathscr{X}$ is a {\em global $S$-quotient} if it is isomorphic to a stack of the form $[T/G]$, where $T$ is an algebraic space of finite type over $S$ and $G$ is an $S$-flat group scheme which is a group subscheme (a locally closed subscheme which is a subgroup) of the general linear group scheme $\mathrm{GL}_{N, S}$ over $S$  for some integer $N.$
\end{definition}
\begin{theorem}[{\cite[Sect.~5]{art:olssonstarr2003}}]\label{thm:generatingsheaf}
$ $
\begin{itemize}\setlength{\itemsep}{2pt}
\item[(i)] A Deligne-Mumford $S$-stack $\mathscr{X}$ which is a global $S$-quotient always has a generating sheaf $\mathcal{G}.$ 
\item[(ii)] Under the same hypothesis of (i), let $\pi\colon \mathscr{X}\to X$ be the coarse moduli space of $\mathscr{X}$ and $f\colon X'\to X$ a morphism of algebraic spaces. Then $p_{\mathscr{X}\times_X X', \mathscr{X}}^\ast \mathcal{G}$ is a generating sheaf for $\mathscr{X}\times_X X'$.
\end{itemize}
\end{theorem}

Now we are ready to give the definition of projective stack.
\begin{defin}[{\cite[Def.~5.5]{art:kresch2009}}]
A Deligne-Mumford stack $\mathscr{X}$ is a \emph{(quasi-)projective stack} if $\mathscr{X}$ admits a (locally) closed embedding into a smooth proper Deligne-Mumford stack which has a projective coarse moduli scheme.\parbox{2mm}{\hfill} \hfill\hfill$\oslash$
\end{defin}
\begin{proposition}[{\cite[Thm.~5.3]{art:kresch2009}}]\label{prop:kresch}
Let $\mathscr{X}$ be a Deligne-Mumford stack. The following statements are equivalent:
\begin{itemize}\setlength{\itemsep}{2pt}
\item[(i)]   $\mathscr{X}$ is (quasi-)projective.
\item[(ii)]   $\mathscr{X}$ has a (quasi-)projective coarse moduli scheme and has a generating sheaf.
\item[(iii)]  $\mathscr{X}$ is a separated global quotient with a coarse moduli space which is a (quasi-)projective scheme.
\end{itemize}
\end{proposition}
\begin{definition}
Let $\mathscr{X}$ be a projective stack with coarse moduli scheme $X.$  A \emph{polarization} for $\mathscr{X}$ is a pair $(\mathcal{G}, \mathcal{O}_X(1))$, where $\mathcal{G}$ is a generating sheaf of $\mathcal{X}$ and $\mathcal{O}_X(1)$ is an ample line bundle on $X.$
\end{definition}

We  give a relative version of the notion of projective stacks. 
\begin{definition}
Let $p\colon \mathscr{X}\to S$ be a Deligne-Mumford $S$-stack which is a global $S$-quotient with a coarse moduli scheme $X$ such that $p$ factorizes as $\pi\colon \mathscr{X}\to X$ followed by a projective morphism $\rho\colon X\to S.$ We  call $p\colon \mathscr{X}\to S$ a \emph{family of projective stacks}.
\end{definition}
\begin{remark}
Let $p\colon \mathscr{X}=[T/G]\to S$ be a family of projective stacks. For any geometric point $s\in S$ we have the following cartesian diagram
\begin{equation}
  \begin{tikzpicture}[xscale=1.5,yscale=-1.2]
    \node (A0_0) at (0, 0) {$\mathscr{X}_s$};
    \node (A0_2) at (2, 0) {$\mathscr{X}$};
    \node (A1_1) at (1, 1) {$\square$};
    \node (A2_0) at (0, 2) {$X_s$};
    \node (A2_2) at (2, 2) {$X$};
    \node (A1_3) at (1, 3) {$\square$};
    \node (A4_0) at (0, 4) {$s$};
    \node (A4_2) at (2, 4) {$S$};
    \path (A0_0) edge [->]node [auto] {$\scriptstyle{\pi_s}$} (A2_0);
    \path (A0_0) edge [->]node [auto] {$\scriptstyle{}$} (A0_2);
    \path (A0_2) edge [->]node [auto] {$\scriptstyle{\pi}$} (A2_2);
    \path (A2_0) edge [->]node [auto] {$\scriptstyle{}$} (A2_2);
    \path (A4_0) edge [->]node [auto] {$\scriptstyle{}$} (A4_2);
    \path (A2_2) edge [->]node [auto] {$\scriptstyle{\rho}$} (A4_2);
    \path (A2_0) edge [->]node [auto] {$\scriptstyle{\rho_s}$} (A4_0);        
  \end{tikzpicture} 
\end{equation}
with $\mathscr{X}_s=[T_s/G_s]$, where $T_s$ and $G_s$ are the fibres of $T$ and $G$, respectively. Since the morphism $\rho$ is projective, the fibres $X_s$ are projective schemes. The property of being coarse moduli spaces is invariant under base change, so that each $X_s$ is a coarse moduli scheme for $\mathscr{X}_s$, and each $\mathscr{X}_s$ is a projective stack. 
\end{remark}
By Theorem \ref{thm:generatingsheaf}, a family of projective stacks $p\colon \mathscr{X}\to S$ has a generating sheaf $\mathcal{G}$ and the fibre of $\mathcal{G}$ at a geometric point $s\in S$ is a generating sheaf for $\mathscr{X}_s.$ This justifies the following definition.
\begin{definition}
Let $p\colon \mathscr{X}\to S$ be a family of projective stacks. A \emph{relative} polarization of $p\colon \mathscr{X}\to S$ is a pair $(\mathcal{G}, \mathcal{O}_X(1))$ where $\mathcal{G}$ is a generating sheaf for $\mathscr{X}$ and $\mathcal{O}_X(1)$ is an ample line bundle relative to $\rho\colon X\to S.$
\end{definition}

\subsection{Coherent sheaves on projective stacks}

In this section we briefly recall the theory of coherent sheaves on projective stacks from \cite{art:nironi2008}, Sect.~3.1. In particular, we shall see that the functor $F_{\mathcal{G}}$ preserves the dimension and the pureness of coherent sheaves on projective stacks. Let us fix a projective stack $\mathscr{X}$ of dimension $d$, with  a coarse moduli scheme $\pi\colon \mathscr{X}\to X$, and a polarization $(\mathcal{G}, \mathcal{O}_X(1))$ on it.

\begin{remark}\label{rem:equiv}
By \cite{art:kresch2009}, Prop.~5.1, the stack $\mathscr{X}$ is of the form $[T/G]$ with $T$ a quasi-projective scheme and $G$ a linear algebraic group acting on $T$. This implies that the category of coherent sheaves on $\mathscr{X}$ is equivalent to the category of coherent $G$-equivariant sheaves on $T$ (cf.\ \cite{book:laumonmoretbailly2000}, Example 12.4.6 and \cite{art:vistoli1989}, Example 7.21). In the following, we shall use this correspondence freely. 
\end{remark}

\begin{definition}
Let $\mathcal{E}$ be a coherent sheaf on $\mathscr{X}.$ The \emph{support} $\mathrm{supp}(\mathcal{E})$ of $\mathcal{E}$ is the closed substack associated with the ideal $\mathcal{I} = \operatorname{ker}(\mathcal{O}_{\mathscr{X}} \to \mathcal{E}nd(\mathcal{E}))$.
The dimension $\dim(\mathcal{E})$ of $\mathcal{E}$ is the dimension of its support.
 We say that $\mathcal{E}$ is a \emph{pure sheaf of dimension} $\dim(\mathcal{E})$ if for any nonzero subsheaf $\mathcal{G}$ of $\mathcal{E}$ the support of $\mathcal{G}$ is pure of dimension $\dim(\mathcal{E}).$ We say that $\mathcal{E}$ is \emph{torsion-free} if it is a pure sheaf of dimension $d.$
\end{definition}
\begin{remark}\label{rem:pureness}
Let $u\colon U\to \mathscr{X}$ be an \'etale presentation of $\mathscr{X}.$ Let $\mathcal{E}$ be a coherent sheaf on $\mathscr{X}$ of dimension $d.$ First note that $u^\ast \mathcal{E}$ is exactly the representative $\mathcal{E}_{U,u}$ of $\mathcal{E}$ on $U$. As explained in \cite{art:nironi2008}, Rem.~3.3, $\mathrm{supp}(u^\ast \mathcal{E}) \to \mathrm{supp}(\mathcal{E})$ is an \'etale presentation of $\mathrm{supp}(\mathcal{E}).$ Moreover, $\dim(\mathcal{E})=\dim(u^\ast \mathcal{E})$ and $\mathcal{E}$ is pure if and only if $u^\ast\mathcal{E}$ is pure. 
\end{remark}

As it was shown in \cite{art:nironi2008}, Sect.~3 (cf.\ also \cite{book:huybrechtslehn2010}, Def.~1.1.4), there exists a unique filtration, the so-called \emph{torsion filtration}, of a coherent sheaf $\mathcal{E}$
\begin{equation}
0\subseteq T_0(\mathcal{E})\subseteq T_1(\mathcal{E})\subseteq \cdots \subseteq T_{\dim(\mathcal{E})-1}(\mathcal{E})\subseteq T_{\dim(\mathcal{E})}(\mathcal{E})=\mathcal{E}\ ,
\end{equation}
where $T_i(\mathcal{E})$ is the maximal subsheaf of $\mathcal{E}$ of dimension $\leq i.$ Note that $T_i(\mathcal{E})/T_{i-1}(\mathcal{E})$ is zero or pure of dimension $i.$ In particular, $\mathcal{E}$ is pure if and only if $T_{\dim(\mathcal{E})-1}(\mathcal{E})=0.$
\begin{definition}
The \emph{saturation} of a subsheaf $\mathcal{E}'\subset \mathcal{E}$ is the minimal subsheaf $\bar{\mathcal{E}}'$ containing $\mathcal{E}'$ such that $\mathcal{E}/\bar{\mathcal{E}}'$ is zero or pure of dimension $\dim(\mathcal{E}).$
\end{definition}
Clearly, the saturation of $\mathcal{E}'$ is the kernel of the surjection
\begin{equation}
\mathcal{E}\to \mathcal{E}/\mathcal{E}'\to \frac{\mathcal{E}/\mathcal{E}'}{T_{\dim(\mathcal{E})-1}(\mathcal{E}/\mathcal{E}')}\ .
\end{equation}
\begin{lemma}[{\cite[Lemma 3.4]{art:nironi2008}}]\label{lem:support}
Let $\mathscr{X}$ be a projective stack with coarse moduli scheme $\pi\colon \mathscr{X}\to X$. Let $\mathcal{E}$ be a coherent sheaf on $\mathscr{X}.$ Then we have
\begin{itemize}\setlength{\itemsep}{2pt}
\item[(i)] $\pi(\mathrm{Supp}(\mathcal{E}))=\pi(\mathrm{Supp}(\mathcal{E}\otimes \mathcal{G}^\vee))\supseteq \mathrm{Supp}(F_{\mathcal{G}}(\mathcal{E}));$
\item[(ii)] $F_{\mathcal{G}}(\mathcal{E})$ is zero if and only if $\mathcal{E}$ is zero.
\end{itemize}
\end{lemma}
\begin{proposition}\label{prop:pureness}
Let $\mathscr{X}$ be a projective stack with coarse moduli scheme $\pi\colon \mathscr{X}\to X$. A coherent sheaf $\mathcal{E}$ on $\mathscr{X}$ and the sheaf $F_{\mathcal{G}}(\mathcal{E})$ on $X$ have the same dimension. Moreover, $\mathcal{E}$ is pure if and only if $F_{\mathcal{G}}(\mathcal{E})$ is pure.
\end{proposition}
\proof
Assume first that $\mathcal{E}$ is pure. Then the necessary part is proved in \cite{art:nironi2008}, Prop.~3.6. For the sufficient part, let us consider the short exact sequence
\begin{equation}\label{eq:torsionpart}
0\to T_{\dim(\mathcal{E})-1}(\mathcal{E})\to \mathcal{E}\to \mathcal{Q}\to 0\ .
\end{equation}
Since the functor $F_{\mathcal{G}}$ is exact, we obtain
\begin{equation}
0\to F_{\mathcal{G}}(T_{\dim(\mathcal{E})-1}(\mathcal{E}))\to F_{\mathcal{G}}(\mathcal{E})\to F_{\mathcal{G}}(\mathcal{Q})\to 0\ .
\end{equation}
By Lemma \ref{lem:support}, $\mathrm{Supp}(F_{\mathcal{G}}(T_{\dim(\mathcal{E})-1}(\mathcal{E})))\subseteq \pi(\mathrm{Supp}(T_{\dim(\mathcal{E})-1}(\mathcal{E})))$, and since $\pi$ preserves the dimensions, $\dim F_{\mathcal{G}}(T_{\dim(\mathcal{E})-1}(\mathcal{E}))\leq \dim\mathcal{E}-1$. As by hypothesis $F_{\mathcal{G}}(\mathcal{E})$ is pure of dimension $\dim\mathcal{E}$, we have $F_{\mathcal{G}}(T_{\dim(\mathcal{E})-1}(\mathcal{E}))=0$ and therefore $T_{\dim(\mathcal{E})-1}(\mathcal{E})=0$ by Lemma \ref{lem:support}.

If $\mathcal{E}$ is not pure, to prove the assertion it is enough to use the short exact sequence \eqref{eq:torsionpart} and a similar argument as before applied to $\mathcal{E}$ and $\mathcal{Q}$.
\qedhere
%Consider the short exact sequence
%\begin{equation}
%0\to T_{\dim(\mathcal{E})-1}(\mathcal{E})\to \mathcal{E}\to \mathcal{Q}\to 0\ .
%\end{equation}
%Since the functor $F_{\mathcal{G}}$ is exact, we obtain
%\begin{equation}
%0\to F_{\mathcal{G}}(T_{\dim(\mathcal{E})-1}(\mathcal{E}))\to F_{\mathcal{G}}(\mathcal{E})\to F_{\mathcal{G}}(\mathcal{Q})\to 0\ .
%\end{equation}
%By Proposition \ref{prop:pureness}, $\dim F_{\mathcal{G}}(\mathcal{Q})=\dim \mathcal{Q}=\dim \mathcal{E}.$ On the other hand, by Lemma \ref{lem:support} we get $\mathrm{Supp}(F_{\mathcal{G}}(\mathcal{Q}))\subseteq \mathrm{Supp}(F_{\mathcal{G}}(\mathcal{E}))\subseteq \pi(\mathrm{Supp}(\mathcal{E})).$ Therefore
%\begin{equation}
%\dim \mathcal{E}=\dim F_{\mathcal{G}}(\mathcal{Q})\leq \dim F_{\mathcal{G}}(\mathcal{E}) \leq \dim \pi(\mathrm{Supp}(\mathcal{E}))=\dim \mathcal{E} \ .
%\end{equation}
%Thus, $\dim F_{\mathcal{G}}(\mathcal{E})=\dim\mathcal{E}.$
\endproof
For pure coherent sheaves on $\mathscr{X}$, the functor $F_{\mathcal{G}}$ \emph{preserves} the supports.
\begin{corollary}[{\cite[Cor.~3.8]{art:nironi2008}}]\label{cor:supportpure}
Let $\mathcal{E}$ be a pure coherent sheaf on $\mathscr{X}.$ Then $\mathrm{Supp}(F_{\mathcal{G}}(\mathcal{E}))=\pi(\mathrm{Supp}(\mathcal{E})).$
\end{corollary}
Further, the functor $F_{\mathcal{G}}$ is compatible with torsion filtrations.
\begin{corollary}[{\cite[Cor.~3.7]{art:nironi2008}}]\label{cor:torsionfiltration}
The functor $F_{\mathcal{G}}$ sends the torsion filtration $0\subseteq T_0(\mathcal{E})\subseteq \cdots \subseteq T_{\dim(\mathcal{E})}(\mathcal{E})$ $=\mathcal{E}$ of $\mathcal{E}$ to the torsion filtration  of $F_{\mathcal{G}}(\mathcal{E})$, that is, $F_{\mathcal{G}}(T_i(\mathcal{E}))=T_i(F_{\mathcal{G}}(\mathcal{E}))$ for $i=0, \ldots, \dim(\mathcal{E}).$
\end{corollary}

\begin{exa}\label{ex:curvesurface}
Let $\mathscr{X}$ be a smooth projective stack and $\pi\colon \mathscr{X}\to X$ its coarse moduli scheme. By \cite{art:nironi2008}, Lemma 6.9, any torsion-free sheaf $\mathcal{E}$ on $\mathscr{X}$ fits into an exact sequence
\begin{equation}
0\to \mathcal{E}\to \mathcal{E}^{\vee\vee} \to \mathcal{Q}\to 0\ .
\end{equation}
Let $u\colon U \to \mathscr{X}$ be an \'etale presentation of $\mathscr{X}.$ In particular, $U$ is a regular scheme of dimension $\dim(\mathscr{X})$ and $u$ is a flat morphism. By applying the functor $u^\ast$, we obtain an exact sequence
\begin{equation}
0\to u^\ast \mathcal{E}\to u^\ast \mathcal{E}^{\vee\vee} \to u^\ast\mathcal{Q}\to 0\ .
\end{equation}
Note that $u^\ast \mathcal{E}^{\vee\vee}\simeq (u^\ast\mathcal{E})^{\vee\vee}$ (cf.\ \cite{book:matsumura1980}, Claim 3.E p.\ 20). Moreover, $\mathrm{codim}\,\mathcal{Q}\geq 2$ and $u^\ast(\mathcal{E})^{\vee\vee}$ is locally free except on a closed subset of $U$ of codimension at least 3 (cf.\ \cite{art:hartshorne1980}, Sect.~1). If $\dim(\mathscr{X})=1$, we obtain $\mathcal{Q}=0$ and $u^\ast \mathcal{E}^{\vee\vee}$ is locally free. Thus $\mathcal{E}^{\vee\vee}$ is locally free and $\mathcal{E}\simeq \mathcal{E}^{\vee \vee}.$ Therefore any torsion-free sheaf on a smooth projective stack of dimension one is locally free. If $\dim(\mathscr{X})=2$, then $\mathcal{Q}$ is a zero-dimensional sheaf and $\mathcal{E}^{\vee\vee}$ is locally free. Thus we obtain the analogue of the usual characterization of torsion-free sheaves on smooth curves and surfaces (cf.\ \cite{book:huybrechtslehn2010}, Example 1.1.16). \hfill $\triangle$
\end{exa}

\subsection{Hilbert polynomial}

We define a polynomial which will be the analogue of the usual Hilbert polynomial for coherent sheaves on projective schemes. Let us fix a projective stack $\mathscr{X}$ of dimension $d$, with coarse moduli space $\pi\colon \mathscr{X}\to X$, and a polarization $(\mathcal{G}, \mathcal{O}_X(1))$ on it. (This was called {\em modified Hilbert polynomial} in \cite{art:nironi2008}).

\begin{definition}
The \emph{Hilbert polynomial} of a coherent sheaf $\mathcal{E}$ on $\mathscr{X}$ is
\begin{equation}
P_{\mathcal{G}}(\mathcal{E},n):=\chi(\mathscr{X},\mathcal{E}\otimes\mathcal{G}^\vee\otimes \pi^\ast\mathcal{O}_X(n))=\chi(X,F_{\mathcal{G}}(\mathcal{E})\otimes\mathcal{O}_X(n))=P(F_{\mathcal{G}}(\mathcal{E}),n)\ .
\end{equation}
\end{definition}

By Proposition \ref{prop:pureness}, $\dim F_{\mathcal{G}}(\mathcal{E})=\dim(\mathcal{E}).$ The function $n\mapsto P_{\mathcal{G}}(\mathcal{E},n)$ is a polynomial with rational coefficients by \cite{book:huybrechtslehn2010}, Lemma 1.2.1, and can be uniquely written in the form
\begin{equation}
P_{\mathcal{G}}(\mathcal{E},n)=\sum_{i=0}^{\dim(\mathcal{E})} \alpha_{\mathcal{G},i}(\mathcal{E}) \frac{n^i}{i!}\in\mathbb{Q}[n]\ .
\end{equation}
Moreover, the Hilbert polynomial is additive on short exact sequences since $F_{\mathcal{G}}$ is an exact functor (cf.\ Remark \ref{rem:exactness}) and the Euler characteristic is additive on short exact sequences.  

Let $\mathcal{E}$ be a coherent sheaf on $\mathscr{X}.$ We call \emph{multiplicity} of $\mathcal{E}$ the leading coefficient $\alpha_{\mathcal{G},\dim(\mathcal{E})}(\mathcal{E})$ of its Hilbert polynomial. The \emph{reduced Hilbert polynomial} of $\mathcal{E}$ is
\begin{equation}
p_{\mathcal{G}}(\mathcal{E},n):=\frac{P_{\mathcal{G}}(\mathcal{E},n)}{\alpha_{\mathcal{G},\dim(\mathcal{E})}(\mathcal{E})}\ .
\end{equation}
The \emph{hat-slope} of $\mathcal{E}$ is
\begin{equation}
\hat{\mu}_{\mathcal{G}}(\mathcal{E}):=\frac{\alpha_{\mathcal{G},\dim(\mathcal{E})-1}(\mathcal{E})}{\alpha_{\mathcal{G},\dim(\mathcal{E})}(\mathcal{E})}\ .
\end{equation}
For a $d$-dimensional coherent sheaf $\mathcal{E}$, its \emph{rank} is
\begin{equation}
\mathrm{rk}_{\mathcal{G}}(\mathcal{E}):=\frac{\alpha_{\mathcal{G},d}(\mathcal{E})}{\alpha_d(\mathcal{O}_X)},
\end{equation}
where $\alpha_d(\mathcal{O}_X)$ is the leading coefficient of the Hilbert polynomial of $\mathcal{O}_X.$

\begin{remark}\label{rem:coeff}
Let $\mathcal{E}$ be a coherent sheaf of dimension $d.$ Let $\mathcal{E}'$ be a $d$-dimensional coherent subsheaf of $\mathcal{E}$ and $\bar{\mathcal{E}}'$ its saturation. Then $\mathrm{rk}_{\mathcal{G}}(\bar{\mathcal{E}}')=\mathrm{rk}_{\mathcal{G}}(\mathcal{E}')$ and $\hat{\mu}_{\mathcal{G}}(\bar{\mathcal{E}}')\geq \hat{\mu}_{\mathcal{G}}(\mathcal{E}').$ 
\end{remark}

\subsubsection{Smooth case}\label{sec:smoothcase}

If $\mathscr{X}$ is smooth one can give another definition of rank of a coherent sheaf. Let $\mathcal{E}$ be a $d$-dimensional coherent sheaf. The \emph{orbifold rank} of $\mathcal{E}$ is
\begin{equation}
\mathrm{ork}(\mathcal E)=\frac{1}{\alpha_d(\mathcal{O}_X)}\int_{\mathscr X}^{et} \mathrm{ch}^{et}(\mathcal E)\,[\pi^\ast c_1^{et}(\mathcal O_X(1))]^d\ ,
\end{equation}
where $\mathrm{ch}^{et}(\mathcal E)$ is the {\em \'etale} Chern character of $\mathcal{E}$ and $\int_{\mathscr X}^{et}$ denotes the pushfoward $p_\ast\colon H_{et}^\bullet(\mathscr{X})\to H_{et}^\bullet(\mathrm{Spec}(k))\simeq  \mathbb{Q}$ of the morphism $p\colon \mathscr{X}\to \mathrm{Spec}(k)$, which is proper since $\mathscr{X}$ is projective. (For a more detailed introduction to the \'etale cohomology of a Deligne-Mumford stack, we refer to \cite{art:borne2007}, App.~C.)

The \emph{degree} of $\mathcal{E}$ is
\begin{equation}
\deg_{\mathcal{G}}(\mathcal{E}):=\alpha_{\mathcal{G},d-1}(\mathcal{E})-\mathrm{ork}(\mathcal{E})\alpha_{\mathcal{G},d-1}(\mathcal{O}_{\mathscr{X}})\ ,
\end{equation}
and its \emph{slope} is
\begin{equation}
\mu_{\mathcal{G}}(\mathcal{E}):=\frac{\deg_{\mathcal{G}}(\mathcal{E})}{\mathrm{ork}(\mathcal{E})}\ .
\end{equation}
In this case the (in)equalities in Remark \ref{rem:coeff} are still valid. 

\begin{remark} \label{rem:orbirank}
Assume moreover that $\mathscr X$ is an orbifold. Then the only codimension zero component of the {\em inertia stack} $\mathcal I(\mathscr X)$ is $\mathscr X$ (which is {\em associated} with the trivial stabilizer), so that, by the T\"oen-Riemann-Roch Theorem \cite{phd:toen1999,art:toen1999} (see also \cite{art:borne2007}, App.~C), we get
\begin{equation}
\mathrm{ork}(\mathcal{E})=\frac{\alpha_d(\mathcal{E})}{\alpha_d(\mathcal{O}_X)}\ ,
\end{equation} 
where $\alpha_d(\mathcal{E})$ is the leading coefficient of the Hilbert polynomial of $\pi_\ast(\mathcal{E}).$ More details about the T\"oen-Riemann-Roch Theorem will be given in Section \ref{sec:stability}.

Let $\mathcal{E}$ be a coherent sheaf on $\mathscr{X}$. Then $\mathrm{ork}(\mathcal{E})$ is the zero degree part $\mathrm{ch}_0^{et}(\mathcal E)$ of the \'etale Chern character of $\mathcal{E}$. This is a trivial check if $\mathcal{E}$ is locally free. In general, we can note that  by \cite{art:kresch2009}, Prop.~5.1, $\mathscr{X}$ has the \emph{resolution property}, i.e., any coherent sheaf on $\mathscr{X}$ admits a surjective morphism from a locally free sheaf. Since $\mathscr{X}$ is also smooth, the Grothendieck group of coherent sheaves on $\mathscr{X}$ is isomorphic to the Grothendieck group of locally free sheaves on $\mathscr{X}$. Therefore $\mathrm{ork}(\mathcal{E})=\mathrm{ch}_0^{et}(\mathcal E)$ for any coherent sheaf $\mathcal{E}$ on $\mathscr{X}$. As a byproduct, we get $\mathrm{rk}_{\mathcal G}(\mathcal E) =\mathrm{ork}(\mathcal G)\,\mathrm{ork}(\mathcal E)$. Moreover, we have the following relation between the hat-slope and the slope of $\mathcal{E}$, which is a generalization of the usual relation in the case of coherent sheaves on projective schemes (cf.\ \cite{book:huybrechtslehn2010}, Sect.~1.6):
\begin{equation}\label{eq:relationslope}
\mu_{\mathcal{G}}(\mathcal{E})=\mathrm{ork}(\mathcal G)\alpha_d(\mathcal{O}_X)\hat{\mu}_{\mathcal{G}}(\mathcal{E})-\alpha_{\mathcal{G},d-1}(\mathcal{O}_{\mathscr{X}})\ .
\end{equation}
\end{remark}

\subsection{Families of coherent sheaves}

We introduce the notions of set-theoretic family of coherent sheaves and bounded family (cf.\ \cite{book:KleimanSGA6}, Sect.~1.12). Let us fix a family of projective stacks $p\colon \mathscr{X}\to S$ and a relative polarization $(\mathcal{G}, \mathcal{O}_X(1))$ on it.

Given a point $s\in S$ with residue field $k(s)$ and an extension field $K$ of $k(s)$, a coherent sheaf on a fibre of $p$ is a coherent sheaf $\mathcal{E}_K$ on $\mathscr{X}_K:=\mathscr{X}\times_S \mathrm{Spec}(K)$. Given two field extensions $K$ and $K'$, two  coherent sheaves $\mathcal{E}_K$ and $\mathcal{E}'_{K'}$ on $\mathscr{X}_K$ and $\mathscr{X}_{K'}$, respectively,
are equivalent if there are $k(s)$-homomorphisms of $K, K'$ to a third extension $K''$ of $k(s)$ such that $\mathcal{E}_K \otimes_{k(s)} K''$ and $\mathcal{E}'_{K'}\otimes_{k(s)} K''$ are isomorphic.
\begin{definition}
A \emph{set-theoretic family} of coherent sheaves on $p\colon \mathscr{X}\to S$ is a set $\mathscr{F}$ of coherent sheaves defined on the fibres of $p.$
A set-theoretic family $\mathscr{F}$ of coherent sheaves on the fibres of $p$ is {\em bounded} if there is a $S$-scheme $T$ of finite type and a coherent sheaf $\mathcal{H}$ on $\mathscr{X}_T:=\mathscr{X}\times_S T$ such that the set $\mathscr{F}$ is contained in $\{\mathcal{H}_{|\mathscr{X}\times_S \mathrm{Spec}(k(t))}\mid t\in T\}.$
\end{definition}
\begin{proposition}[{\cite[Cor.~4.17]{art:nironi2008}}]\label{prop:boundedness}
A set-theoretic family $\mathscr{F}$ of coherent sheaves on the fibres of $p$ is bounded if and only if the set-theoretic family $F_{\mathcal{G}}(\mathscr{F})$ is bounded.
\end{proposition}
 
In \cite{art:nironi2008}, Sect.~4.1, Nironi proved a stacky version of Kleiman's criterion (\cite{art:nironi2008}, Thm.~4.12). In particular, for a set-theoretic family $\mathscr{F}$ of coherent sheaves on the fibres of $p$   the set of   Hilbert polynomials $P_{\mathcal{G}_K}(\mathcal{E}_K)$, for $\mathcal{E}_K\in\mathscr{F}$, is finite and there exist integers $N,m$ such that any coherent sheaf $\mathcal{E}_K\in \mathscr{F}$ is a quotient of $\left(\mathcal{G}^{\oplus N}\otimes \pi^\ast(\mathcal{O}_X(-m))\right)_K$ for any $K$-point of $S.$ The integer $m$ is exactly the regularity of $\mathcal{E}_K$ (the regularity of a coherent sheaf on $\mathscr{X}$ is by definition the regularity of its image on $X$ via the functor $F_{\mathcal{G}}$). 

Let $\mathcal{E}$ be an $S$-flat coherent sheaf on $\mathscr{X}.$ We can look at $\mathcal{E}$ as a bounded set-theoretic family on the fibres of $p\colon \mathscr{X}\to S.$ Moreover, if $S$ is connected, the   Hilbert polynomials of the fibres of $\mathcal{E}$ are constant as a function of $s\in S$ (cf.\ \cite{art:nironi2008}, Lemma 3.16).
\begin{proposition}\label{prop:openpure}
Let $p\colon \mathscr{X}\to S$ be a family of projective stacks with relative polarization $(\mathcal{G},$ $\mathcal{O}_X(1)).$ Let $\mathcal{E}$ be an $S$-flat $d$-dimensional coherent sheaf on $\mathscr{X}$ with fixed Hilbert polynomial $P$ of degree $d.$ The set $\{s\in S\mid \mathcal{E}_s \mbox{ is pure of dimension } d\}$ is open in $S.$
\end{proposition}
\proof
The proof is the same as in the case of projective schemes (\cite{book:huybrechtslehn2010}, Prop.~2.3.1): one uses the stacky version of the Grothendieck Lemma (\cite{art:nironi2008}, Lemma 4.13) and the projectivity of the Quot scheme for coherent sheaves on stacks \cite{art:olssonstarr2003}.
\endproof

\bigskip\section{Framed sheaves on projective stacks}\label{Framed sheaves on projective stacks}

In this section we start the study of $\delta$-(semi)stable framed sheaves on projective stacks. Most of our results are straightforward generalizations of those holding for framed sheaves on smooth projective varieties \cite{art:huybrechtslehn1995-I, art:huybrechtslehn1995-II}. We refer to these papers as main references for framed sheaves on schemes.

\subsection{Preliminaries}

Let $\mathscr{X}$ be a projective stack of dimension $d$ with coarse moduli scheme $\mathscr{X}\xrightarrow{\pi} X.$ Let $(\mathcal{G}, \mathcal{O}_X(1))$ be a polarization on $\mathscr{X}.$
Fix a coherent sheaf $\mathcal{F}$ on $\mathscr{X}$ and a polynomial 
\begin{equation}
\delta(n):=\delta_1 \frac{n^{d-1}}{(d-1)!}+\delta_2 \frac{n^{d-2}}{(d-2)!} +\cdots+\delta_d \in\mathbb{Q}[n]
\end{equation} 
with $\delta_1>0.$ We call $\mathcal{F}$ a \emph{framing sheaf} and $\delta$ a \emph{stability polynomial}.
\begin{definition}
A \emph{framed sheaf} on $\mathscr{X}$ is a pair $\mathfrak{E}:=(\mathcal{E}, \phi_{\mathcal{E}})$, where $\mathcal{E}$ is a coherent sheaf on $\mathscr{X}$ and $\phi_{\mathcal{E}}\colon  \mathcal{E}\to \mathcal{F}$ is a morphism of sheaves. We call $\phi_{\mathcal{E}}$ a \emph{framing} of $\mathcal{E}.$
\end{definition}
First note that the pair $F_{\mathcal{G}}(\mathfrak{E}):=(F_{\mathcal{G}}(\mathcal{E}), F_{\mathcal{G}}(\phi_{\mathcal{E}})\colon F_{\mathcal{G}}(\mathcal{E})\to F_{\mathcal{G}}(\mathcal{F}))$ is a framed sheaf on $X.$ Moreover, since $F_{\mathcal{G}}$ is an exact functor (cf.\ Remark \ref{rem:exactness}), we have $\ker (F_{\mathcal{G}}(\phi_{\mathcal{E}}))=F_{\mathcal{G}}(\ker(\phi_{\mathcal{E}}))$ and $\mathrm{Im}\,(F_{\mathcal{G}}(\phi_{\mathcal{E}}))=F_{\mathcal{G}}(\mathrm{Im}\,(\phi_{\mathcal{E}})).$ Therefore by Lemma \ref{lem:support}, $F_{\mathcal{G}}(\phi_{\mathcal{E}})$ is zero if and only if $\phi_{\mathcal{E}}$ is zero.

For any framed sheaf $\mathfrak{E}=(\mathcal{E}, \phi_{\mathcal{E}})$, its dimension, Hilbert polynomial, multiplicity, rank and hat-slope are just the corresponding quantities for its underlying coherent sheaf $\mathcal{E}.$

Define the function $\epsilon(\phi_{\mathcal{E}})$ by
\begin{equation*}
\epsilon(\phi_{\mathcal{E}}) := \left\{ \begin{array}{ll}
         1 & \mbox{if $\phi_{\mathcal{E}}\neq 0$}\ ,\\
        0 & \mbox{if $\phi_{\mathcal{E}}=0$}\ .\end{array} \right.
\end{equation*}
The \emph{framed Hilbert polynomial} of $\mathfrak{E}$ is
\begin{equation}
P_{\mathcal{G}}(\mathfrak{E},n):= P_{\mathcal{G}}(\mathcal{E},n)-\epsilon(\phi_{\mathcal{E}})\delta(n)\ ,
\end{equation}
and its \emph{reduced framed Hilbert polynomial} is
\begin{equation}
p_{\mathcal{G}}(\mathfrak{E},n):=\frac{P_{\mathcal{G}}(\mathfrak{E},n)}{\alpha_{\mathcal{G},\dim(\mathfrak{E})}(\mathcal{E})}\ .
\end{equation}
The \emph{framed hat-slope} of $\mathfrak{E}=(\mathcal{E}, \phi_{\mathcal{E}})$ is
\begin{equation}
\hat{\mu}_{\mathcal{G}}(\mathfrak{E}):=\hat{\mu}_{\mathcal{G}}(\mathcal{E})-\frac{\epsilon(\phi_{\mathcal{E}})\delta_1}{\alpha_{\mathcal{G},\dim(\mathfrak{E})}(\mathcal{E})}\ .
\end{equation} 
If $\mathcal{E}'$ is a subsheaf of $\mathcal{E}$ with quotient $\mathcal{E}'':=\mathcal{E}/\mathcal{E}'$, the framing $\phi_{\mathcal{E}}$ induces framings $\phi_{\mathcal{E}'}:={\phi_{\mathcal{E}}}_{|\mathcal E'}$ on $\mathcal{E}'$ and $\phi_{\mathcal{E}''}$ on $\mathcal{E}''$, where the framing $\phi_{\mathcal{E}''}$ is defined as  $\phi_{\mathcal{E}''}=0$ if $\phi_{\mathcal{E}'}\neq 0$; otherwise, $\phi_{\mathcal{E}''}$ is the induced morphism on $\mathcal{E}''.$ 
If $\mathfrak{E}=(\mathcal{E},\phi_{\mathcal{E}})$ is a framed sheaf on $\mathscr{X}$ and $\mathcal{E}'$ is a subsheaf of $\mathcal{E}$, we denote by $\mathfrak{E}'$ the framed sheaf $(\mathcal{E}', \phi_{\mathcal{E}'})$ and by $\mathfrak{E}''$ the framed sheaf $(\mathcal{E}'', \phi_{\mathcal{E}''}).$ 
With this convention the framed Hilbert polynomial of $\mathfrak{E}$ behaves additively:
\begin{equation}
P_{\mathcal{G}}(\mathfrak{E})=P_{\mathcal{G}}(\mathfrak E')+P_{\mathcal{G}}(\mathfrak E'')\ .
\end{equation}
The same property holds for the framed hat-slope.

Thus  there are canonical framings on subsheaves and quotients. The same happens for subquotients, indeed we have the following result
(\cite{art:huybrechtslehn1995-II}, Lemma 1.12).
\begin{lemma}
Let $\mathcal{E}_2\subset \mathcal{E}_1\subset \mathcal{E}$ be coherent sheaves and $\phi_{\mathcal{E}}$ a framing of $\mathcal{E}.$ Then the framings induced on $\mathcal{E}_1/\mathcal{E}_2$ as a quotient of $\mathcal{E}_1$ and as a subsheaf of $\mathcal{E}/\mathcal{E}_2$ agree.
\end{lemma}
Now we introduce the notion of a morphism of framed sheaves.
\begin{definition}\label{def:morphism}
 A \emph{morphism of framed sheaves} $f\colon \mathfrak{E}\to \mathfrak{H}$ is a morphism of the underlying coherent sheaves $f\colon  \mathcal{E} \to \mathcal{H}$ for which there is an element $\lambda\in k$ such that $\phi_{\mathcal{H}} \circ f=\lambda\phi_{\mathcal{E}}.$ We say that $f$ is injective (resp.\ surjective) if the morphism $f\colon  \mathcal{E} \to \mathcal{H}$ is injective (resp.\ surjective). If $f$ is injective, we  call $\mathfrak{E}$ a \emph{framed submodule} of $\mathfrak{H}.$ If $f$ is surjective, we  call $\mathfrak{H}$ a \emph{framed quotient module} of $\mathfrak{E}.$
\end{definition}
\begin{lemma}[{\cite[Lemma 1.5]{art:huybrechtslehn1995-II}}]
\label{lem:morphisms}
 The set $\operatorname{Hom}(\mathfrak{E}, \mathfrak{H})$ of morphisms of framed sheaves is a linear subspace of $\operatorname{Hom} (\mathcal{E},\mathcal{H}).$ If $f\colon  \mathfrak{E}\to \mathfrak{H}$ is an isomorphism,   the factor $\lambda$  can be taken in $k^{\ast}.$ In particular, the isomorphism $f_0=\lambda^{-1}f$ satisfies $\phi_{\mathcal{H}}\circ f_0=\phi_{\mathcal{E}}.$
\end{lemma}
\begin{remark}\label{rem:morphisms}
Let us consider the cartesian diagram
\begin{equation}
  \begin{tikzpicture}[xscale=1.5,yscale=-1.2]
    \node (A0_0) at (0, 0) {$W$};
    \node (A0_2) at (2, 0) {$k$};
    \node (A1_1) at (1, 1) {$\square$};
    \node (A2_0) at (0, 2) {$\operatorname{Hom} (\mathcal{E},\mathcal{H})$};
    \node (A2_2) at (2, 2) {$\operatorname{Hom} (\mathcal{E},\mathcal{F})$};
    \path (A0_0) edge [->]node [auto] {$\scriptstyle{}$} (A2_0);
    \path (A0_0) edge [->]node [auto] {$\scriptstyle{}$} (A0_2);
    \path (A0_2) edge [->]node [auto] {$\scriptstyle{\cdot \phi_{\mathcal{E}}}$} (A2_2);
    \path (A2_0) edge [->]node [auto] {$\scriptstyle{\phi_{\mathcal{H}}\circ}$} (A2_2);
  \end{tikzpicture} 
\end{equation}
Then if $\phi_{\mathcal{E}}\neq 0$, one has $W\simeq\operatorname{Hom}(\mathfrak{E}, \mathfrak{H})$; otherwise $W\simeq\operatorname{Hom}(\mathfrak{E}, \mathfrak{H})\times k.$ 
\end{remark}

\subsection{Semistability}

We  use the following convention: if the word ``(semi)\-stable'' occurs in any statement in combination with the symbol $(\leq)$,   two variants of the statement are understood at the same time: a ``semistable'' one involving the relation ``$\leq$'' and a ``stable'' one involving the relation ``$<$''.

We  give a definition of $\delta$-(semi)stability for $d$-dimensional framed sheaves.
\begin{definition}\label{def:semi}
A $d$-dimensional framed sheaf $\mathfrak{E}=(\mathcal{E}, \phi_{\mathcal{E}})$ is said to be $\delta$-\emph{(semi)stable}  if and only if the following conditions are satisfied:
\begin{itemize}\setlength{\itemsep}{2pt}
\item[(i)] $P_{\mathcal{G}}(\mathcal{E}')\; (\leq)\; \alpha_{\mathcal{G},d}(\mathcal{E}') p_{\mathcal{G}}(\mathfrak{E})$ for all subsheaves $\mathcal{E}'\subseteq \ker\phi_{\mathcal{E}}$,
\item[(ii)] $(P_{\mathcal{G}}(\mathcal{E}')-\delta)\; (\leq)\; \alpha_{\mathcal{G},d}(\mathcal{E}') p_{\mathcal{G}}(\mathfrak{E})$ for all subsheaves $\mathcal{E}'\subset \mathcal{E}.$
\end{itemize}
\end{definition}
By using the same arguments as in the proof of Lemma 1.2 in \cite{art:huybrechtslehn1995-II}, one can prove the following.
\begin{lemma}
Let $\mathfrak{E}=(\mathcal{E},\phi_{\mathcal{E}})$ be a $d$-dimensional framed sheaf. If $\mathfrak{E}$ is $\delta$-semistable,  then $\ker\phi_{\mathcal{E}}$ is torsion-free.
\end{lemma}
%\proof
%Let $T:=T_{d-1}(\ker\phi_{\mathcal{E}}).$ By the semistability condition (i), we get
%\begin{equation}
% P_{\mathcal{G}}(T)\leq \alpha_{\mathcal{G},d}(T) p_{\mathcal{G}}(\mathfrak{E}).
%\end{equation}
%Since $\dim(T)\leq d-1$, we get $\alpha_{\mathcal{G},d}(T)=0$, therefore $P_{\mathcal{G}}(T)\le 0.$ On the other hand, if $T\neq 0$, the leading coefficient of $P_{\mathcal{G}}(T)$ is positive. Thus we get a contradiction and therefore $T=0$.
%\endproof
\begin{definition}
Let $\mathfrak{E}=(\mathcal{E},\phi_{\mathcal{E}})$ be a framed sheaf with $\alpha_{\mathcal{G},d}(\mathcal{E})=0.$ If $\phi_{\mathcal{E}}$ is injective, we say that $\mathfrak{E}$ is \emph{semistable} (indeed, in this case, the semistability of the framed sheaf $\mathfrak{E}$ does not depend on $\delta$). Moreover, if $P_{\mathcal{G}}(\mathcal{E})=\delta$ we say that $\mathfrak{E}$ is $\delta$-\emph{stable}.\end{definition}
\begin{definition}\label{def:musemi}
A $d$-dimensional framed sheaf $\mathfrak{E}=(\mathcal{E},\phi_{\mathcal{E}})$ is $\hat{\mu}$-(semi)stable with respect to $\delta_1$ if and only if $\ker\phi_{\mathcal{E}}$ is torsion-free and the following conditions are satisfied:
\begin{itemize}\setlength{\itemsep}{2pt}
\item[(i)] $\alpha_{\mathcal{G},d-1}(\mathcal{E}')\; (\leq)\; \alpha_{\mathcal{G},d}(\mathcal{E}')\hat{\mu}_{\mathcal{G}}(\mathfrak{E})$ for all subsheaves $\mathcal{E}'\subseteq \ker\phi_{\mathcal{E}}$,
\item[(ii)] $\alpha_{\mathcal{G},d-1}(\mathcal{E}')-\delta_1\; (\leq)\; \alpha_{\mathcal{G},d}(\mathcal{E}')\hat{\mu}_{\mathcal{G}}(\mathfrak{E})$ for all subsheaves $\mathcal{E}'\subset \mathcal{E}$ with $\alpha_{\mathcal{G},d}(\mathcal{E}')< \alpha_{\mathcal{G},d}(\mathcal{E}).$
\end{itemize}
\end{definition}
\begin{definition}
Let $\mathfrak{E}=(\mathcal{E},\phi_{\mathcal{E}})$ be a framed sheaf with $\alpha_{\mathcal{G},d}(\mathcal{E})=0.$ If $\phi_{\mathcal{E}}$ is injective, we say that $\mathfrak{E}$ is $\hat{\mu}$-semistable (for sheaves of dimension $\leq d-1$, the definition of $\mu$-semistability of the corresponding framed sheaves does not depend on $\delta_1$). Moreover, if $\alpha_{\mathcal{G},d-1}(\mathcal{E})=\delta_1$, we say that $\mathfrak{E}$ is $\hat{\mu}$-stable with respect to $\delta_1.$
\end{definition}

One has the usual implications among different stability properties of a $d$-dimensional framed sheaf:
\begin{equation}
\hat{\mu}\mbox{-stable}\Rightarrow \mbox{stable} \Rightarrow \mbox{semistable}\Rightarrow \hat{\mu}\mbox{-semistable}\ .
\end{equation}
The following result can be proved as Lemma 1.6 in \cite{art:huybrechtslehn1995-I}.
\begin{lemma}\label{lem:morphism-stability}
Let $\mathfrak{E}=(\mathcal{E},\phi_{\mathcal{E}})$ and $\mathfrak{H}=(\mathcal{H},\phi_{\mathcal{H}})$ be two $d$-dimensional framed sheaves with the same reduced framed   Hilbert polynomial $p.$
\begin{itemize}\setlength{\itemsep}{2pt}
\item[(i)] If $\mathfrak{E}$ is semistable and $\mathfrak{H}$ is stable,   any nonzero morphism $f\colon \mathfrak{E}\to \mathfrak{H}$ is surjective. 
\item[(ii)] If $\mathfrak{E}$ is stable and $\mathfrak{H}$ is semistable,   any nonzero morphism $f\colon \mathfrak{E}\to \mathfrak{H}$ is injective.
\item[(iii)] If $\mathfrak{E}$ and $\mathfrak{H}$ are stable,   any nonzero morphism $f\colon \mathfrak{E}\to \mathfrak{H}$ is an isomorphism. Moreover, $\operatorname{Hom}(\mathfrak{E},\mathfrak{H})\simeq k.$ If in addition $\phi_{\mathcal{E}}\neq 0$ or, equivalently, $\phi_{\mathcal{H}}\neq 0$,   there is a unique isomorphism $f_0$ with $\phi_{\mathcal{H}}\circ f_0=\phi_{\mathcal{E}}.$ 
\end{itemize}
\end{lemma}
\begin{remark}\label{rem:trick}
In the unframed case, (semi)stable $d$-dimensional sheaves are torsion free. On the other hand, $\delta$-(semi)stable framed sheaves may contain torsion subsheaves. However, an easy trick allows one to make use of results about torsion-free sheaves. Since $\mathscr{X}$ is a projective stack, it has the resolution property. Fix a locally free sheaf $\hat{\mathcal{F}}$ and a surjective morphism $\phi\colon \hat{\mathcal{F}}\to \mathcal{F}$ (here $\mathcal F$ is the framing sheaf). Let $\mathcal{B}$  be the corresponding kernel.  With any $d$-dimensional framed sheaf $\mathfrak{E}=(\mathcal{E},\phi_{\mathcal{E}})$ we can associate a commutative diagram with exact rows and columns:
\begin{equation}
  \begin{tikzpicture}[xscale=2.2,yscale=-1.3]
  \node (A0_2) at (2, 0) {$0$};
  \node (A0_3) at (3, 0) {$0$};
  \node (A1_2) at (2, 1) {$\ker\phi_{\mathcal{E}}$};
  \node (A1_3) at (3, 1) {$\ker\phi_{\mathcal{E}}$}; 
\node (A2_0) at (0, 2) {$0$};
\node (A2_1)at (1, 2) {$\mathcal{B}$};
\node (A2_2) at (2, 2) {$\hat{\mathcal{E}}$};
\node (A2_3) at (3, 2) {$\mathcal{E}$};
\node (A2_4) at (4, 2) {$0$}; 
\node (A3_0) at (0, 3) {$0$};
\node (A3_1)at (1, 3) {$\mathcal{B}$};
\node (A3_2) at (2, 3) {$\hat{\mathcal{F}}$};
\node (A3_3) at (3, 3) {$\mathcal{F}$};
\node (A3_4) at (4, 3) {$0$};
    \path (A0_2) edge [->]node [left] {$\scriptstyle{}$} (A1_2);
    \path (A0_3) edge [->]node [auto] {$\scriptstyle{}$} (A1_3);
    \path (A1_2) edge [thin, double distance=1.5pt]node [auto] {$\scriptstyle{}$} (A1_3);
    \path (A1_2) edge [->]node [auto] {$\scriptstyle{}$} (A2_2);
    \path (A1_3) edge [->]node [auto] {$\scriptstyle{}$} (A2_3);
    \path (A2_0) edge [->]node [auto] {$\scriptstyle{}$} (A2_1);
    \path (A2_1) edge [->]node [auto] {$\scriptstyle{}$} (A2_2);    
    \path (A2_2) edge [->]node [auto] {$\scriptstyle{}$} (A2_3);
    \path (A2_3) edge [->]node [auto] {$\scriptstyle{}$} (A2_4);
    \path (A3_0) edge [->]node [auto] {$\scriptstyle{}$} (A3_1);
    \path (A3_1) edge [->]node [auto] {$\scriptstyle{}$} (A3_2);    
    \path (A3_2) edge [->]node [auto] {$\scriptstyle{}$} (A3_3);
    \path (A3_3) edge [->]node [auto] {$\scriptstyle{}$} (A3_4);
    \path (A2_1) edge [thin, double distance=1.5pt]node [auto] {$\scriptstyle{}$} (A3_1);
    \path (A2_2) edge [->]node [auto] {$\scriptstyle{\phi_{\hat{\mathcal{E}}}}$} (A3_2);
    \path (A2_3) edge [->]node [auto] {$\scriptstyle{\phi_{\mathcal{E}}}$} (A3_3);    
  \end{tikzpicture}
\end{equation}
If $\ker\phi_{\mathcal{E}}$ is torsion-free, the torsion subsheaf of $\hat{\mathcal{E}}$ injects into $\hat{\mathcal{F}}.$ Since $\hat{\mathcal{F}}$ is locally free,  $\hat{\mathcal{E}}$ is torsion-free. Obviously, if $\hat{\mathcal{E}}$ is torsion-free,   $\ker\phi_{\mathcal{E}}$ is torsion-free as well. We denote by $\hat{\mathfrak{E}}$ the framed sheaf $(\hat{\mathcal{E}}, \phi_{\hat{\mathcal{E}}}\colon \hat{\mathcal{E}}\to \hat{\mathcal{F}}).$ 
\end{remark}

\subsection{Jordan-H\"older filtration}\label{sec:jh}
The construction of the Jordan-H\"older filtrations does not differ from the case of framed sheaves on smooth projective varieties. Their existence in the case of projective stacks is granted by the fact that $F_{\mathcal{G}}$ is an exact functor and is \emph{compatible} with the torsion filtration (cf.\ Corollary \ref{cor:torsionfiltration}).
\begin{definition}
Let $\mathfrak{E}=(\mathcal{E},\phi_{\mathcal{E}})$ be a $\delta$-semistable $d$-dimensional framed sheaf. A \emph{Jordan-H\"older filtration} of $\mathfrak{E}$ is a filtration
\begin{equation}
\mathcal{E}_{\bullet}: 0=\mathcal{E}_0\subset \mathcal{E}_1\subset\cdots \subset \mathcal{E}_l=\mathcal{E}\ ,
\end{equation}
such that all the factors $\mathcal{E}_i/\mathcal{E}_{i-1}$ together with the induced framings $\phi_i$ are $\delta$-stable with framed   Hilbert polynomial $P_{\mathcal{G}}(\mathcal{E}_i/\mathcal{E}_{i-1},\phi_i)=\alpha_{\mathcal{G},d}(\mathcal{E}_i/\mathcal{E}_{i-1})p_{\mathcal{G}}(\mathfrak{E}).$
\end{definition}
A straightforward generalization of \cite{art:huybrechtslehn1995-II}, Prop.~1.13, yields the following result.
\begin{proposition} Every $\delta$-semistable framed sheaf $\mathfrak E$ admits a 
Jordan-H\"older filtration. The framed sheaf
\begin{equation}
gr(\mathfrak{E})=(gr(\mathcal{E}),gr(\phi_{\mathcal{E}})):=\bigoplus_i (\mathcal{E}_i/\mathcal{E}_{i-1},\phi_i)
\end{equation}
does not depend, up to isomorphism, on the choice of the Jordan-H\"older filtration.
\end{proposition}
%\begin{remark}
%By construction, for $i>0$ all subsheaves $\mathcal{E}_i$ are framed saturated and the framed sheaves $\mathfrak{E}_i$ are semistable with framed   Hilbert polynomial $\alpha_{\mathcal{G}, d}(\mathcal{E}_i)p_{\mathcal{G}}(\mathfrak{E}).$ In particular $\mathfrak{E}_1$ is a stable framed sheaf. Moreover at most one of the framings $\phi_i$ is nonzero and all but possibly one of the factors $\mathcal{E}_i/\mathcal{E}_{i-1}$ are torsion-free sheaves. 
%\end{remark}
\begin{definition}
Two $\delta$-semistable framed sheaves $\mathfrak{E}=(\mathcal{E},\phi_{\mathcal{E}})$ and $\mathfrak{H}=(\mathcal{H},\phi_{\mathcal{H}})$ with the same reduced framed Hilbert polynomial are called \emph{S-equivalent} if their associated graded objects $gr(\mathfrak{E})$ and $gr(\mathfrak{H})$ are isomorphic.
\end{definition}

\subsection{Boundedness}
We introduce the notion of family of framed sheaves and we prove a related boundedness result, which is a stacky version of the one for framed sheaves on smooth projective varieties (cf.\ \cite{art:huybrechtslehn1995-II}, Sect.~2). All base schemes $S$ are of finite type.
 
\begin{definition}
A \emph{flat family} $\mathfrak{E}=(\mathcal{E},L_{\mathcal{E}},\phi_{\mathcal{E}})$ of framed sheaves on $\mathscr  X$  parameterized by a scheme $S$   consists of a coherent sheaf $\mathcal{E}$ on $\mathscr{X}\times S$, flat over $S$, a line bundle $L_{\mathcal{E}}$ on $S$, and a  morphism $\phi_{\mathcal{E}}\colon L_{\mathcal{E}}\to {p_S}_\ast\mathcal{H}om(\mathcal{E},  p_{\mathscr{X}}^\ast\mathcal{F})$ called a \emph{framing} of $\mathcal{E}.$ Two families $\mathfrak{E}=(\mathcal{E},L_{\mathcal{E}},\phi_{\mathcal{E}})$ and $\mathfrak{E}'=(\mathcal{E}',L_{\mathcal{E}'},\phi_{\mathcal{E}'})$ are isomorphic if there exist isomorphisms $g\colon \mathcal{E}\to \mathcal{E}'$ and $h\colon L_{\mathcal{E}}\to L_{\mathcal{E}'}$ such that
\begin{equation}\label{eq:compatibility-isomorphism}
\tilde{g}\circ \phi_{\mathcal{E}}=\phi_{\mathcal{E}'}\circ h\ ,
\end{equation}
where
\begin{equation}
\tilde{g}\colon {p_S}_\ast\mathcal{H}om(\mathcal{E}, p_{\mathscr{X}}^\ast\mathcal{F})\to {p_S}_\ast\mathcal{H}om(\mathcal{E}',   p_{\mathscr{X}}^\ast\mathcal{F})
\end{equation}
is the isomorphism induced by $g$.
\end{definition}
\begin{remark}
We may look at a framing $\phi_{\mathcal{E}}\colon L_{\mathcal{E}}\to {p_S}_\ast\mathcal{H}om(\mathcal{E}, p_{\mathscr{X}}^\ast\mathcal{F})$ as a nowhere vanishing morphism
\begin{equation}
\tilde{\phi}_{\mathcal{E}}\colon p_S^\ast L_{\mathcal{E}}\otimes\mathcal{E}\to p_{\mathscr{X}}^\ast\mathcal{F}\ ,
\end{equation}
defined as the composition
\begin{equation}
p_S^\ast L_{\mathcal{E}}\otimes\mathcal{E}\to p_S^\ast{p_S}_\ast\mathcal{H}om(\mathcal{E}, p_{\mathscr{X}}^\ast\mathcal{F})\otimes \mathcal{E}
\to\mathcal{H}om(\mathcal{E}, p_{\mathscr{X}}^\ast\mathcal{F})\otimes \mathcal{E}\xrightarrow{\mathrm{ev}} p_{\mathscr{X}}^\ast\mathcal{F}\ .
\end{equation}

\end{remark}
We say that the flat family $\mathfrak{E}=(\mathcal{E},L_{\mathcal{E}},\phi_{\mathcal{E}})$ has the property \textbf{P} if for any closed point $s\in S$ the framed sheaf $(\mathcal{E}_s, (\tilde{\phi}_{\mathcal{E}})_s\colon p_s^\ast((L_{\mathcal{E}})_{s})\otimes\mathcal{E}_s\to p_{\mathscr{X}}^\ast(\mathcal{F})_s)$ has the property \textbf{P}, where $p_s\colon \mathrm{Spec}(k(s))\times \mathscr{X}\to \mathrm{Spec}(k(s))$ is the projection.
\begin{definition}\label{def:familyquotients}
Let $\mathfrak{H}=(\mathcal{H},L_{\mathcal{H}},\phi_{\mathcal{H}})$ be a flat family of framed sheaves on $\mathscr{X}$  parameterized by $S.$ A \emph{flat family of framed quotients} of $\mathfrak{H}$ is a flat family of framed sheaves $\mathfrak{E}=(\mathcal{E},L_{\mathcal{E}},\phi_{\mathcal{E}})$ on $\mathscr{X}$  parameterized by $S$ with an epimorphism $q\colon \mathcal{H}\to \mathcal{E}$ and a morphism $\sigma\in\operatorname{Hom}(L_{\mathcal{E}}, L_{\mathcal{H}})$ such that the diagram
\begin{equation}
  \begin{tikzpicture}[xscale=3.5,yscale=-1.3]
\node (A0_0) at (0,0) {$p_S^\ast L_{\mathcal{E}}\otimes \mathcal{H}$}; 
  \node (A0_2) at (2, 0) {$p_S^\ast L_{\mathcal{E}}\otimes \mathcal{E}$};
  \node (A1_1) at (1, 1) {$p_S^\ast L_{\mathcal{H}}\otimes \mathcal{H}$};
\node (A2_2) at (2, 2) {$  p_{\mathscr{X}}^\ast \mathcal{F}$};
    \path (A0_0) edge [->]node [auto] {$\scriptstyle{\mathrm{id}_{p_S^\ast L_{\mathcal{E}}}\otimes q}$} (A0_2);
    \path (A0_0) edge [->]node [left] {$\scriptstyle{p_S^\ast\sigma\otimes \mathrm{id}_{\mathcal{H}}}\ \ $} (A1_1);
    \path (A1_1) edge [->]node [left] {$\scriptstyle{\tilde{\phi}_{\mathcal{H}}}\ \ $} (A2_2);
    \path (A0_2) edge [->]node [auto] {$\scriptstyle{\tilde{\phi}_{\mathcal{E}}}$} (A2_2);
  \end{tikzpicture}
\end{equation}
commutes.
\end{definition}
\begin{remark}
Let $\mathfrak{H}=(\mathcal{H},\phi_{\mathcal{H}})$ be a framed sheaf on $\mathscr{X}$. Given a scheme $S$, by pulling  $\mathfrak{H}$ back to $\mathscr X\times S$ one defines a flat family $(p^\ast_{\mathscr{X}}(\mathcal{H}), \mathcal{O}_S, p^\ast_{\mathscr{X}}(\phi_{\mathcal{H}}))$ parameterized by    $S$. A  flat family of framed quotients of $\mathfrak{H}$ is a flat family of framed sheaves $\mathfrak{E}=(\mathcal{E},L_{\mathcal{E}},\phi_{\mathcal{E}})$ on $\mathscr{X}$  parameterized by $S$ with an epimorphism $q\colon p^\ast_{\mathscr{X}}(\mathcal{H})\to \mathcal{E}$ and a section $\sigma\in\Gamma(S, L_{\mathcal{E}}^{\vee})$ such that the previous diagram commutes. 
\end{remark}

By Proposition \ref{prop:openpure} torsion-freeness is an open property and we get the following Corollary.
\begin{corollary}\label{cor:torsionfreekernel}
Let $\mathfrak{E}=(\mathcal{E},L_{\mathcal{E}},\phi_{\mathcal{E}})$ be a flat family of framed sheaves on $\mathscr{X}$  parameterized by $S.$ The subset of points $s\in S$ for which $\ker\tilde{\phi}_{\mathcal{E}}$ is torsion-free is open in $S.$
\end{corollary}

By arguing along the lines of the proofs of \cite{book:huybrechtslehn2010}, Prop.~2.3.1 and \cite{art:bruzzomarkushevichtikhomirov2010}, Prop.~3.1 and using the stacky version of Grothendieck's lemma (\cite{art:nironi2008}, Lemma 4.13), one can prove the following result.
\begin{proposition}\label{prop:openesshat-mu}
The property of being $\hat{\mu}$-(semi)stable with respect to $\delta_1$ is open in flat families.
\end{proposition}

\subsubsection{Framed version of two technical Lemmas of Le Potier}

We describe here  generalizations to the framed case of two results, which will be useful to compare the semistability of framed sheaves on $\mathscr{X}$ with GIT semistability (cf.\ Thm.~4.4.1 and Prop.~4.4.2 in \cite{book:huybrechtslehn2010}). The first statement allows one to relate the notion of semistability to the number of global sections of framed submodules or framed quotient modules.
As usual, we denote by $h^0(E)$ the dimension of the vector space of global sections of a coherent sheaf $E$ on $X. $ For a framed sheaf $(E,\phi_E)$ on $X$ we denote by $h^0((E,\phi_E)(m))$ the difference $h^0(E (m))-\epsilon(\phi_{E})\delta(m).$

Let $P$ be a numerical polynomial of degree $d$, $r>0$ its leading coefficient, and $\hat{\mu}_P$ the corresponding hat-slope.
\begin{theorem}\label{thm:boundednessframed}
Let  $\mathfrak{E}=(\mathcal{E}, \phi_{\mathcal{E}})$ be a framed sheaf with $\phi_{\mathcal{E}}\neq 0$,   Hilbert polynomial $P$ and $\ker\phi_{\mathcal{E}}$ torsion-free. There is an integer $m_0$ such that the following three properties of a framed sheaf are equivalent for $m\geq m_0$:
\begin{itemize}\setlength{\itemsep}{2pt}
\item[(i)] $\mathfrak{E}=(\mathcal{E}, \phi_{\mathcal{E}})$ is $\delta$-(semi)stable;
\item[(ii)] $P(m)-\delta(m)\leq h^0(F_{\mathcal{G}}(\mathfrak{E})(m))$ and $h^0(F_{\mathcal{G}}(\mathfrak{E}')(m))\;(\leq)\;\frac{r'}{r}(P(m)-\delta(m))$ for all framed submodules $\mathfrak{E}'$ of $\mathfrak{E}$ of multiplicity $r'$, $0\neq \mathcal{E}'\neq \mathcal{E}$;
\item[(iii)] $h^0(F_{\mathcal{G}}(\mathfrak{E}'')(m))\;(\geq)\; \frac{r''}{r}(P(m)-\delta(m))$ for all framed quotient modules $\mathfrak{E}''$ of $\mathfrak{E}$ of multiplicity $r''$, $\mathcal{E}\neq \mathcal{E}''\neq 0.$
\end{itemize}
Moreover, $\mathcal{E}$ is $m$-regular for all $m\geq m_0$.
\end{theorem}
The set-theoretic families of framed sheaves having torsion-free kernel and satisfying the  version of one of the conditions (i)-(iii) with the weak inequality are denoted by $\mathscr{F}^s$, $\mathscr{F}_m'$ and $\mathscr{F}_m''$, respectively. To prove Theorem \ref{thm:boundednessframed} we need a stacky version of the Le Potier-Simpson boundedness theorem (\cite{book:huybrechtslehn2010}, Thm.~3.3.1). By \cite{art:nironi2008}, Rem.~4.6, its proof   is straightforward.
\begin{theorem}\label{thm:lepotiersimpson}
Let $\mathscr{X}$ be a projective stack with polarization $(\mathcal{G}, \mathcal{O}_X(1)).$ For any pure coherent sheaf $\mathcal{E}$ there is an $F_{\mathcal{G}}(\mathcal{E})$-regular sequence of hyperplane sections $\sigma_1, \ldots, \sigma_{\dim(\mathcal{E})}$ such that
\begin{equation}
h^0(\mathscr{X}_\nu, \mathcal{E}\otimes \mathcal{G}^\vee_{|\mathscr{X}_\nu})\leq \frac{1}{\nu!}\left(\left[\hat{\mu}_{max}(F_{\mathcal{G}}(\mathcal{E}))+r^2+\frac{1}{2}(r+\dim(\mathcal{E}))-1\right]_+\right)^\nu,
\end{equation}
for all $\nu=\dim(\mathcal{E}), \ldots, 0$ and $\mathscr{X}_\nu=Z(\pi^\ast(\sigma_1))\cap \cdots \cap Z(\pi^\ast(\sigma_{\dim(\mathcal{E})-\nu}))$, where $Z(\pi^\ast(\sigma_i))\subset\mathscr X $ is the zero locus of $\pi^\ast(\sigma_i)\in H^0(\mathscr{X}, \pi^\ast(\mathcal{O}_X(1)))$ for $i=1, \ldots, \dim(\mathcal{E})$, $r$ is the multiplicity of $F_{\mathcal{G}}(\mathcal{E})$ and $[x]_+:=\max\{0,x\}.$
\end{theorem}
The next Lemma can be proved by the same arguments as in the proof of \cite{art:nironi2008}, Prop.~4.24. 
\begin{lemma}\label{lem:4.24}
Let $\mathscr{X}$ be a projective stack with polarization $(\mathcal{G}, \mathcal{O}_X(1)).$ Let $\mathfrak{E}=(\mathcal{E}, \phi_{\mathcal{E}})$ be a framed sheaf on $\mathscr{X}$ with $\mathcal{E}$ torsion-free. Let $\bar{E}'$ be the maximal $\hat{\mu}$-destabilizing sheaf of $F_{\mathcal{G}}(\mathcal{E}).$ Then there exists a subsheaf $\bar{\mathcal{E}}'$ of $\mathcal{E}$ such that  \begin{equation}
\hat{\mu}_{max}(F_{\mathcal{G}}(\mathcal{E}))=\hat{\mu}(\bar{E}')\leq \hat{\mu}_{\mathcal{G}}(\bar{\mathcal{E}}')+\tilde{m}\deg(X),
\end{equation}
where $\tilde{m}$ is an integer such that $\pi_\ast(\mathcal{E}nd_{\mathcal{O}_{\mathscr{X}}}(\mathcal{G})) (\tilde{m})$ is generated by global sections.
\end{lemma}
\begin{lemma}\label{lem:inequality}
There are integers $C$ and $m_1$ such that for all $d$-dimensional framed sheaves $\mathfrak{E}=(\mathcal{E}, \phi_{\mathcal{E}})$ in the family $\mathscr{F}:=\mathscr{F}^s\cup \bigcup_{m\geq m_1}\mathscr{F}_m''$ and for all framed saturated subsheaves $\mathcal{E}'$ the following holds: $\alpha_{\mathcal{G},d-1}(\mathcal{E}')-\epsilon(\phi_{\mathcal{E'}})\delta_1\leq r'\left(\hat{\mu}_P-\frac{\delta_1}{r}\right)+C$, and either $\alpha_{\mathcal{G},d-1}(\mathcal{E}')-\epsilon(\phi_{\mathcal{E'}})\delta_1\geq r'\left(\hat{\mu}_P-\frac{\delta_1}{r}\right)-C$ or
\begin{equation*}
\left\{
\begin{array}{ll}
h^0(F_{\mathcal{G}}(\mathfrak{E}')(m))\;<\;\frac{r'}{r}(P(m)-\delta(m)) & \mbox{if } \mathfrak{E}\in \mathscr{F}^s \mbox{ and } m\geq m_1; \mbox{ and }\\[5pt]
\frac{r''}{r}(P-\delta)\;<\;P_{\mathcal{G}}(\mathfrak{E}'') & \mbox{if } \mathfrak{E}\in \mathscr{F}_m'' \mbox{ for some } m\geq m_1.
\end{array}\right.
\end{equation*}
Here $r'$ and $r''$ denote the multiplicity of $\mathcal{E}'$ and $\mathcal E''=\mathcal{E}/\mathcal{E}'$ respectively.
\end{lemma}
\proof
The proof is a straightforward generalization of that for framed sheaves on smooth projective varieties (cf.\ \cite{art:huybrechtslehn1995-II}, Lemma 2.4) by combining Theorem \ref{thm:lepotiersimpson} and Lemma \ref{lem:4.24}.
\endproof
\begin{lemma}\label{lem:boundednessframed}
Let $\mathscr{F}_{\rm ker}$ be the family of kernels of framed sheaves in $\mathscr{F}.$ The families $\mathscr{F}$ and $\mathscr{F}_{\rm ker}$ are bounded.
\end{lemma}
\proof
Assume that $\mathfrak{E}=(\mathcal{E}, \phi_{\mathcal{E}})$ belongs to $\mathscr{F}.$ Let $\hat{\mathfrak{E}}=(\hat{\mathcal{E}}, \phi_{\hat{\mathcal{E}}})$ be the framed sheaf obtained as it is explained in Remark \ref{rem:trick}. Then $P_{\mathcal{G}}(\hat{\mathcal{E}})=P_{\mathcal{G}}(\mathcal{E})+P_{\mathcal{G}}(\mathcal{B})$ does not depend of $\mathfrak{E}.$ Moreover, if $\hat{\mathcal{E}}'$ is a nontrivial subsheaf of $\hat{\mathcal{E}}$, let $\mathcal{E}'$ denote its image in $\mathcal{E}$ and $\mathcal{E}_B'=\hat{\mathcal{E}}\cap \mathcal{B}.$ Then 
\begin{equation}
\hat{\mu}_{\mathcal{G}}(\hat{\mathcal{E}}')\leq C+\hat{\mu}_P-\frac{\delta_1}{r}+\hat{\mu}_{\mathcal G max}(\mathcal{B})\ ,
\end{equation}
and the quantity on the right hand side is independent of $\mathfrak{E}.$ Thus by Lemma \ref{lem:4.24}, 
\begin{equation}
\hat{\mu}_{max}(F_{\mathcal{G}}(\hat{\mathcal{E}}))=\hat{\mu}(\bar{E}')\leq \hat{\mu}_{\mathcal{G}}(\bar{\mathcal{E}}')+\tilde{m}\deg(X)\leq C+\hat{\mu}_P-\frac{\delta_1}{r}+\hat{\mu}_{\mathcal G max}(\mathcal{B})+\tilde{m}\deg(X)\ ,
\end{equation}
where $\bar{E}'$ is the maximal $\hat{\mu}$-destabilizing sheaf of $F_{\mathcal{G}}(\hat{\mathcal{E}})$ and $\bar{\mathcal{E}}'\subset \hat{\mathcal{E}}$ is  the corresponding sheaf. Thus by \cite{art:nironi2008}, Thm.~4.27-(1), the family of coherent sheaves $\hat{\mathcal{E}}$ is bounded. Since the sheaves $\mathcal{E}$ are quotients of $\hat{\mathcal{E}}$ with Hilbert polynomial $P$, they form a bounded family, too. Finally, the family of kernels $\ker(\phi_{\mathcal{E}})$ of framed sheaves $\mathfrak{E}=(\mathcal{E}, \phi_{\mathcal{E}})$ is bounded because all the morphisms $\phi_{\mathcal{E}}$ are morphisms between elements of two bounded families (cf.\ \cite{art:grothendieck1995-SB6}, Prop.~1.2-(i)).
\endproof
\proof[Proof of Theorem \ref{thm:boundednessframed}]
By using Lemma \ref{lem:inequality} and Lemma \ref{lem:boundednessframed}, the proof is a straightforward generalization of  the one for framed sheaves on smooth projective varieties (cf.\ \cite{art:huybrechtslehn1995-II}, Thm.~2.1).
\endproof

Now we prove a Lemma which allows us to deal with possibly framed sheaves with non torsion-free kernels. We need an assumption of normality on the Deligne-Mumford stack $\mathscr{X}.$
\begin{lemma}\label{lem:nontorsionfree} 
Let $\mathscr{X}$ be a normal projective stack. If $(\mathcal{E}, \phi_{\mathcal{E}})$ is a framed sheaf on $\mathscr{X}$ that can be deformed to a framed sheaf with torsion-free kernel, there is a morphism $f\colon (\mathcal{E}, \phi_{\mathcal{E}}) \to (\mathcal{H}, \phi_{\mathcal{H}})$ of framed sheaves such that:
 \begin{itemize}\setlength{\itemsep}{3pt}
 \item[(i)]  $(\mathcal{H}, \phi_{\mathcal{H}})$ has torsion-free kernel; 
\item[(ii)]   $P_{\mathcal{G}}(\mathcal{E})=P_{\mathcal{G}}(\mathcal{H})$ and $P_{\mathcal{G}}(\mathcal{E}, \phi_{\mathcal{E}})=P_{\mathcal{G}}(\mathcal{H},\phi_{\mathcal{H}})$; 
\item[(iii)] $\ker(f)=T_{d-1}(\ker\phi_{\mathcal{E}}).$
\end{itemize}
\end{lemma}
\proof
One combines the arguments of \cite{art:huybrechtslehn1995-II}, Lemma 1.11, and \cite{art:nironi2008}, Lemma 6.10.
\endproof

\bigskip\section{Moduli spaces of framed sheaves on projective stacks}\label{Moduli spaces of framed sheaves on projective stacks}

In this section we shall describe a construction of the moduli spaces of $\delta$-(semi)stable framed sheaves on a normal projective stack $\mathscr{X}$. If the framing vanishes, these are just the moduli spaces of (semi)stable torsion-free sheaves, for which we refer to Nironi's paper~\cite{art:nironi2008}. From now on we shall always assume that the framings are nonzero unless the contrary is explicitly stated.  

Let $\mathscr{X}$ be a {$d$-dimensional} projective stack with coarse moduli scheme $\pi\colon \mathscr{X}\to X$. In this section we make the following assumptions on $\mathscr{X}$:
\begin{itemize}\setlength{\itemsep}{2pt}
\item $\mathscr{X}$ is normal (this hypothesis is necessary only to use Lemma \ref{lem:nontorsionfree} in our construction);
\item $\mathscr{X}$ is irreducible. By \cite{art:vistoli1989}, Lem. 2.3, also the coarse moduli scheme $X$ is irreducible. We shall use this hypothesis in the proof of Proposition \ref{prop:Gm-invariance}, which is in turn used to prove that the moduli space of $\delta$-stable framed sheaves is fine.
\end{itemize}

\subsection{GIT}

The construction of the moduli spaces of $\delta$-(semi)stable framed sheaves on $\mathscr{X}$ is  quite involved, hence, for the sake of clarity, we divide it into several steps.
\subsubsection*{Step 1: construction of a ``Quot-like'' scheme that also takes the framing into account.}
By \cite{art:nironi2008}, Prop.~4.20, the functor $F_{\mathcal{G}}$ defines a closed embedding of $\mathrm{Quot}_{\mathscr{X}/k}(\mathcal{E}, P_0)$ into $\mathrm{Quot}_{X/k}(F_{\mathcal{G}}(\mathcal{E}), P_0)$, for any coherent sheaf $\mathcal{E}$ on $\mathscr{X}$ and numerical polynomial $P_0$ of degree $d$. In particular, $\mathrm{Quot}_{\mathscr{X}/k}(\mathcal{E},$ $P_0)$ is a projective scheme.

Let $P_0$ denote a numerical polynomial of degree $d$, $P=P_0-\delta.$ Fix an integer $m\geq m_0$ (notations of Theorem \ref{thm:boundednessframed}) and let $V$ be a vector space of dimension $P_0(m).$ For every sheaf $E$ on $X$ we shall denote $E(-m)=E\otimes\mathcal{O}_X(-m).$ 

Set $\tilde{\mathbf{Q}}:=\mathrm{Quot}_{\mathscr{X}/k}(G_{\mathcal{G}}(V(-m)), P_0)$ and $\mathbf{P}:=\mathbb{P}\left(\operatorname{Hom}(V, H^0(\mathcal{F}\otimes \mathcal{G}^\vee\otimes\pi^\ast \mathcal{O}_X(m)))^\vee\right)\simeq \mathbb{P}\left(\operatorname{Hom}(V, H^0(F_{\mathcal{G}}(\mathcal{F}) (m)))^\vee\right).$ Given a point $[a\colon V\to H^0(F_{\mathcal{G}}(\mathcal{F}) (m))]$ in $\mathbf{P}$ we can define a framing on $G_{\mathcal{G}}(V(-m))$ as follows. Let us consider the composition
\begin{equation}
V(-m)\xrightarrow{a\circ \mathrm{id}} H^0(F_{\mathcal{G}}(\mathcal{F}) (m)) (-m) \xrightarrow{\mathrm{ev}} F_{\mathcal{G}}(\mathcal{F})\ .
\end{equation}
By applying the functor $G_{\mathcal{G}}$  and   composing on the right with $\theta_{\mathcal{G}}(\mathcal{F})$, we obtain
\begin{equation}
\phi_a\colon G_{\mathcal{G}}(V(-m))\xrightarrow{G_{\mathcal{G}}(a\circ \mathrm{id})} H^0(F_{\mathcal{G}}(\mathcal{F}) (m))\otimes G_{\mathcal{G}}(\mathcal{O}_X(-m)) 
\xrightarrow{G_{\mathcal{G}}(\mathrm{ev})} G_{\mathcal{G}}(F_{\mathcal{G}}(\mathcal{F}))\xrightarrow{\theta_{\mathcal{G}}(\mathcal{F})} \mathcal{F}\ .
\end{equation}
Let $i\colon Z'\hookrightarrow \tilde{\mathbf{Q}}\times \mathbf{P}$ be the closed subscheme of points 
\begin{equation}
\left([\tilde{q}\colon G_{\mathcal{G}}(V(-m))\to \mathcal{E}], [a\colon V\to H^0(F_{\mathcal{G}}(\mathcal{F}) (m))]\right)
\end{equation}
such that the framing $\phi_a$ factors through $\tilde{q}$ and induces a framing $\phi_{\mathcal{E}}\colon \mathcal{E}\to \mathcal{F}.$ 

We explain how to define a flat family of framed sheaves on $\mathscr{X}$  parameterized by $Z'\subset \tilde{\mathbf{Q}}\times \mathbf{P}.$  Let $\tilde{\mathbf{q}}\colon p_{\tilde{\mathbf{Q}}\times \mathscr{X}, \mathscr{X}}^\ast G_{\mathcal{G}}(V(-m))\to \tilde{\mathcal{U}}$ be the universal quotient family on $\mathscr{X}$  parameterized by $\tilde{\mathbf{Q}}.$ Set
\begin{equation}
\mathcal{H}:=\left(p_{\tilde{\mathbf{Q}}\times \mathscr{X}, \mathscr{X}}\circ p_{\tilde{\mathbf{Q}}\times \mathbf{P}\times \mathscr{X}, \tilde{\mathbf{Q}}\times \mathscr{X}}\right)^\ast G_{\mathcal{G}}(V(-m))\ .
\end{equation}
Then we have a quotient morphism
\begin{equation}
p_{\tilde{\mathbf{Q}}\times \mathbf{P}\times \mathscr{X}, \tilde{\mathbf{Q}}\times \mathscr{X}}^\ast \tilde{\mathbf{q}}\colon\mathcal{H}\to p_{\tilde{\mathbf{Q}}\times \mathbf{P}\times \mathscr{X}, \tilde{\mathbf{Q}}\times \mathscr{X}}^\ast \tilde{\mathcal{U}} \to 0\ .
\end{equation}
Consider now the universal quotient sheaf of $\mathbf{P}$, that is,
\begin{equation}
\rho\colon \mathcal{H}om(V\otimes \mathcal{O}_{\mathbf{P}}, H^0(F_{\mathcal{G}}(\mathcal{F}) (m))\otimes \mathcal{O}_{\mathbf{P}})\to \mathcal{O}_{\mathbf{P}}(1)\to 0\ .
\end{equation}
By an argument similar to the one used earlier to construct $\phi_a$ from a point $[a]\in\mathbf{P}$, we can define a morphism
\begin{equation}
\phi_{\mathcal{H}}\colon L_{\mathcal{H}}\to {p_{\tilde{\mathbf{Q}}\times \mathbf{P}\times \mathscr{X}, \tilde{\mathbf{Q}}\times\mathbf{P}}}_\ast\mathcal{H}om(\mathcal{H}, p_{\tilde{\mathbf{Q}}\times \mathbf{P}\times \mathscr{X}, \mathscr X}^\ast \mathcal{F})\ ,
\end{equation}
where $L_{\mathcal{H}}:=p_{\tilde{\mathbf{Q}}\times \mathbf{P}, \mathbf{P}}^\ast\mathcal{O}_{\mathbf{P}}(-1)$. In this way, $(\mathcal H,L_{\mathcal H},\phi_{\mathcal H})$ is a flat family of framed sheaves on the stack $\mathscr X$  parameterized by $\tilde{\mathbf{Q}}\times \mathbf{P}$.

We can endow the universal quotient family $\mathcal{U}:=(i\times \mathrm{id}_{\mathscr{X}})^\ast \tilde{\mathcal{U}}$ on $\mathscr{X}$  parameterized by $Z'$ with a framed sheaf structure in the following way. By the definition of $Z'$ there exists a morphism
\begin{equation}
\phi_{\mathcal{U}}\colon L_{\mathcal{U}}\to {p_{Z'\times \mathscr{X}, Z'}}_\ast \mathcal{H}om(\mathcal{U}, p_{Z'\times \mathscr{X}, \mathscr X}^\ast \mathcal{F})\ ,
\end{equation}
where $L_{\mathcal{U}}:=\left(p_{\tilde{\mathbf{Q}}\times \mathbf{P}, \mathbf{P}}\circ i\right)^\ast \mathcal{O}_{\mathbf{P}}(-1)=i^\ast L_{\mathcal{H}}$. 

Set $\mathfrak{U}:=(\mathcal{U}, L_{\mathcal{U}}, \phi_{\mathcal{U}})$. Then by choosing the morphism $\sigma$ in Definition \ref{def:familyquotients} to be $\mathrm{id}_{L_{\mathcal{U}}}$, we obtain the following result.
\begin{proposition}
$\mathfrak{U}$ is a flat family of framed sheaves on $\mathscr{X}$ parameterized by $Z'$, and is formed by framed quotients of the flat family $\mathfrak{H}:=(i^\ast \mathcal{H}, i^\ast L_{\mathcal{H}}, i^\ast \phi_{\mathcal{H}})$ of framed sheaves on $\mathscr{X}$, which is also  parameterized by $Z'.$
\end{proposition}

The schemes  $\tilde{\mathbf{Q}}$ and $\mathbf{P}$ enjoy universality properties so that the same happens for  the scheme $Z'$.
This is proved as in \cite{art:bruzzomarkushevichtikhomirov2010} and \cite{art:huybrechtslehn1995-II}. 
\begin{proposition}\label{prop:factorization}
Let $[a]$ be a point in $\mathbf{P}$, and  let $\mathfrak{E}=(\mathcal{E},L_{\mathcal{E}},\phi_{\mathcal{E}})$ be a flat family of framed quotients of $(G_{\mathcal{G}}(V(-m)), \phi_a).$ Assume that the   Hilbert polynomial of $\mathcal{E}_s$ is independent of $s\in S.$ There is a morphism $f\colon S\to Z'$ (unique up to a unique isomorphism) such that $\mathfrak{E}$ is isomorphic to the pull-back of $\mathfrak{U}$ via $f\times \mathrm{id}.$
\end{proposition}
%\begin{remark}
%Note that the pullback of $L_{\mathcal{H}}$ by $f$ is the trivial line bundle on $S.$ 
%\end{remark}

\subsubsection*{Step 2: $\mathrm{GL}(V)$-action on $Z'$.}

Until now, we constructed a projective scheme $Z'$ which parameterizes a flat family of framed quotients of $G_{\mathcal{G}}(V(-m))$, with its framed sheaf structure. To use the GIT machinery we need to define an action of a reductive group on $Z'.$ We shall endow $Z'$ of a $\mathrm{GL}(V)$-action induced by $\mathrm{GL}(V)$-actions on $\tilde{\mathbf{Q}}$ and $\mathbf{P}.$

Let $\tau\colon V\otimes\mathcal{O}_{\mathrm{GL}(V)}\to V\otimes \mathcal{O}_{\mathrm{GL}(V)}$ be the \emph{universal automorphism} of $V$  parameterized by $\mathrm{GL}(V).$ The composition
\begin{multline}
p_{\tilde{\mathbf{Q}}\times \mathrm{GL}(V)\times \mathscr{X}, \mathscr{X}}^\ast G_{\mathcal{G}}(V(-m))
\xrightarrow{p_{\tilde{\mathbf{Q}}\times \mathrm{GL}(V)\times \mathscr{X}, \mathrm{GL}(V)}^\ast \tau}
 p_{\tilde{\mathbf{Q}}\times \mathrm{GL}(V)\times \mathscr{X}, \mathscr{X}}^\ast G_{\mathcal{G}}(V(-m))\\ \xrightarrow{p_{\tilde{\mathbf{Q}}\times \mathrm{GL}(V)\times \mathscr{X}, \tilde{\mathbf{Q}}\times \mathscr{X}}^\ast \tilde{\mathbf{q}}} p_{\tilde{\mathbf{Q}}\times \mathrm{GL}(V)\times \mathscr{X}, \tilde{\mathbf{Q}}\times \mathscr{X}}^\ast \tilde{\mathcal{U}}\to 0
\end{multline}
is a flat family of quotients on $\mathscr{X}$  parameterized by $\tilde{\mathbf Q}\times \mathrm{GL}(V).$ Therefore, by the universal property of $\tilde{\mathbf{Q}}$ we get a classifying morphism $\xi_1\colon \tilde{\mathbf{Q}}\times \mathrm{GL}(V)\to \tilde{\mathbf{Q}}$, which is just the right $\mathrm{GL}(V)$-action on $\tilde{\mathbf{Q}}$ pointwise defined as $[\tilde{q}]\cdot g:=[\tilde{q}\circ (g\otimes \mathrm{id})]$, where $[\tilde{q}\colon G_{\mathcal{G}}(V(-m))\to \mathcal{E}]$ is a closed point in $\tilde{\mathbf{Q}}$ and $g\in\mathrm{GL}(V).$

On the other hand, consider 
\begin{multline}
\mathcal{H}om(V\otimes \mathcal{O}_{\mathbf{P}\times \mathrm{GL(V)}}, H^0(F_{\mathcal{G}}(\mathcal{F}) (m))\otimes \mathcal{O}_{\mathbf{P}\times\mathrm{GL}(V)})\xrightarrow{p_{\mathbf{P}\times \mathrm{GL}(V), \mathrm{GL}(V)}^\ast \tau}\\ 
\mathcal{H}om(V\otimes \mathcal{O}_{\mathbf{P}\times \mathrm{GL(V)}}, H^0(F_{\mathcal{G}}(\mathcal{F}) (m))\otimes \mathcal{O}_{\mathbf{P}\times\mathrm{GL}(V)}) \xrightarrow{p_{\mathbf{P}\times\mathrm{GL}(V), \mathbf{P}}^\ast \rho} p_{\mathbf{P}\times\mathrm{GL}(V), \mathbf{P}}^\ast \mathcal{O}_{\mathbf{P}}(1)\to 0\ .
\end{multline}
This induces a classifying morphism $\xi_2\colon \mathbf{P}\times \mathrm{GL}(V)\to \mathbf{P}$, which is the right action of $\mathrm{GL}(V)$ on $\mathbf{P}$ given by $[a]\circ g:=[a\circ g]$ for any closed point $[a\colon V\to H^0(F_{\mathcal{G}}(\mathcal{F}) (m))]$ and $g\in\mathrm{GL}(V).$ By the definition of classifying morphisms there are isomorphisms
\begin{align}
&\Lambda_1 \colon (\xi_1\times\mathrm{id}_{\mathscr{X}})^\ast \tilde{\mathcal{U}}\xrightarrow{\sim} p_{\tilde{\mathbf{Q}}\times \mathrm{GL}(V)\times \mathscr{X}, \tilde{\mathbf{Q}}\times \mathscr{X}}^\ast \tilde{\mathcal{U}}\ ,  \\
&\Lambda_2 \colon \xi_2^\ast \mathcal{O}_{\mathbf{P}}(1)\xrightarrow{\sim} p_{\mathbf{P}\times\mathrm{GL}(V), \mathbf{P}}^\ast \mathcal{O}_{\mathbf{P}}(1)\ ,
\end{align}
such that the following diagrams commute:
\begin{equation}
  \begin{tikzpicture}[xscale=9,yscale=-2]
\node (A0_0) at (0,0) {$p_{\tilde{\mathbf{Q}}\times \mathrm{GL}(V)\times \mathscr{X}, \mathscr{X}}^\ast G_{\mathcal{G}}(V(-m))$}; 
  \node (A0_1) at (1, 0) {$(\xi_1\times\mathrm{id}_{\mathscr{X}})^\ast \tilde{\mathcal{U}}$};
  \node (A1_0) at (0, 1) {$p_{\tilde{\mathbf{Q}}\times \mathrm{GL}(V)\times \mathscr{X}, \mathscr{X}}^\ast G_{\mathcal{G}}(V(-m))$};
\node (A1_1) at (1, 1) {$p_{\tilde{\mathbf{Q}}\times \mathrm{GL}(V)\times \mathscr{X}, \tilde{\mathbf{Q}}\times \mathscr{X}}^\ast \tilde{\mathcal{U}}$};
\node (A) at (1.25,0.5) {,};

\node (A0_0') at (0,1.5) {$\mathcal{H}om(V\otimes \mathcal{O}_{\mathbf{P}\times \mathrm{GL(V)}}, H^0(F_{\mathcal{G}}(\mathcal{F}) (m))\otimes \mathcal{O}_{\mathbf{P}\times\mathrm{GL}(V)})$}; 
  \node (A0_1') at (1, 1.5) {$\xi_2^\ast \mathcal{O}_{\mathbf{P}}(1)$};
  \node (A1_0') at (0, 2.5) {$\mathcal{H}om(V\otimes \mathcal{O}_{\mathbf{P}\times \mathrm{GL(V)}}, H^0(F_{\mathcal{G}}(\mathcal{F}) (m))\otimes \mathcal{O}_{\mathbf{P}\times\mathrm{GL}(V)})$};
\node (A1_1') at (1, 2.5) {$p_{\mathbf{P}\times\mathrm{GL}(V), \mathbf{P}}^\ast \mathcal{O}_{\mathbf{P}}(1)$};
\node (A') at (1.25, 2) {.};

    \path (A0_0) edge [->]node [auto] {$\scriptstyle{(\xi_1\times\mathrm{id}_{\mathscr{X}})^\ast \tilde{\mathbf{q}}}$} (A0_1);
    \path (A0_0) edge [->]node [left] {$\scriptstyle{p_{\tilde{\mathbf{Q}}\times \mathrm{GL}(V)\times \mathscr{X}, \mathrm{GL}(V)}^\ast \tau}$} (A1_0);
    \path (A0_1) edge [->]node [left] {$\scriptstyle{\Lambda_1}$} (A1_1);
    \path (A1_0) edge [->]node [auto] {$\scriptstyle{p_{\tilde{\mathbf{Q}}\times \mathrm{GL}(V)\times \mathscr{X}, \tilde{\mathbf{Q}}\times \mathscr{X}}^\ast \tilde{\mathbf{q}}}$} (A1_1);

    \path (A0_0') edge [->]node [auto] {$\scriptstyle{\xi_2^\ast(\rho)}$} (A0_1');
    \path (A0_0') edge [->]node [left] {$\scriptstyle{p_{\mathbf{P}\times \mathrm{GL}(V), \mathrm{GL}(V)}^\ast \tau}$} (A1_0');
    \path (A0_1') edge [->]node [left] {$\scriptstyle{\Lambda_2}$} (A1_1');
    \path (A1_0') edge [->]node [auto] {$\scriptstyle{p_{\mathbf{P}\times\mathrm{GL}(V), \mathbf{P}}^\ast \rho}$} (A1_1');
      \end{tikzpicture}
\end{equation}
One can check that $\Lambda_1$ is a $\mathrm{GL}(V)$-linearization for $\tilde{\mathcal{U}}$ in the sense of Romagny, cf.\ \cite{art:romagny2005}, Example 4.3. In the same way, $\Lambda_2$ is a $\mathrm{GL}(V)$-linearization for $\mathcal{O}_{\mathbf{P}}(1)$. (For coherent sheaves on Deligne-Mumford stacks we shall always use the term ``linearization" in this sense.)

The classifying morphisms $\xi_1, \xi_2$ induce a right action
\begin{equation}\label{eq:action}
\xi\colon \tilde{\mathbf{Q}}\times \mathbf{P}\times \mathrm{GL}(V)\to \tilde{\mathbf{Q}}\times \mathbf{P}\ .
\end{equation}
The closed subscheme $Z'$ is invariant with respect to this action. Thus we have an induced $\mathrm{GL}(V)$-linearization of $\mathcal{U}$ and an induced $\mathrm{GL}(V)$-linearization of $L_{\mathcal{U}}$ which are compatible with $\phi_{\mathcal{U}}$, i.e., the two linearizations satisfy an equation of the form \eqref{eq:compatibility-isomorphism}.

\subsubsection*{Step 3: comparison between GIT (semi)stability and the $\delta$-(semi)stability condition for framed sheaves}

We need to define suitable $\mathrm{SL}(V)$-linearized ample line bundles on $Z'$ which will allow us to deal with GIT (semi)stable points on $Z'$ and compare them to $\delta$-(semi)stable framed sheaves on $\mathscr{X}.$ From now on we consider $\mathrm{SL}(V)$ instead of $\mathrm{GL}(V)$ because the study of the GIT (semi)stable points is easier for the first group.

As it is described in~\cite{art:nironi2008}, Sect.~6.1, one can define line bundles on $\tilde{\mathbf{Q}}$
\begin{equation}
L_\ell:=\det({p_{\tilde{\mathbf{Q}}}}_\ast F_{\mathcal{G}}(\tilde{\mathcal{U}}) (\ell))\ .
\end{equation}
By~\cite{art:nironi2008}, Prop.~6.2, for $\ell$ sufficiently large the line bundles $L_\ell$ are very ample. Moreover, they carry natural  $\mathrm{SL}(V)$-linearizations (cf.\ \cite{art:nironi2008}, Lemma 6.3). Then the ample line bundles 
\begin{equation}
\mathcal{O}_{Z'}(n_1, n_2):=q^\ast_{\tilde{\mathbf{Q}}} L_\ell^{\otimes n_1} \otimes q_{\mathbf{P}}^\ast \mathcal{O}_{\mathbf{P}}(n_2)
\end{equation}
carry natural $SL(V)$-linearizations, where $q_{\tilde{\mathbf{Q}}}$ and $q_{\mathbf{P}}$ are the natural projections from $Z'$ to $\tilde{\mathbf{Q}}$ and $\mathbf{P}$ respectively. As explained in \cite{art:huybrechtslehn1995-II}, Sect.~3, only the ratio $n_2/n_1$ matters, and we choose it to be
\begin{equation}
\frac{n_2}{n_1}:=P(\ell)\frac{\delta(m)}{P(m)}-\delta(\ell)\ ,
\end{equation}
assuming, of course, that $\ell$ is chosen large enough to make this term positive.

To use the GIT machinery we need to compare the GIT (semi)stability with the $\delta$-(semi)sta\-bil\-i\-ty condition for framed sheaves. The results we show in the following are generalizations of those proved in \cite{art:huybrechtslehn1995-II}, Sect.~3, for framed sheaves on smooth projective varieties. The proofs are rather straightforward due to the properties of the functors $F_{\mathcal{G}}$ and $G_{\mathcal{G}}$.

Let $\tilde{q}\colon G_{\mathcal{G}}(V(-m))\to \mathcal{E}$ be a surjective morphism and $V'$ a nontrivial proper linear subspace of $V.$ The vector space $V'$ defines a subsheaf of $\mathcal{E}$ in the following way: if we apply the functor $G_{\mathcal{G}}$ to the inclusion map $i\colon V'(-m)\hookrightarrow V(-m)$, we get an injective morphism $\tilde{i}\colon G_{\mathcal{G}}(V'(-m))\hookrightarrow G_{\mathcal{G}}(V(-m))$ and the image of $\tilde{q}\circ \tilde{i}$ gives a subsheaf of $\mathcal{E}.$ Now we state the following stacky version of \cite{art:huybrechtslehn1995-II}, Prop.~3.1.
\begin{proposition}\label{prop:previous}
For a sufficiently large $\ell$, the point $([\tilde{q}], [a])\in Z'$ is (semi)stable with respect to the linearization of $\mathcal{O}_{Z'}(n_1, n_2)$ if and only if
\begin{equation}
\dim(V')\cdot \left(n_1 P_0(\ell)+n_2\right)\,(\leq)\,\dim(V)\cdot\left(n_1 P_{\mathcal{G}}(\mathcal{E}',\ell)+n_2\epsilon(\phi_{\mathcal{E}'})\right)\ .
\end{equation}
for every  nontrivial proper linear subspace  $V'$ of $V$. Here $\mathcal{E}'\subset \mathcal{E}$ is the subsheaf given by $V'$.
\end{proposition}
By Corollary \ref{cor:torsionfreekernel}  there is an open subscheme $U\subset Z'$ whose  points   represent framed sheaves with torsion-free kernel. We assume that $U$ is nonempty and denote by $Z$ its closure in $Z'.$

Let $\tilde{q}\colon G_{\mathcal{G}}(V(-m))\to \mathcal{E}$ be a morphism representing a point $[\tilde{q}]\in\tilde{\mathbf{Q}}.$ By applying the functor $F_{\mathcal{G}}$ to $\tilde{q}$ and then composing on the left by $\varphi_{\mathcal{G}}(V(-m))$, we obtain
\begin{equation}
V(-m)\xrightarrow{\varphi_{\mathcal{G}}(V(-m))}F_{\mathcal{G}}(G_{\mathcal{G}}(V(-m)))\to F_{\mathcal{G}}(\mathcal{E})\ ,
\end{equation}
and in cohomology we get $q\colon V\to H^0(F_{\mathcal{G}}(\mathcal{E}) (m)).$
\begin{proposition}\label{prop:semistableZ}
For sufficiently large $\ell$, a point $([\tilde{q}],[a])\in Z$ is (semi)stable with respect to the $\mathrm{SL}(V)$-action on $Z$ if and only if the corresponding framed sheaf $(\mathcal{E}, \phi_{\mathcal{E}})$ is $\delta$-(semi)stable and the map $q\colon V\to H^0(F_{\mathcal{G}}(\mathcal{E}) (m))$ induced by $\tilde{q}$ is an isomorphism.
\end{proposition}
\proof
We have seen that given a point $([\tilde{q}],[a])\in Z$ one can  construct a framed sheaf $(\mathcal{E},\phi_{\mathcal{E}})$ and a map $q\colon V\to H^0(F_{\mathcal{G}}(\mathcal{E}) (m))$. On the other hand, if we fix a framed sheaf $(\mathcal{E},\phi_{\mathcal{E}})$ with an isomorphism $q\colon V\to H^0(F_{\mathcal{G}}(\mathcal{E}) (m))$, we obtain a surjective morphism
\begin{equation}
V(-m)\xrightarrow{q\otimes\mathrm{id}} H^0(F_{\mathcal{G}}(\mathcal{E}) (m)) (-m)\xrightarrow{\mathrm{ev}} F_{\mathcal{G}}(\mathcal{E})
\end{equation}
(note that  $F_{\mathcal{G}}(\mathcal{E})(m)$ is $m$-regular, hence it is globally generated). Since $G_{\mathcal{G}}$ is a right-exact functor and $\theta_{\mathcal{G}}(\mathcal{E})$ is surjective by definition,
the morphism 
\begin{equation}
\tilde{q}\colon G_{\mathcal{G}}(V(-m))\to G_{\mathcal{G}}(F_{\mathcal{G}}(\mathcal{E}))\xrightarrow{\theta_{\mathcal{G}}(\mathcal{E})} \mathcal{E}
\end{equation}
is surjective as well. Furthermore, the framing $\phi=\phi_{\mathcal{E}}\circ \tilde{q}$ defines a morphism $a\colon V\to H^0(F_{\mathcal{G}}(\mathcal{F}) (m)).$ Thus we get a point in $Z'$.  

Now the proof is obtained by combining the arguments in \cite{art:huybrechtslehn1995-II}, Prop.~3.2, with those in \cite{art:nironi2008}, Thm.~5.1. The former allow us to compare the two different (semi)stability conditions (thanks also to Theorem \ref{thm:boundednessframed}), while the latter show that   from a $\delta$-(semi)stable framed sheaf $(\mathcal{E},\phi_{\mathcal{E}})$ with an isomorphism $q\colon V\to H^0(F_{\mathcal{G}}(\mathcal{E}) (m))$ we  construct, by means of the procedure described above, a (semi)stable point in $Z$,  and {\em vice versa}. In view of Lemma \ref{lem:nontorsionfree} we need only to deal with framed sheaves with a torsion-free kernel.
\endproof
The following Lemma is obtained by arguing along the lines of the analogous result in the unframed case (cf.\ \cite{book:huybrechtslehn2010}, Lemma 4.3.2).
\begin{lemma}\label{lem:automorphism}
Let $\left([\tilde{q}\colon G_{\mathcal{G}}(V(-m))\to \mathcal{E}], [a\colon V\to H^0(F_{\mathcal{G}}(\mathcal{F}) (m))]\right)$ be a closed point of $Z'$ such that $F_\mathcal{G}(\mathcal{E})(m)$ is globally generated and   $q\colon V\to  H^0(F_{\mathcal{G}}(\mathcal{E}) (m))$ induced by $\tilde{q}$ is an isomorphism. There is a natural injective homomorphism $i\colon \mathrm{Aut}(\mathcal{E},\phi_{\mathcal{E}})\to \mathrm{GL}(V)$ whose image is precisely the stabilizer subgroup $\mathrm{GL}(V)_{([\tilde{q}],[a])}$ of the point $([\tilde{q}],[a])$.
\end{lemma}
%\proof
%Let us consider an automorphism $f\in\mathrm{Aut}(\mathcal{E},\phi_{\mathcal{E}}).$ Since $F_{\mathcal{G}}$ is an exact functor, $F_{\mathcal{G}}(f)$ is an automorphism of $F_{\mathcal{G}}(\mathcal{E}).$ Therefore, we can define the homomorphism $i\colon \mathrm{Aut}(\mathcal{E},\phi_{\mathcal{E}})\to \mathrm{GL}(V)$ by
%\begin{equation}
%V\xrightarrow{q} H^0(F_{\mathcal{G}}(\mathcal{E}) (m))\xrightarrow{H^0(F_{\mathcal{G}}(f)(m))} H^0(F_{\mathcal{G}}(\mathcal{E}) (m))\xrightarrow{q^{-1}} V\ .
%\end{equation}
%As $F_{\mathcal{G}}(\mathcal{E}) (m)$ is globally generated, $i(f)=\operatorname{id}_V$ implies $F_{\mathcal G}(f) = \operatorname{id}_{F_{\mathcal G}(\mathcal E)}$. By Lemma \ref{lem:support} $f=\operatorname{id}_{\mathcal E}$, i.e., $i$ in injective.
% 
%Let $g$ be an element of the stabilizer subgroup $\mathrm{GL}(V)_{([\tilde{q}],[a])}$ of the point $([\tilde{q}],[a]).$ Then there exist $f\colon \mathcal{E}\to \mathcal{E}$ and $\lambda\in k^\ast$ such that
%\begin{equation}
%f\circ \tilde{q}=\tilde{q}\circ(g\otimes \mathrm{id})\qquad\mbox{and}\qquad a\circ g=\lambda a\ .
%\end{equation}
%Since $([\tilde{q}],[a])$ is a point in $Z'$, we get $\phi_{\mathcal{E}}\circ f=\lambda \phi_{\mathcal{E}}.$
%\endproof

\subsubsection*{Step 4: Good and geometric quotients and (semi)stable locus.}

Thanks to the results we proved before, we are ready to use \cite{book:huybrechtslehn2010}, Thm.~4.2.10, which allows us to construct a (quasi-)projective scheme {\em parameterizing} (semi)stable points of $Z$.

We  denote by $\mathfrak{U}^{(s)s}=(\mathcal{U}^{(s)s}, L_{\mathcal{U}^{(s)s}}, \phi_{\mathcal{U}^{(s)s}})$ the universal family of $\delta$-(semi)stable framed sheaves on $\mathscr{X}$  parameterized by $Z^{(s)s}$ induced, through pull-back, by the one  parameterized by $Z'.$

By using \cite{book:huybrechtslehn2010}, Thm.~4.2.10, we get directly the following.
\begin{theorem}
There exists a projective scheme $\mathrm{M}^{ss}=\mathrm{M}^{ss}_{\mathscr{X}/k}(\mathcal{G}, \mathcal{O}_X(1); P_0, \mathcal{F}, \delta)$ and a morphism $\tilde{\pi}\colon Z^{ss}\to \mathrm{M}^{ss}$ such that $\tilde{\pi}$ is a universal good quotient for the $SL(V)$-action on $Z^{ss}.$ Moreover, there is an open subscheme $\mathrm{M}^{s}=\mathrm{M}^{s}_{\mathscr{X}/k}(\mathcal{G}, \mathcal{O}_X(1); P_0, \mathcal{F}, \delta)\subset \mathrm{M}^{ss}$ such that $Z^s=\tilde{\pi}^{-1}(\mathrm{M}^{s})$ and $\tilde{\pi}\colon Z^s\to \mathrm{M}^{s}$ is a universal geometric quotient. Finally, there is a positive integer $\ell$ and a very ample line bundle $\mathcal{O}_{\mathrm{M}^{ss}}(1)$ on $\mathrm{M}^{ss}$ such that $\mathcal{O}_{Z'}(n_1, n_2)^{\otimes \ell}_{|Z^{ss}}\simeq \tilde{\pi}^{}(\mathcal{O}_{\mathrm{M}^{ss}}(1)).$
\end{theorem}

By using the same arguments as in the proof of \cite{art:huybrechtslehn1995-II}, Prop.~3.3, and the semicontinuity theorem for Hom groups of flat families of framed sheaves (Proposition \ref{prop:semicontinuity}), we get the following result.
\begin{proposition}\label{prop:s-equivalence}
Two points $([\tilde{q}], [a])$ and $([\tilde{q}'], [a'])$ in $Z^{ss}$ are mapped to the same point in $\mathrm{M}^{ss}$ if and only if the corresponding framed sheaves are $S$-equivalent.
\end{proposition}

\subsection{The moduli stacks of $\delta$-(semi)stable framed sheaves}

In the previous section we used GIT machinery to construct a good (geometric) quotient $\mathrm M^{(s)s}$ of $Z^{(s)s}$. Now we introduce a {\em moduli} stack associated with $Z^{(s)s}$ and describe its relation with $\mathrm M^{(s)s}$. Let us define the {algebraic stack of finite type}
\begin{equation}
\mathfrak{SM}^{(s)s}=\mathfrak{SM}^{(s)s}_{\mathscr{X}/k}(\mathcal{G}, \mathcal{O}_X(1); P_0, \mathcal{F}, \delta):=[Z^{(s)s}/\mathrm{SL}(V)] \ .
\end{equation}
Note that $\mathfrak{SM}^{s}$ is an open substack of $\mathfrak{SM}^{ss}.$

%\begin{proposition}\label{prop:artinstack}
%$\mathfrak{SM}^{(s)s}$ is an algebraic stack of finite type.
%\end{proposition}
%\proof
%By \cite{book:laumonmoretbailly2000}, Example 4.6.1, since the group scheme $\mathrm{SL}(V)$ is smooth, separated and of finite presentation,  $\mathfrak{SM}^{(s)s}$ is an algebraic stack and $Z^{(s)s}\to \mathfrak{SM}^{(s)s}$ is a smooth presentation. Moreover, since $Z^{(s)s}$ is of finite type,   $\mathfrak{SM}^{(s)s}$ is of finite type as well.
%\endproof

We explain the relation between $\mathfrak{SM}^{(s)s}$ and $\mathrm{M}^{(s)s}.$ First we recall the notion of \emph{good moduli space} for algebraic stacks.
\begin{definition}[{\cite[Def.~3.1]{art:alper2008}}]
A morphism of algebraic stacks  $f\colon \mathscr{X}\to \mathscr{Y}$ is \emph{cohomologically affine} if it is quasi-compact and the functor $f_\ast\colon \mathrm{QCoh}(\mathscr{X})\to\mathrm{QCoh}(\mathscr{Y})$ is exact.
\end{definition}
\begin{definition}[{\cite[Def.~4.1 and 7.1]{art:alper2008}}]\label{def:goodmoduli}
Let $f\colon \mathscr{X}\to Y$ be a morphism where $\mathscr{X}$ is an algebraic stack and $Y$ an algebraic space. We say that $f$ is a \emph{good moduli space} if the following properties are satisfied:
\begin{itemize}\setlength{\itemsep}{2pt}
\item $f$ is cohomologically affine,
\item the natural morphism $\mathcal{O}_Y\to f_\ast(\mathcal{O}_{\mathscr{X}})$ is an isomorphism.
\end{itemize}
Moreover, a good moduli space $f$  is a \emph{tame moduli space} if the map $[\mathscr{X}(\mathrm{Spec}(k))]\to Y(\mathrm{Spec}(k))$ is a bijection of sets, where $[\mathscr{X}(\mathrm{Spec}(k))]$ denotes the set of isomorphism classes of objects of $\mathscr{X}(\mathrm{Spec}(k)).$
\end{definition}

Since the ample line bundle $\mathcal{O}_{Z'}(n_1, n_2)_{|Z^{ss}}$ is $\mathrm{SL}(V)$-equivariant, it descends to a line bundle $\mathcal{O}(n_1,n_2)$ on $\mathfrak{SM}^{ss}.$ The morphism $\tilde{\pi}$ induces a morphism $\pi_{\mathfrak{S}}\colon \mathfrak{SM}^{ss}\to \mathrm{M}^{ss}.$ By \cite{art:alper2008}, Thm.~13.6 (which is a stacky version of \cite{book:huybrechtslehn2010}, Thm.~4.2.10), we get the following  result.
\begin{theorem}\label{thm:maintheorem}
The morphism $\pi_{\mathfrak{S}}\colon \mathfrak{SM}^{ss}\to \mathrm{M}^{ss}$ is a good moduli space and $\pi_{\mathfrak{S}}^\ast(\mathcal{O}_{\mathrm{M}^{ss}}(1))\simeq \mathcal{O}(n_1,n_2)^{\otimes \ell}.$ Moreover, the morphism
$\pi_{\mathfrak{S}}\colon \mathfrak{SM}^{s}\to \mathrm{M}^{s}$
is a tame moduli space. 
\end{theorem}
Furthermore, by \cite{art:alper2008}, Thm.~6.6, we can state the following universal property for $\pi_{\mathfrak{S}}\colon $ $\mathfrak{SM}^{ss}\to \mathrm{M}^{ss}$.
\begin{proposition}\label{prop:universalproperty}
Let $T$ be an algebraic space and $f\colon \mathfrak{SM}^{ss}\to T$ a morphism. There exists a unique morphism $g\colon \mathrm{M}^{ss}\to T$ such that $f=g\circ \pi_{\mathfrak{S}}.$
\end{proposition}

We introduce two more algebraic stacks of finite type
\begin{align}
\mathfrak{M}^{(s)s}&=\mathfrak{M}^{(s)s}_{\mathscr{X}/k}(\mathcal{G}, \mathcal{O}_X(1); P_0, \mathcal{F}, \delta):=[Z^{(s)s}/\mathrm{GL}(V)]\ ,\\ \mathfrak{PM}^{(s)s}&=\mathfrak{PM}^{(s)s}_{\mathscr{X}/k}(\mathcal{G}, \mathcal{O}_X(1); P_0, \mathcal{F}, \delta):=[Z^{(s)s}/\mathrm{PGL}(V)]\ .
\end{align}
%For the same arguments as in Proposition \ref{prop:artinstack}, these are algebraic stacks of finite type. 
Note that the stack $\mathfrak{PM}^{(s)s}$ is well defined as the multiplicative group $\mathbb{G}_m$ is contained in the stabilizer of every point of $Z^{ss}$ (cf.\ Lemma \ref{lem:automorphism}).

A natural question  is if there is a relation between the stacks $\mathfrak{SM}^{(s)s}$, $\mathfrak{M}^{(s)s}$ and $\mathfrak{PM}^{(s)s}.$ First, note that the smooth groupoid of the smooth presentation $Z^{(s)s}\to \mathfrak{M}^{(s)s}$ is
\begin{equation}
  \begin{tikzpicture}[xscale=2.8,yscale=-.7, ->, bend angle=25]
\node (A0_1) at (0,1) {$Z^{(s)s}\times \mathrm{GL}(V)$};
\node (A1_1) at (1.15,1) {$Z^{(s)s}$}; 
\node (B1) at (0.5,0.8) {$ $};
\node (B2) at (1,0.8) {$ $};
\node (C1) at (0.5,1.2) {$ $};
\node (C2) at (1,1.2) {$ $};
\node (Comma) at (1.4,1.25) {,};
    \path (B1) edge [->]node [auto] {$\scriptstyle{a}$} (B2);
        \path (C1) edge [->]node [below] {$\scriptstyle{p_{Z^{(s)s}}}$} (C2);
  \end{tikzpicture}
\end{equation}
where $a$ is the action morphism of $\mathrm{GL}(V)$ on $Z^{(s)s}$. Since $\mathbb{G}_m$ acts on $Z^{(s)s}\times \mathrm{GL}(V)$ by leaving $a$ and $p_{Z^{(s)s}}$ invariant, we can \emph{rigidify} the smooth groupoid (the notion of \emph{rigidification} is explained in \cite{art:abramovichcortivistoli2003}, Sect.~5) to get
\begin{equation}
  \begin{tikzpicture}[xscale=2.8,yscale=-.7, ->, bend angle=25]
\node (A0_1) at (0,1) {$Z^{(s)s}\times \mathrm{PGL}(V)$};
\node (A1_1) at (1.15,1) {$Z^{(s)s}$}; 
\node (B1) at (0.5,0.8) {$ $};
\node (B2) at (1,0.8) {$ $};
\node (C1) at (0.5,1.2) {$ $};
\node (C2) at (1,1.2) {$ $};
\node (Period) at (1.4,1.25) {.};
    \path (B1) edge [->]node [auto] {$\scriptstyle{a}$} (B2);
        \path (C1) edge [->]node [below] {$\scriptstyle{p_{Z^{(s)s}}}$} (C2);
  \end{tikzpicture}
\end{equation}
This is the smooth groupoid  of $\mathfrak{PM}^{(s)s}.$ In particular, $\mathfrak{M}^{(s)s}\to \mathfrak{PM}^{(s)s}$ is a $\mathbb{G}_{m}$-gerbe {(i.e.,  it is a gerbe for which the automorphism group of every section is  $\mathbb{G}_{m}$, and this isomorphism is compatible with the fibred structure of the gerbe)}. On the other hand, we can rigidify the stack $\mathfrak{SM}^{(s)s}$ with respect to the group $\mu(V)\subset \mathrm{SL}(V)$, where $\mu(V)$ is the group of $\dim(V)$-roots of unity, and we get that the rigidification is isomorphic to $\mathfrak{PM}^{(s)s}.$ Hence $\mathfrak{SM}^{(s)s}\to \mathfrak{PM}^{(s)s}$ is a $\mu(V)$-gerbe.

The morphism $\pi_{\mathfrak{S}}\colon \mathfrak{SM}^{(s)s}\to \mathrm{M}^{(s)s}$ induces a morphism $\pi_{\mathfrak{P}}\colon \mathfrak{PM}^{(s)s}\to \mathrm{M}^{(s)s}$ (cf.\ \cite{art:abramovichcortivistoli2003}, Thm.~5.1.5-(2)), so that we get a morphism $\pi\colon \mathfrak{M}^{(s)s}\to \mathrm{M}^{(s)s}$ and the following commutative diagram
%\begin{equation}
%\xymatrix{
%\mathscr{M}^{(s)s} \ar[rr]^\chi \ar[dr] \ar@/_10pt/[ddr]_\pi&& \mathfrak{M}^{(s)s}\ar[dl]\ar@/^10pt/[ddl]^{\bar\pi}\\ 
%& [Z^{(s)s}/PGL(V)] \ar[d]^{\pi'}& \\
%& \mathrm{M}^{(s)s} & 
%}
%\end{equation} 
\begin{equation}
  \begin{tikzpicture}[xscale=2.8,yscale=-.7, ->, bend angle=25]
\node (A0_0) at (0,0) {$\mathfrak{SM}^{(s)s}$};
\node (A0_2) at (2,0) {$\mathfrak{M}^{(s)s}$}; 
  \node (A1_1) at (1, 0) {$\mathfrak{PM}^{(s)s}$};
  \node (A2_1) at (1, 2) {$\mathrm{M}^{(s)s}$};
  \node (A) at (2.2,1) {$.$} ;
    \path (A0_0) edge [->]node [auto] {$\scriptstyle{}$} (A1_1);
    \path (A0_2) edge [->]node [auto] {$\scriptstyle{}$} (A1_1);
    \path (A0_0) edge [bend left]node [below left] {$\scriptstyle{\pi_{\mathfrak{S}}}$} (A2_1);
    \path (A0_2) edge [bend right]node [below right] {$\scriptstyle{\pi}$} (A2_1);
    \path (A1_1) edge [->] node [auto] {$\scriptstyle{\pi_{\mathfrak{P}}}$} (A2_1);
  \end{tikzpicture}
\end{equation}
Statements as those in Theorem \ref{thm:maintheorem} hold also for $\pi$ and $\pi_{\mathfrak{P}}$, cf.\ \cite{question:vistoli2013}. Moreover, according to the proof of \cite{art:nironi2008}, Thm.~6.22-(1), the universal property stated in Proposition \ref{prop:universalproperty} also holds  for $\pi$ and $\pi_{\mathfrak{P}}$.

Let us denote by $[\mathfrak{M}^{(s)s}]$ the contravariant functor which associates with any scheme $S$ of finite type the set $[\mathfrak{M}^{(s)s}](S)$ of isomorphism classes of objects of $\mathfrak{M}^{(s)s}(S)$. The morphism $\pi$ factors through  $\mathfrak{M}^{(s)s}\to [\mathfrak{M}^{(s)s}]$. To conclude this section we show that the contravariant functors $[\mathfrak{M}^{(s)s}]$ is isomorphic to the {\em moduli} functor $\underline{\mathcal{M}}^{(s)s}$ of $\delta$-(semi)stable framed sheaves on $\mathscr{X}$, i.e., the contravariant functor
\begin{equation}
\underline{\mathcal{M}}^{(s)s}=\underline{\mathcal{M}}^{(s)s}_{\mathscr{X}/k}(\mathcal{G}, \mathcal{O}_X(1); P_0, \mathcal{F}, \delta)\colon (Sch/k)^{\circ}\to (Sets)
\end{equation}
which associates with any scheme $S$ of finite type the set of isomorphism classes of flat families of $\delta$-(semi)stable framed sheaves on $\mathscr{X}$ with Hilbert polynomial $P_0$ parameterized by $S$.

\begin{theorem}\label{thm:functors}
The functor $\underline{\mathcal{M}}^{(s)s}$ is isomorphic to $[\mathfrak{M}^{(s)s}]$. 
\end{theorem}
The proof of this Theorem requires a preliminary result. Let $G$ be an affine algebraic group which acts on an algebraic stack $\mathscr M$ (the action is given in the sense of Romagny \cite{art:romagny2005}, Def.~2.1). Let $(\mathscr{M}/G)^\ast$ be the stack introduced in \cite{art:romagny2005}, Sect.~4, whose objects are --- roughly speaking --- $G$-torsors $P$ over a base scheme $T$ with a $G$-equivariant morphism $P\to \mathscr{M}$ (cf.\ \cite{art:romagny2005}, Sect.~1). Assume that  $(\mathscr{M}/G)^\ast$ is an algebraic stack and  denote by $p\colon \mathscr{M}\to (\mathscr{M}/G)^\ast$ the morphism which associates with any object $x\colon T\to \mathscr{M}$ in $\mathscr{M}(T)$ the trivial $G$-torsor $G\times T$ over $T$ with the $G$-equivariant morphism $G\times T\xrightarrow{\mathrm{id}\times x}G\times \mathscr{M}\to\mathscr{M}$, where the last morphism is the $G$-action on $\mathscr{M}$. The notion of $G$-linearized coherent sheaf is introduced in \cite{art:romagny2005}, Example~4.3. We say that a $G$-linearized coherent sheaf $\mathcal F$ on $\mathscr M$ {\em descends to $(\mathscr{M}/G)^\ast$} if there is a coherent sheaf $\mathcal E$ on $(\mathscr{M}/G)^\ast$ such that there is an isomorphism $\mathcal F\simeq p^\ast \mathcal E$ of $G$-linearized sheaves. Assume  that the following is true:
 \begin{itemize}\setlength{\itemsep}{2pt} 
 \item every $G$-linearized coherent sheaf $\mathcal F$ on $\mathscr M$   descends to $(\mathscr{M}/G)^\ast$;
 \item  given two $G$-linearized  coherent sheaves $\mathcal F$ and $\mathcal F'$ on $\mathscr{M}$ that descend to $\mathcal E$ and $\mathcal E'$, and a morphism $f\colon \mathcal F \to \mathcal F'$, then $f$ descends to $(\mathscr{M}/G)^\ast$, i.e., there is morphism $g\colon \mathcal E \to \mathcal E'$ such that $p^\ast g$ corresponds to $f$ under the isomorphisms $\mathcal F\simeq p^\ast \mathcal E$, $\mathcal F'\simeq p^\ast \mathcal E'$;
 \end{itemize}
then we say that the category of $G$-linearized coherent sheaves on $\mathscr{M}$ {\em descends} to $(\mathscr{M}/G)^\ast$. 
 
Consider now the following particular situation: let $p\colon P\to S$ be a $\mathrm{GL}(V)$-torsor. Then there is an induced $\mathrm{GL}(V)$-action on $P\times \mathscr{X}$, where the $\mathrm{GL}(V)$-action on $\mathscr{X}$ is the trivial one. Therefore, $\left((P\times \mathscr{X})/\mathrm{GL}(V)\right)^\ast$ is $S\times \mathscr{X}$. We have the following characterization of the category of $\mathrm{GL}(V)$-linearized coherent sheaves on $P\times \mathscr{X}$.
\begin{lemma}\label{lem:descendresult}
Let $p\colon P\to S$ be a $\mathrm{GL}(V)$-torsor. The category of $\mathrm{GL}(V)$-linearized coherent sheaves on $P\times \mathscr{X}$ (where the $\mathrm{GL}(V)$-action on $\mathscr{X}$ is the trivial one) descends to $S\times \mathscr{X}$. 
\end{lemma}
\proof
Since $\mathscr{X}$ is of the form $[T/G]$, with $T$ a scheme and $G$ a linear algebraic group, by \cite{art:gholampourjiangkool2012}, Prop.~2.5, the category of $\mathrm{GL}(V)$-linearized coherent sheaves on $P\times \mathscr{X}$ is equivalent to the category of coherent sheaves on $P\times T$ with commuting $\mathrm{GL}(V)-$ and $G$-equivariant structures, where the action of $\mathrm{GL}(V)$ on $T$ and the action of $G$ on $P$ are trivial. By \cite{book:huybrechtslehn2010}, Thm.~4.2.14, this category descends to the category of coherent sheaves on $S\times T$ with $G$-equivariant structure, where the $G$-action on $S$ is trivial. As we already noticed, this category is equivalent to the category of coherent sheaves on $S\times \mathscr{X}.$
\endproof
In the following we shall freely use this Lemma without referring to it explicitly.
\proof[Proof of   Theorem \ref{thm:functors}]
We combine the arguments of \cite{book:huybrechtslehn2010}, Lemma 4.3.1, and the methods described in \cite{art:diaconescu2012}, as done in \cite{art:sheshmani2010}, Thm.~6.2. First we need to define a natural transformation
\begin{equation}
\eta \colon [\mathfrak{M}^{(s)s}]\to \underline{\mathcal{M}}^{(s)s}\ .
\end{equation}

Let $S$ be a scheme of finite type. An object in $[\mathfrak{M}^{(s)s}](S)$ is an isomorphism class $[(q\colon P\to S, \varphi\colon P\to Z^{(s)s})]$, where $q\colon P\to S$ is a $GL(V)$-torsor over $S$ and $\varphi$ is a $GL(V)$-equivariant morphism. The pullback $\mathfrak{U}_P^{(s)s}=(\mathcal{U}_P^{(s)s}, L_{\mathcal{U}_P^{(s)s}}, \phi_{\mathcal{U}_P^{(s)s}})$ via $\varphi\times\mathrm{id}_{\mathscr{X}}$ of the universal family $\mathfrak{U}^{(s)s}=(\mathcal{U}^{(s)s}, $ $ L_{\mathcal{U}^{(s)s}}, \phi_{\mathcal{U}^{(s)s}})$ parameterized by $Z^{(s)s}$ is a $GL(V)$-equivariant family on $P\times \mathscr{X}.$ Since $q\colon P\to S$ is a $GL(V)$-torsor, this family descends to a flat family of $\delta$-(semi)stable framed sheaves on $\mathscr{X}$  parameterized by $S.$
 
Let $(q\colon P\to S, \varphi\colon P\to Z^{(s)s})$ and $(q'\colon P'\to S, \varphi'\colon P'\to Z^{(s)s})$ be two isomorphic pairs, i.e., $P$ and $P'$ fit into a commutative diagram
\begin{equation}
  \begin{tikzpicture}[xscale=3,yscale=-1.3]
  \node (A0_0) at (0, 0) {$P$};
\node (A1_0) at (1, 0) {$P'$};
\node (A0_1)at (0, 1) {$S$};
\node (A1_1) at (1, 1) {$S$};
\node (A) at (1.5,.5) {$,$};
    \path (A0_0) edge [->]node [auto] {$\scriptstyle{\nu}$} (A1_0);
    \path (A0_0) edge [->]node [auto] {$\scriptstyle{q}$} (A0_1);
    \path (A1_0) edge [->]node [auto] {$\scriptstyle{q'}$} (A1_1);
    \path (A0_1) edge [->]node [auto] {$\scriptstyle{\mathrm{id}_S}$} (A1_1);
  \end{tikzpicture}
\end{equation}
where $\nu$ is an isomorphism of $\mathrm{GL}(V)$-torsors compatible with $\varphi$ and $\varphi'.$ Then $\nu$ induces an isomorphism between $\mathfrak{U}_P^{(s)s}$ and the pullback of $\mathfrak{U}_{P'}^{(s)s}$ via $\nu\times \mathrm{id}_{\mathscr{X}}.$ Thus, $\nu$ induces an isomorphism between the corresponding flat families of framed sheaves  parameterized by $S.$ Therefore, $\eta(S)$ sends isomorphism classes $[(q\colon P\to S, \varphi\colon P\to Z^{(s)s})]$ to isomorphism classes $[\mathfrak{U}_P^{(s)s}].$

We want to define a natural transformation $\gamma \colon \underline{\mathcal{M}}^{(s)s} \to [\mathfrak{M}^{(s)s}]$. Note that $p_{S\times \mathscr{X}, \mathscr{X}}^\ast \mathcal{G}$ is a generating sheaf for $S\times \mathscr{X}$ by Theorem \ref{thm:generatingsheaf}. By the same reason, for any point $s\in S$ the locally free sheaf $p_{\mathrm{Spec}(k(s))\times\mathscr{X}, \mathscr{X}}^\ast \mathcal{G}$ is a generating sheaf for $\mathrm{Spec}(k(s))\times \mathscr{X}$. 

%Let $\mathfrak{E}=(\mathcal{E},L_{\mathcal{E}},\phi_{\mathcal{E}})$ be a flat family of $\delta$-(semi)stable framed sheaves on $\mathscr{X}$  parameterized by $S.$ By Theorem \ref{thm:boundednessframed}, $\mathcal{E}_s$ is $m$-regular, hence $F_{p_{\mathrm{Spec}(k(s))\times\mathscr{X},\mathscr{X}}^\ast(\mathcal{G})}(\mathcal{E}_s)$ is $m$-regular for every closed $s\in S$. Consider the following cartesian diagram
%\begin{equation}
%  \begin{tikzpicture}[xscale=3.5,yscale=-1.3]
%  \node (A0_0) at (0, 0) {$\mathrm{Spec}(k(s))\times \mathscr{X}$};
%\node (A1_0) at (1, 0) {$S\times \mathscr{X}$};
%\node (A0_1)at (0, 1) {$\mathrm{Spec}(k(s))\times X$};
%\node (A1_1) at (1, 1) {$S\times X$};
%    \path (A0_0) edge [->]node [auto] {$\scriptstyle{s\times\mathrm{id}_{\mathscr{X}}}$} (A1_0);
%    \path (A0_0) edge [->]node [auto] {$\scriptstyle{\mathrm{id}\times \pi}$} (A0_1);
%    \path (A1_0) edge [->]node [auto] {$\scriptstyle{\mathrm{id}\times \pi}$} (A1_1);
%    \path (A0_1) edge [->]node [auto] {$\scriptstyle{s\times\mathrm{id}_{X}}$} (A1_1);
%  \end{tikzpicture}
%\end{equation}
By \cite{art:nironi2008}, Prop.~1.5, one has
\begin{align}
F_{p_{S\times\mathscr{X}, \mathscr{X}}^\ast \mathcal{G}}(\mathcal{E})_s&=(s\times\mathrm{id}_{X})^\ast(\mathrm{id}\times \pi)_\ast(\mathcal{E}\otimes p_{S\times\mathscr{X}, \mathscr{X}}^\ast(\mathcal{G})^\vee)\\
&\simeq (\mathrm{id}\times \pi)_\ast(s\times\mathrm{id}_{\mathscr{X}})^\ast(\mathcal{E}\otimes p_{S\times\mathscr{X}, \mathscr{X}}^\ast(\mathcal{G})^\vee)=F_{p_{\mathrm{Spec}(k(s))\times\mathscr{X}, \mathscr{X}}^\ast \mathcal{G}}(\mathcal{E}_s)\ .
\end{align}
Thus $F_{p_{S\times\mathscr{X}, \mathscr{X}}^\ast \mathcal{G}}(\mathcal{E})_s$ is $m$-regular for any closed $s\in S.$ Therefore 
\begin{equation}
\mathcal{B}:={p_{S\times X, S}}_\ast\left(F_{p_{S\times\mathscr{X}, \mathscr{X}}^\ast \mathcal{G}}(\mathcal{E})\otimes p_{S\times X, X}^\ast\mathcal{O}_X(m)\right)
\end{equation}
is a locally free $\mathcal{O}_S$-module of rank $P_0(m)$. Consider the frame bundle $P:=\mathrm{Isom}(V\otimes \mathcal{O}_S, \mathcal{B})\xrightarrow{q} S.$ It is a $GL(V)$-torsor over $S.$ Moreover, there is a \emph{universal} trivialization $\theta_P\colon V\otimes \mathcal{O}_P \xrightarrow{\sim} q^\ast \mathcal{B}$.

So far, starting from the flat family $\mathfrak E$ parameterized by $S$, we constructed a $\mathrm{GL}(V)$-torsor $P$ over $S$. Now we need to defined a morphism $P\to Z^{(s)s}$ induced by $\mathfrak E$. To obtain this we shall build a flat family of $\delta$-(semi)stable framed sheaves parameterized by $P$; by using the universal property of $Z^{(s)s}$ we shall obtain the morphism. Since $\mathcal{B}$ is locally free, there is a surjective morphism
\begin{equation}
p_{S\times X, S}^\ast \mathcal{B}\otimes p_{S\times X, X}^\ast\mathcal{O}_X(-m)\to F_{p_{S\times\mathscr{X}, \mathscr{X}}^\ast \mathcal{G}}(\mathcal{E})\to 0\ .
\end{equation}
By pulling backing by $q\times \mathrm{id}_X\colon P\times X\to S\times X$ we obtain
\begin{equation}
(q\times \mathrm{id}_X)^\ast(p_{S\times X, S}^\ast \mathcal{B}\otimes p_{S\times X, X}^\ast \mathcal{O}_X(-m))\to (q\times \mathrm{id}_X)^\ast F_{p_{S\times\mathscr{X}, \mathscr{X}}^\ast \mathcal{G}}(\mathcal{E})\to 0\ .
\end{equation}
On the other hand we have a morphism
\begin{equation}
(q\times \mathrm{id}_X)^\ast(p_{S\times X, S}^\ast \mathcal{B}\otimes p_{S\times X, X}^\ast \mathcal{O}_X(-m))\xrightarrow{\theta_P\otimes \mathrm{id}} V\otimes (q\times \mathrm{id}_X)^\ast p_{S\times X, X}^\ast \mathcal{O}_X(-m)\ ,
\end{equation}
so we get the quotient on $P\times X$
\begin{equation}
V\otimes p_{P\times X, X}^\ast \mathcal{O}_X(-m)\to (q\times \mathrm{id}_X)^\ast F_{p_{S\times\mathscr{X}, \mathscr{X}}^\ast \mathcal{G}}(\mathcal{E})\to 0\ .
\end{equation}
By applying~\cite{art:nironi2008}, Prop.~1.5, 
%to the cartesian diagram
%\begin{equation}
%  \begin{tikzpicture}[xscale=3.5,yscale=-1.3]
%  \node (A0_0) at (0, 0) {$P\times \mathscr{X}$};
%\node (A1_0) at (1, 0) {$S\times \mathscr{X}$};
%\node (A0_1)at (0, 1) {$P\times X$};
%\node (A1_1) at (1, 1) {$S\times X$};
%    \path (A0_0) edge [->]node [auto] {$\scriptstyle{q\times\mathrm{id}_{\mathscr{X}}}$} (A1_0);
%    \path (A0_0) edge [->]node [auto] {$\scriptstyle{\mathrm{id}\times \pi}$} (A0_1);
%    \path (A1_0) edge [->]node [auto] {$\scriptstyle{\mathrm{id}\times \pi}$} (A1_1);
%    \path (A0_1) edge [->]node [auto] {$\scriptstyle{q\times\mathrm{id}_{X}}$} (A1_1);
%  \end{tikzpicture}
%\end{equation}
we obtain
\begin{multline}
(q\times \mathrm{id}_X)^\ast F_{p_{S\times\mathscr{X}, \mathscr{X}}^\ast \mathcal{G}}(\mathcal{E}) = (q\times \mathrm{id}_X)^\ast(\mathrm{id}\times \pi)_\ast(\mathcal{E}\otimes p_{S\times\mathscr{X}, \mathscr{X}}^\ast(\mathcal{G})^\vee)\\
 \simeq  (\mathrm{id}\times \pi)_\ast(q\times \mathrm{id}_{\mathscr{X}})^\ast(\mathcal{E}\otimes p_{S\times\mathscr{X}, \mathscr{X}}^\ast(\mathcal{G})^\vee)=F_{p_{P\times\mathscr{X}, \mathscr{X}}^\ast \mathcal{G}}((q\times \mathrm{id}_{\mathscr{X}})^\ast \mathcal{E})\ .
\end{multline}
So we have the quotient on $P\times X$
\begin{equation}
V\otimes p_{P\times X, X}^\ast \mathcal{O}_X(-m)\to F_{p_{P\times\mathscr{X}, \mathscr{X}}^\ast \mathcal{G}}((q\times \mathrm{id}_{\mathscr{X}})^\ast \mathcal{E})\to 0\ .
\end{equation}
By applying the functor $G_{p_{P\times\mathscr{X}, \mathscr{X}}^\ast \mathcal{G}}$ and composing on the right with $\theta_{p_{P\times\mathscr{X}, \mathscr{X}}^\ast \mathcal{G}}((q\times \mathrm{id}_{\mathscr{X}})^\ast \mathcal{E})$, we obtain a surjective morphism
\begin{equation}
G_{p_{P\times\mathscr{X}, \mathscr{X}}^\ast \mathcal{G}}(V\otimes p_{P\times X, X}^\ast \mathcal{O}_X(-m))\to \mathcal{E}_P\ ,
\end{equation}
where $\mathcal{E}_P:=(q\times \mathrm{id}_{\mathscr{X}})^\ast \mathcal{E}$. By the universal property of $\tilde{\mathbf{Q}}$, we obtain a morphism $P\to \tilde{\mathbf{Q}}.$ We need to define a framing for $\mathcal{E}_P$; this will
give a morphism $P\to \mathbf P$.
Let us consider the morphism
\begin{equation}
\tilde{\phi}_{\mathcal{E}}\colon p_{S\times \mathscr{X}, S}^\ast L_{\mathcal{E}}\otimes \mathcal{E}\to p_{S\times \mathscr{X}, \mathscr{X}}^\ast \mathcal{F}\ .
\end{equation}
By applying the functors $\otimes p_{S\times \mathscr{X}, \mathscr{X}}^\ast(\mathcal{G})^\vee$ and  $(\mathrm{id}\times \pi)_\ast$ we get
\begin{equation}
 p_{S\times X, S}^\ast L_{\mathcal{E}}\otimes F_{p_{S\times \mathscr{X}, \mathscr{X}}^\ast \mathcal{G}}(\mathcal{E})\to (\mathrm{id}\times \pi)_\ast(p_{S\times \mathscr{X}, \mathscr{X}}^\ast \mathcal{F}\otimes p_{S\times \mathscr{X}, \mathscr{X}}^\ast(\mathcal{G})^\vee)\ ,
\end{equation}
%Here we use the fact $p_{S\times \mathscr{X}, S}=p_{S\times X, S}\circ (\mathrm{id}\times \pi)$ 
where we used the projection formula for $\pi$ (Lemma \ref{lem:projectionformula}).

By \cite{art:nironi2008}, Prop.~1.5, we have $(\mathrm{id}\times \pi)_\ast(p_{S\times \mathscr{X}, \mathscr{X}}^\ast \mathcal{F}\otimes p_{S\times \mathscr{X}, \mathscr{X}}^\ast(\mathcal{G})^\vee)\simeq p_{S\times X, X}^\ast F_{\mathcal{G}}(\mathcal{F})$. So, we get
\begin{equation}
 p_{S\times X, S}^\ast L_{\mathcal{E}}\otimes F_{p_{S\times \mathscr{X}, \mathscr{X}}^\ast \mathcal{G}}(\mathcal{E})\to p_{S\times X, X}^\ast F_{\mathcal{G}}(\mathcal{F})\ .
\end{equation}
By applying the functors $\otimes p_{S\times X, X}^\ast \mathcal{O}_X(m)$ and ${p_{S\times X, S}}_\ast$ we obtain
\begin{equation}
\hat{\phi}_{\mathcal{E}}\colon L_{\mathcal{E}}\otimes {p_{S\times X, S}}_\ast(F_{p_{S\times \mathscr{X}, \mathscr{X}}^\ast \mathcal{G}}(\mathcal{E})(m))\to \mathcal{O}_S\otimes H^0(F_{\mathcal{G}}(\mathcal{F})(m))\ .
\end{equation}

Let $L_{\mathcal{E}_P}:=q^\ast L_{\mathcal{E}}$. By considering the composition
\begin{equation}
L_{\mathcal{E}_P}\otimes V\otimes \mathcal{O}_P\xrightarrow{\mathrm{id}\otimes \theta_P} L_{\mathcal{E}_P}\otimes q^\ast{p_{S\times X, S}}_\ast (F_{p_{S\times\mathscr{X}, \mathscr{X}}^\ast \mathcal{G}}(\mathcal{E})(m))\xrightarrow{q^\ast \hat{\phi}_{\mathcal{E}}} \mathcal{O}_P\otimes H^0(F_{\mathcal{G}}(\mathcal{F})(m))
\end{equation}
we obtain the morphism
\begin{equation}
\hat{\phi}_{\mathcal{E}_P}\colon L_{\mathcal{E}_P}\to \mathcal{H}om(V\otimes \mathcal{O}_P, \mathcal{O}_P\otimes H^0(F_{\mathcal{G}}(\mathcal{F})(m)))
\simeq \operatorname{Hom}(V, H^0(F_{\mathcal{G}}(\mathcal{F})(m)))\otimes \mathcal{O}_P\ .
\end{equation}
By the universal property of the projective space $\mathbf{P}$, this defines a morphism $P\to \mathbf{P}.$ Note that $\hat{\phi}_{\mathcal{E}_P}$ is induced by the following morphism
\begin{equation}
\phi_{\mathcal{E}_P}:=q^\ast \phi_{\mathcal{E}} \colon L_{\mathcal{E}_P}\to {p_{P\times \mathscr{X}, \mathscr{X}}}_\ast\mathcal{H}om(\mathcal{E}_P, p_{P\times \mathscr{X}, \mathscr{X}}^\ast \mathcal{F})\ .
\end{equation}
Since $\mathfrak{E}=(\mathcal{E},L_{\mathcal{E}},\phi_{\mathcal{E}})$ is a flat family of $\delta$-(semi)stable framed sheaves, by Proposition \ref{prop:factorization} the morphism $P\to \tilde{\mathbf{Q}}\times \mathbf{P}$ factorizes through a morphism $P\to Z'$, hence by Proposition \ref{prop:semistableZ} also through a morphism $\Psi_{\mathfrak{E}}\colon P \to Z^{(s)s}.$ Moreover, the morphism $\Psi_{\mathfrak{E}}$ is $\mathrm{GL}(V)$-equivariant by construction. Thus $\mathfrak{E}$ defines a $\mathrm{GL}(V)$-torsor $P$ over $S$ with a $\mathrm{GL}(V)$-equivariant morphism $\Psi_{\mathfrak{E}}.$ Therefore $\mathfrak{E}$ defines an object in $\mathfrak{M}^{(s)s}(S).$

By the previous constructions, two isomorphic families $\mathfrak{E}$ and $\mathfrak{E}'$ define  two isomorphic $\mathrm{GL}(V)$-torsors 
$P$ and $P'$ over $S$, and the isomorphism 
$\nu\colon P\xrightarrow{\sim} P'$  is compatible with $\Psi_{\mathfrak{E}}$ and $\Psi_{\mathfrak{E}'}.$ So we get a natural transformation
\begin{equation}
\gamma \colon \underline{\mathcal{M}}^{(s)s} \to [\mathfrak{M}^{(s)s}] \ .
\end{equation}
Moreover, it is easy to prove that $\gamma \circ \eta$ and $\eta\circ \gamma$ are the identity.
\endproof
We have obtained the following factorization of the structure morphism $\pi$:
%\begin{equation}
%\xymatrix{
%\mathfrak{M}^{(s)s} \ar[r]^{\bar{\pi}}\ar[d]_\eta& \mathrm{M}^{(s)s} \\
%\underline{\mathcal{M}}^{(s)s}   \ar[ur]_{\Psi^{(s)s}}& 
%}
%\end{equation}
\begin{equation}\label{eq:modulispace}
\begin{aligned}
  \begin{tikzpicture}[xscale=3.5,yscale=-1.3,->, bend angle=40]
  \node (A0_0) at (0, 0) {$\mathfrak{M}^{(s)s}$};
\node (A1_0) at (1, 0) {$[\mathfrak{M}^{(s)s}]$};
\node (A2_0) at (2, 0) {$\mathrm{M}^{(s)s}$};
\node (A3_0) at (2.3, 0.165) {$.$};
\node (A0_1)at (1, 1) {$\underline{\mathcal{M}}^{(s)s}$};
    \path (A0_0) edge [bend right]node [auto] {$\scriptstyle{\pi}$} (A2_0);
        \path (A0_0) edge [->]node [auto] {$\scriptstyle{}$} (A1_0);
        \path (A1_0) edge [->]node [auto] {$\scriptstyle{}$} (A2_0);
    \path (A1_0) edge [->]node [auto] {$\scriptstyle{\eta}$} (A0_1);
        \path (A0_1) edge [->]node [below right] {$\scriptstyle{\Psi^{(s)s}}$} (A2_0);
  \end{tikzpicture}
\end{aligned}
\end{equation}
 
\subsection{The moduli spaces of $\delta$-(semi)stable framed sheaves}

In this section we prove that $\mathrm{M}^{(s)s}$ is a {\em moduli space} for the functor $\underline{\mathcal{M}}^{(s)s}$, i.e., $\mathrm{M}^{(s)s}$ corepresents $\underline{\mathcal{M}}^{(s)s}$ (cf.\ \cite{book:huybrechtslehn2010}, Def.~2.2.1). In addition, thanks to the next Proposition, we can prove that $\mathrm{M}^{s}$ is a {\em fine} moduli space for $\underline{\mathcal{M}}^{s}$, i.e., $\mathrm{M}^{s}$ represents $\underline{\mathcal{M}}^{s}$.
%\begin{lemma}\label{lem:vanishing}
%Let $T$ and $Y$ be schemes. Let $F$ be a coherent sheaf on $T\times Y.$ If for all closed points $t\in T$ the restriction $F_t:=F_{|\{t\}\times Y}$ is the zero sheaf, then $F$ is the zero sheaf as well.
%\end{lemma}
%\proof
%Let $x\in T\times Y$ be a closed point and $t=p_{T\times Y, T}(x)\in T$. By using the commutative diagram
%\begin{equation}
%  \begin{tikzpicture}[xscale=2.5,yscale=-1,->, bend angle=23]
%  \node (A0_0) at (0, 0) {$\mathrm{Spec}(k(x))$};
%\node (A1_1) at (1,1) {$Y_t$};
%\node (A2_1) at (2, 1) {$T\times Y$};
%\node (A1_2) at (1, 2) {$\mathrm{Spec}(k(t))$};
%\node (A2_2)at (2, 2) {$T$};
%    \path (A0_0) edge [bend right]node [auto] {$\scriptstyle{p}$} (A2_1);
%   \path (A0_0) edge [bend left]node [auto] {$\scriptstyle{}$} (A1_2);
%    \path (A0_0) edge [->]node [auto] {$\scriptstyle{\phi}$} (A1_1);
%    \path (A1_1) edge [->]node [auto] {$\scriptstyle{}$} (A2_1);
%    \path (A1_1) edge [->]node [auto] {$\scriptstyle{}$} (A1_2);
%    \path (A2_1) edge [->]node [auto] {$\scriptstyle{}$} (A2_2);
%    \path (A1_2) edge [->]node [auto] {$\scriptstyle{t}$} (A2_2);
%  \end{tikzpicture}
%\end{equation}
%where $Y_t$ is the fibre product,  
%we obtain $F\otimes k(x)\simeq  \phi^\ast(F_t)$. Since by hypothesis $F_t=0$, also $F\otimes k(x)=0$. By Nakayama's lemma,  the stalk of $F$ at $x$ is zero. By \cite{book:eisenbud1995}, Lem. 2.8, this holds true at all points $x\in T\times Y$. Thus $F=0$.
%\endproof 
\begin{proposition}\label{prop:Gm-invariance}
Let $\mathfrak{U}^{s}=(\mathcal{U}^{s}, L_{\mathcal{U}^{s}}, \phi_{\mathcal{U}^{s}})$ be the universal family of stable framed sheaves on $\mathscr{X}$  parameterized by $Z^s.$ Then $\mathcal{U}^{s}$ and $L_{\mathcal{U}^{s}}$ are invariant with respect to the action of the centre $\mathbb{G}_m$ of $\mathrm{GL}(V).$
\end{proposition}
\proof
By Lemma \ref{lem:automorphism} and Proposition \ref{prop:semistableZ}, the restriction $\tilde{\xi}$ of the action $\xi$  to the centre $\mathbb{G}_m$ of $\mathrm{GL}(V)$, defined in \eqref{eq:action}, is the trivial action on $Z^s.$ Therefore, we have induced isomorphisms
\begin{align}
&\tilde{\Lambda}_1 \colon p_{Z^s\times \mathbb{G}_m\times \mathscr{X}, Z^s\times \mathscr{X}}^\ast \mathcal{U}^s \xrightarrow{\sim} p_{Z^s\times \mathbb{G}_m\times \mathscr{X}, Z^s\times \mathscr{X}}^\ast \mathcal{U}^s\ ,  \\
&\tilde{\Lambda}_2 \colon p_{Z^s\times\mathbb{G}_m, Z^s}^\ast L_{\mathcal{U}^s}\xrightarrow{\sim} p_{Z^s\times\mathbb{G}_m, Z^s}^\ast L_{\mathcal{U}^s}\ ,
\end{align}
compatibly with $\phi_{\mathcal{U}^s}.$ In particular, we have the following commutative diagram
%%%%%
%\newcommand{\fe}{\kern-40mm\mbox{$p^\ast( p_{Z^s\times\mathbb{G}_m, Z^s}^\ast(L_{\mathcal{U}^s}))\otimes p_{Z^s\times \mathbb{G}_m\times \mathscr{X}, Z^s\times \mathscr{X}}^\ast(\mathcal{U}^s) $}}
%
%\newcommand{\te}{\kern-40mm\mbox{$p^\ast(p_{Z^s\times\mathbb{G}_m, Z^s}^\ast(L_{\mathcal{U}^s}))\otimes p_{Z^s\times \mathbb{G}_m\times \mathscr{X}, Z^s\times \mathscr{X}}^\ast(\mathcal{U}^s) $} }
%\newcommand{\se}{p_{Z^s\times \mathbb{G}_m\times \mathscr{X}, Z^s\times \mathscr{X}}^\ast(p_{Z^s\times \mathscr{X}, Z^s}^\ast(\mathcal{O}_{Z^s})\otimes p_{Z^s\times \mathscr{X}, \mathscr{X}}^\ast(\mathcal{F}))}
%\newcommand{\aruno}{p_{Z^s\times \mathbb{G}_m\times \mathscr{X}, Z^s\times \mathscr{X}}^\ast(\tilde{\phi}_{\mathcal{U}^s})}
%\newcommand{\ardue}{p^\ast(\tilde{\Lambda}_2)\otimes \Lambda_1}
%\begin{equation}
%\xymatrix{
%& \fe \ar[dd]^{\ardue}\ar[dl]_(.6){\aruno\ \ \ \ \ \ }\\
%\se & \\
%& \te \ar[ul]^(.6){\aruno\ \ \ \ \ \ \ \ \ }
%}
%\end{equation}
\begin{equation}
  \begin{tikzpicture}[xscale=3.2,yscale=-2.5]
\node (A0_2) at (2,0) {$p_{Z^s\times\mathbb{G}_m\times\mathscr X, Z^s}^\ast L_{\mathcal{U}^s}\otimes p_{Z^s\times \mathbb{G}_m\times \mathscr{X}, Z^s\times \mathscr{X}}^\ast \mathcal{U}^s$}; 
  \node (A1_0) at (0, 1) {$p_{Z^s\times \mathbb{G}_m\times \mathscr{X}, \mathscr{X}}^\ast \mathcal{F}$};
  \node (A2_2) at (2, 2) {$p_{Z^s\times\mathbb{G}_m\times \mathscr X, Z^s}^\ast L_{\mathcal{U}^s}\otimes p_{Z^s\times \mathbb{G}_m\times \mathscr{X}, Z^s\times \mathscr{X}}^\ast \mathcal{U}^s$};
 \node (A1_2) at (3.5,1) {$.$} ;
    \path (A0_2) edge [->]node [auto] {$\scriptstyle{p_{Z^s\times \mathbb{G}_m\times \mathscr{X}, Z^s\times \mathbb{G}_m}^\ast \tilde{\Lambda}_2\otimes \tilde{\Lambda}_1}$} (A2_2);
    \path (A0_2) edge [->]node [above left] {$\scriptstyle{p_{Z^s\times \mathbb{G}_m\times \mathscr{X}, Z^s\times \mathscr{X}}^\ast \tilde{\phi}_{\mathcal{U}^s}}$} (A1_0);
    \path (A2_2) edge [->]node [below left] {$\scriptstyle{p_{Z^s\times \mathbb{G}_m\times \mathscr{X}, Z^s\times \mathscr{X}}^\ast \tilde{\phi}_{\mathcal{U}^s}}$} (A1_0);
  \end{tikzpicture}
\end{equation}

Let $\mathcal{H}:=\mathrm{Im}\,(p_{Z^s\times \mathbb{G}_m\times \mathscr{X}, Z^s\times \mathbb{G}_m}^\ast \tilde{\Lambda}_2\otimes \tilde{\Lambda}_1-\mathrm{id})$. By Lemma \ref{lem:morphism-stability}-(iii), the restriction $\mathcal{H}_{(z,\lambda)}$ of $\mathcal{H}$ to the fibre $\{(z,\lambda)\}\times \mathscr{X}$ is zero for any closed point $(z,\lambda)\in Z^s\times \mathbb{G}_m$. So $F_{p_{\{(z,\lambda)\}\times \mathscr{X}, \mathscr{X}}^\ast \mathcal{G}}(\mathcal{H}_{(z,\lambda)})$ is zero as well by Lemma \ref{lem:support}. On the other hand, by \cite{art:nironi2008}, Prop.~1.5, we get 
\begin{equation}
F_{p_{\{(z,\lambda)\}\times \mathscr{X}, \mathscr{X}}^\ast \mathcal{G}}(\mathcal{H}_{(z,\lambda)})\simeq  \left(F_{p_{Z^s\times \mathbb{G}_m\times \mathscr{X}, \mathscr{X}}^\ast \mathcal{G}}(\mathcal{H})\right)_{(z,\lambda)}\  ,
\end{equation}
where $F_{p_{Z^s\times \mathbb{G}_m\times \mathscr{X}, \mathscr{X}}^\ast \mathcal{G}}(\mathcal{H})_{(z,\lambda)}$ is the restriction of the sheaf $F_{p_{Z^s\times \mathbb{G}_m\times \mathscr{X}, \mathscr{X}}^\ast \mathcal{G}}(\mathcal{H})$ to $\{(z,\lambda)\}\times X$. Thus $F_{p_{Z^s\times \mathbb{G}_m\times \mathscr{X}, \mathscr{X}}^\ast \mathcal{G}}(\mathcal{H})_{(z,\lambda)}$ is zero for any closed point $(z,\lambda)\in Z^s\times \mathbb{G}_m$. So $F_{p_{Z^s\times \mathbb{G}_m\times \mathscr{X}, \mathscr{X}}^\ast(\mathcal{G})}(\mathcal{H})$ is zero
%by Lemma \ref{lem:vanishing} 
and therefore by Lemma \ref{lem:support} the sheaf $\mathcal{H}$ is zero as well. 

Thus the morphism $p_{Z^s\times \mathbb{G}_m\times \mathscr{X}, Z^s\times \mathbb{G}_m}^\ast \tilde{\Lambda}_2\otimes \tilde{\Lambda}_1$ is the identity morphism 
\begin{equation}
\mathrm{id}_{p_{Z^s\times \mathbb{G}_m\times \mathscr{X}, Z^s\times \mathbb{G}_m}^\ast( p_{Z^s\times\mathbb{G}_m, Z^s}^\ast(L_{\mathcal{U}^s}))}\otimes \mathrm{id}_{p_{Z^s\times \mathbb{G}_m\times \mathscr{X}, Z^s\times \mathscr{X}}^\ast(\mathcal{U}^s)}\ ,
\end{equation}
and $\tilde{\Lambda}_1$ and $p_{Z^s\times \mathbb{G}_m\times \mathscr{X}, Z^s\times \mathbb{G}_m}^\ast \tilde{\Lambda}_2$ are the identity morphisms.

Now it remains to prove that $\tilde{\Lambda}_2$ is the identity morphism as well. Since $Z^s\times \mathbb{G}_m\times X$ is a coarse moduli scheme for $Z^s\times \mathbb{G}_m\times \mathscr{X}$, the morphism $\mathrm{id}_{Z^s\times \mathbb{G}_m}\times \pi$ is a good moduli space (cf.\ Definition \ref{def:goodmoduli}). Therefore, by applying \cite{art:alper2008}, Prop.~4.5, we get for $L:=p_{Z^s\times\mathbb{G}_m, Z^s}^\ast L_{\mathcal{U}^s}$
\begin{multline*}
{p_{Z^s\times \mathbb{G}_m\times \mathscr{X}, Z^s\times \mathbb{G}_m}}_\ast p_{Z^s\times \mathbb{G}_m\times \mathscr{X}, Z^s\times \mathbb{G}_m}^\ast L\\ \simeq  {p_{Z^s\times \mathbb{G}_m\times X, Z^s\times \mathbb{G}_m}}_\ast (\mathrm{id}_{Z^s\times \mathbb{G}_m}\times \pi)_\ast(\mathrm{id}_{Z^s\times \mathbb{G}_m}\times \pi)^\ast p_{Z^s\times \mathbb{G}_m\times X, Z^s\times \mathbb{G}_m}^\ast L\\ 
\simeq  {p_{Z^s\times \mathbb{G}_m\times X, Z^s\times \mathbb{G}_m}}_\ast p_{Z^s\times \mathbb{G}_m\times X, Z^s\times \mathbb{G}_m}^\ast L\ .
\end{multline*}
Here we used that $p_{Z^s\times \mathbb{G}_m\times \mathscr{X}, Z^s\times \mathbb{G}_m}=p_{Z^s\times \mathbb{G}_m\times X, Z^s\times \mathbb{G}_m}\circ (\mathrm{id}_{Z^s\times \mathbb{G}_m}\times \pi)$.

Since $X$ is a projective irreducible scheme, by using \cite{book:hartshorne1977}, Prop.~III-9.3, we have 
\begin{equation}
{p_{Z^s\times \mathbb{G}_m\times X, Z^s\times \mathbb{G}_m}}_\ast p_{Z^s\times \mathbb{G}_m\times X, Z^s\times \mathbb{G}_m}^\ast L\simeq  L\ .
\end{equation}
Thus $\tilde{\Lambda}_1$ and $\tilde{\Lambda}_2$ are the identity morphisms and therefore $\mathcal{U}^s$ and $L_{\mathcal{U}^s}$ are $\mathbb{G}_m$-invariant.
\endproof

\begin{theorem}
Let $\mathscr{X}$ be a $d$-dimensional normal projective irreducible stack with coarse moduli scheme $\pi\colon \mathscr{X}\to X$ and $(\mathcal{G}, \mathcal{O}_X(1))$ a polarization on it. For any framing sheaf $\mathcal{F}$, stability polynomial $\delta$ and numerical polynomial $P_0$ of degree $d$, the projective scheme $\mathrm{M}^{ss}=\mathrm{M}^{ss}_{\mathscr{X}/k}(\mathcal{G}, \mathcal{O}_X(1); P_0,$ $\mathcal{F}, \delta)$ is a moduli space for the moduli functor $\underline{\mathcal{M}}^{ss}=\underline{\mathcal{M}}^{ss}_{\mathscr{X}/k}(\mathcal{G}, \mathcal{O}_X(1); P_0, \mathcal{F}, \delta)$. Moreover the quasi-projective scheme $\mathrm{M}^s=\mathrm{M}^{s}_{\mathscr{X}/k}(\mathcal{G}, \mathcal{O}_X(1); P_0, \mathcal{F}, \delta)$ is a fine moduli space for the moduli functor $\underline{\mathcal{M}}^{s}=\underline{\mathcal{M}}^{s}_{\mathscr{X}/k}(\mathcal{G}, \mathcal{O}_X(1); P_0, \mathcal{F}, \delta)$.
\end{theorem}
\proof
Let $\Psi^{(s)s}\colon\underline{\mathcal{M}}^{(s)s}\to \mathrm{M}^{(s)s}$ be the natural transformation defined in \eqref{eq:modulispace}. Let $N$ be a scheme and $\psi\colon\underline{\mathcal{M}}^{ss}\to N$ a natural transformation. Then the universal family $\mathfrak{U}^{ss}$ of $\delta$-semistable framed sheaves on $\mathscr{X}$ parameterized by $Z^{ss}$ defines a morphism $f\colon Z^{ss}\to N$ which is $\mathrm{SL}(V)$-invariant due to the  $\mathrm{SL}(V)$-equivariance of $\mathfrak{U}^{ss}$. Since $\mathrm{M}^{ss}$ is a categorical quotient, the morphism $f$ factors through a morphism $\mathrm{M}^{ss}\to N$, therefore the natural transformation $\psi$ factors through $\Psi^{ss}.$

By Lemma \ref{lem:automorphism} the stabilizer in $\mathrm{PGL}(V)$ of a closed point in $Z^s$ is trivial. Hence, by Proposition \ref{prop:s-equivalence} and Luna's \'etale slice Theorem (\cite{book:huybrechtslehn2010}, Thm.~4.2.12), $Z^s\to M^s$ is a $\mathrm{PGL}(V)$-torsor. Since the universal family $\mathfrak{U}^s$ of  $\delta$-stable framed sheaves on $\mathscr{X}$ parameterized by $Z^s$ is $\mathrm{PGL}(V)$-linearized by Proposition \ref{prop:Gm-invariance}, it descends to a universal family of $\delta$-stable framed sheaves on $\mathscr{X}$ parameterized by $\mathrm{M}^s.$
\endproof
\begin{corollary}
The algebraic stack $\mathfrak{M}^s$ is a $\mathbb{G}_{m}$-gerbe over its coarse moduli scheme $\mathrm{M}^s$.
\end{corollary}

We conclude this section by stating a theorem about the tangent space and the obstruction to the smoothness of the moduli spaces of $\delta$-stable framed sheaves. The proof is just a straightforward generalization of the same result for $\delta$-framed sheaves on smooth projective varieties (cf.\ \cite{art:huybrechtslehn1995-II}, Thm.~4.1), thanks to the result of Olsson and Starr about the tangent space of the Quot scheme for Deligne-Mumford stacks (cf.\ \cite{art:olssonstarr2003}, Lemma 2.5).
\begin{theorem}\label{thm:deformations}
Let $[(\mathcal{E}, \phi_\mathcal{E})]$ be a point in the moduli space $\mathrm{M}^s_{\mathscr{X}/k}(\mathcal{G}, \mathcal{O}_X(1); P_0, \mathcal{F}, \delta)$ of $\delta$-stable framed sheaves on $\mathscr{X}$. Consider $\mathcal{E}$ and $\phi_\mathcal{E}\colon\mathcal{E}\to\mathcal{F}$ as complexes  concentrated in degree zero, and  zero and one, respectively. 
\begin{itemize}\setlength{\itemsep}{2pt}
\item[(i)] The Zariski tangent space of $\mathrm{M}^s_{\mathscr{X}/k}(\mathcal{G}, \mathcal{O}_X(1); P_0, \mathcal{F}, \delta)$ at a point $[(\mathcal{E}, \phi_\mathcal{E})]$ is naturally isomorphic to the first hyper-Ext group $\mathbb{E}\mathrm{xt}^1(\mathcal{E}, \mathcal{E}\xrightarrow{\phi_\mathcal{E}}\mathcal{F})$.
\item[(ii)] If the second hyper-Ext group $\mathbb{E}\mathrm{xt}^2(\mathcal{E}, \mathcal{E}\xrightarrow{\phi_\mathcal{E}}\mathcal{F})$ vanishes, $\mathrm{M}^s_{\mathscr{X}/k}(\mathcal{G}, \mathcal{O}_X(1); P_0, \mathcal{F}, \delta)$ is smooth at $[(\mathcal{E}, \phi_\mathcal{E})]$.
\end{itemize}
\end{theorem} 
 
\subsection{Framed sheaves on projective orbifolds}

In this section we assume that $\mathscr{X}$ is a projective orbifold of dimension $d$. Let $\mathfrak{E}=(\mathcal{E},\phi_{\mathcal{E}})$ be a $d$-dimensional framed sheaf on $\mathscr{X}.$  The \emph{orbifold rank} (resp.\ the \emph{degree}) of $\mathfrak{E}$ is the orbifold rank (resp.\ the degree) of $\mathcal{E}.$ The \emph{framed degree} of a $d$-dimensional framed sheaf $\mathfrak{E}$ is
\begin{equation}
\deg_{\mathcal{G}}(\mathfrak{E}):=\deg_{\mathcal{G}}(\mathcal{E})-\epsilon(\phi_{\mathcal{E}})\delta_1\ ,
\end{equation}
while its \emph{framed slope} is
\begin{equation}
\mu_{\mathcal{G}}(\mathfrak{E}):=\frac{\deg_{\mathcal{G}}(\mathfrak{E})}{\mathrm{ork}(\mathcal{E})}\ .
\end{equation}
Let $\mathcal{E}'$ be a subsheaf of $\mathcal{E}$ with quotient $\mathcal{E}''=\mathcal{E}/\mathcal{E}'.$ If $\mathcal{E}$, $\mathcal{E}'$ and $\mathcal{E}''$ are $d$-dimensional, the framed degree of $\mathfrak{E}$ behaves additively, i.e., $\deg_{\mathcal{G}}(\mathfrak{E})=\deg_{\mathcal{G}}(\mathfrak{E}')+\deg_{\mathcal{G}}(\mathfrak{E}'')$.

\begin{definition}\label{def:stab}
A $d$-dimensional framed sheaf $\mathfrak{E}=(\mathcal{E},\phi_{\mathcal{E}})$ is {\em $\mu$-(semi)stable with respect to $\delta_1$} if and only if $\ker\phi_{\mathcal{E}}$ is torsion-free and the following conditions are satisfied:
\begin{itemize}\setlength{\itemsep}{2pt}
\item[(i)] $\deg_{\mathcal{G}}(\mathcal{E}')\; (\leq)\; \mathrm{ork}(\mathcal{E}')\mu_{\mathcal{G}}(\mathfrak{E})$ for all subsheaves $\mathcal{E}'\subseteq \ker\phi_{\mathcal{E}}$,
\item[(ii)] $(\deg_{\mathcal{G}}(\mathcal{E}')-\delta_1)\; (\leq)\; \mathrm{ork}(\mathcal{E}')\mu_{\mathcal{G}}(\mathfrak{E})$ for all subsheaves $\mathcal{E}'\subset \mathcal{E}$ with $\mathrm{ork}(\mathcal{E}')< \mathrm{ork}(\mathcal{E})$.
\end{itemize}
\end{definition}
\begin{definition}
Let $\mathfrak{E}=(\mathcal{E},\phi_{\mathcal{E}})$ be a framed sheaf of orbifold rank zero. If $\phi_{\mathcal{E}}$ is injective, we say that $\mathfrak{E}$ is {\em $\mu$-semistable} (indeed, in this case the  $\mu$-semistability of the corresponding framed sheaf does not depend on $\delta_1$). Moreover, if the degree of $\mathcal{E}$ is $\delta_1$, we say that $\mathfrak{E}$ is {\em $\mu$-stable with respect to $\delta_1$}.
\end{definition} 
Recall that the $\hat{\mu}$-(semi)stability is an open property in flat families of framed sheaves (cf.\ Proposition \ref{prop:openesshat-mu}). Moreover, on a projective orbifold the definition of $\hat{\mu}$-(semi)stability is equivalent to the definition of $\mu$-(semi)stability, thanks to the relations between the rank and the orbifold rank and between the hat-slope and the slope (cf.\ Remark \ref{rem:orbirank}). Thus we can apply the results of the previous sections. In particular we get the following result.
\begin{theorem}
Let $\mathscr{X}$ be a $d$-dimensional projective irreducible orbifold with coarse moduli scheme $\pi\colon \mathscr{X}\to X$ and $(\mathcal{G}, \mathcal{O}_X(1))$ a polarization on it. For any framing sheaf $\mathcal{F}$, stability polynomial $\delta$ and numerical polynomial $P_0$ of degree $d$, there exists a fine moduli space parameterizing isomorphism classes of $\mu$-stable framed sheaves on $\mathscr{X}$ with Hilbert polynomial $P_0$, which is a quasi-projective scheme.
\end{theorem}

\bigskip\section{$(\mathscr{D}, \mathcal{F}_{\mathscr{D}})$-framed sheaves on two-dimensional projective orbifolds}\label{sec:framedsheaves}

In this section we introduce the theory of $(\mathscr{D}, \mathcal{F}_{\mathscr{D}})$-framed sheaves on two-dimensional smooth projective irreducible stacks. Our main reference for the corresponding theory in the case of smooth projective irreducible surfaces is \cite{art:bruzzomarkushevich2011}. 

Let $\mathscr X$ be a two-dimensional smooth projective irreducible stack with coarse moduli scheme $\pi\colon \mathscr{X}\to X$ a normal projective surface. By \cite{art:vistoli1989}, Prop.~2.8, $X$ only has finite quotient (hence rational, cf.\ \cite{art:kovacs2000}) singularities.

Fix a one-dimensional integral closed substack $\mathscr{D}\subset \mathscr{X}$ and a locally free sheaf $\mathcal{F}_{\mathscr{D}}$ on it. We  call $\mathscr D$ a \emph{framing divisor} and $\mathcal{F}_{\mathscr{D}}$ a \emph{framing sheaf}.
\begin{definition}
A $(\mathscr{D}, \mathcal{F}_{\mathscr{D}})$-framed sheaf on $\mathscr{X}$ is a pair $(\mathcal{E}, \phi_{\mathcal{E}})$, where $\mathcal{E}$ is a torsion-free sheaf on $\mathscr{X}$, locally free in a neighbourhood of $\mathscr{D}$, and $\phi_{\mathcal{E}}\colon \mathcal{E}_{|\mathscr{D}}\xrightarrow{\sim} \mathcal{F}_{\mathscr{D}}$ is an isomorphism. We  call $\phi_{\mathcal{E}}$ a \emph{framing} on $\mathcal{E}.$
\end{definition}
A morphism between $(\mathscr{D}, \mathcal{F}_{\mathscr{D}})$-framed sheaves on $\mathscr{X}$ is a morphism between framed sheaves as stated in Definition \ref{def:morphism}.

The assumption of locally freeness of the underlying coherent sheaf $\mathcal{E}$ of a $(\mathscr{D}, \mathcal{F}_{\mathscr{D}})$-framed sheaf $(\mathcal{E},\phi_{\mathcal{E}})$ in a neighbourhood of $\mathscr{D}$ allows one to prove the next Lemma, which will be useful later on.

\begin{lemma}\label{lem:locallyfree}
Let $\mathcal{E}$ be a torsion-free sheaf on $\mathscr{X}$ which is locally free in a neighbourhood of $\mathscr{D}$. If $\mathcal{E}'$ is a saturated coherent subsheaf of $\mathcal{E}$, the restriction $\mathcal{E}'_{|\mathscr{D}}$ is a subsheaf of $\mathcal{E}_{|\mathscr{D}}$.
\end{lemma}
\proof
By hypothesis $\mathcal{E}$ is locally free in a neighbourhood of $\mathscr{D}$, hence on this neighbourhood $\mathcal{E}'$ is a saturated subsheaf of a locally free sheaf, hence $\mathcal{E}'$ is locally free as well (cf.\ Example \ref{ex:saturatedloc}). Thus we get straightforward the assertion.
\endproof

\subsection{Boundedness}

The first result we prove concerns the boundedness of the family of $(\mathscr{D}, \mathcal{F}_{\mathscr{D}})$-framed sheaves on $\mathscr{X}$ with fixed Hilbert polynomial. In order to prove it, we need to impose some conditions.

As explained for instance in \cite{art:toen1999}, Section 2.2., the structure morphism $\pi\colon \mathscr{X}\to X$ induces a one-to-one correspondence between integral closed substacks of $ \mathscr{X}$ and integral closed subschemes of $X$ in the following way:
for any integral closed substack $\mathscr{V}$ of $\mathscr{X}$, $\pi(\mathscr{V})$ is a closed integral subscheme of $X$, and {\em vice versa}, for any integral closed subscheme $V\subset X$, the fibre product $(V\times_X \mathscr{X})_{red}$ is an integral closed substack of $\mathscr{X}$. Moreover, $V$ is the coarse moduli scheme of $\mathscr{V}$ (cf.\ \cite{art:abramovichvistoli2002}, Lemma 2.3.3).

Let $D:=\pi(\mathscr{D})$ be the coarse moduli scheme of $\mathscr{D}$. In the following we assume that {\em  $D$ is a smooth curve}. Furthermore,  we fix a polarization $(\mathcal{G}, \mathcal{O}_X(1))$ on $\mathscr{X}$ such that {\em $\mathcal{G}$ is a direct sum of powers of a $\pi$-ample locally free sheaf}.
 
Note that the maximum $N_{\mathscr{D}}$ of the numbers of the conjugacy classes of any geometric stabilizer group of $\mathscr{D}$ is less or equal to the corresponding number $N_{\mathscr{X}}$ for $\mathscr{X}$, so that $\mathcal{G}_{|\mathscr{D}}$ is a generating sheaf for $\mathscr{D}$ (cf.\ Remark \ref{rem:generatingsheaf}). Thus, also using part (ii) of Proposition \ref{prop:kresch}, we obtain that $\mathscr D$ is a projective stack. 

Our strategy consists in proving that the family $\mathscr{C}_{\mathscr{X}}$ of torsion-free sheaves on $\mathscr{X}$ whose restriction to $\mathscr{D}$ is isomorphic to a fixed locally free sheaf is contained in the family $\mathscr{C}_{X}$ of torsion-free sheaves on $X$ whose restriction to $D$ is isomorphic to a fixed locally free sheaf. Then, using Proposition \ref{prop:boundedness}, the boundedness of the family $\mathscr{C}_{X}$ ensures the boundedness of the family $\mathscr{C}_{\mathscr{X}}$.
\begin{lemma}
Let $\mathscr X$ be a two-dimensional smooth projective irreducible stack with coarse moduli scheme $\pi\colon \mathscr{X}\to X$ a normal projective surface and $(\mathcal{G},\mathcal{O}_X(1))$ a polarization on it, where $\mathcal{G}$ is a direct sum of powers of a $\pi$-ample locally free sheaf. Fix a one-dimensional integral closed substack $\mathscr{D}\subset \mathscr{X}$, whose coarse moduli space $\mathscr D\to D$ is a smooth curve, and a locally free sheaf $\mathcal{F}_{\mathscr{D}}$ on it. Let $\mathcal{E}$ be a torsion-free sheaf on $\mathscr{X}$ such that $\mathcal{E}_{|\mathscr{D}} \simeq  \mathcal{F}_{\mathscr{D}}$. Then $F_{\mathcal{G}}(\mathcal{E})$ is a torsion-free sheaf on $X$ and $F_{\mathcal{G}}(\mathcal{E})_{| D} \simeq F_{\mathcal{G}_{|\mathscr{D}}}(\mathcal{F}_{\mathscr{D}})$ is an isomorphism, where $F_{\mathcal{G}_{\vert \mathscr{D}}}(\mathcal{F}_{\mathscr{D}})$ is a locally free sheaf on $D$.
\end{lemma}
\proof
Let us consider the short exact sequence
\begin{equation}
0 \to \mathcal{E}\otimes\mathcal{O}_{\mathscr{X}}(-\mathscr{D})\to \mathcal{E} \to i_\ast(\mathcal{F}_{\mathscr{D}}) \to 0\ .
\end{equation}
Since the functor $F_{\mathcal{G}}$ is exact, we get
\begin{equation}
0 \to F_{\mathcal{G}}(\mathcal{E}\otimes\mathcal{O}_{\mathscr{X}}(-\mathscr{D}))\to F_{\mathcal{G}}(\mathcal{E}) \to \iota_\ast(F_{\mathcal{G}_{|\mathscr{D}}}(\mathcal{F}_{\mathscr{D}})) \to 0\ ,
\end{equation}
where $i\colon\mathscr{D}\hookrightarrow \mathscr{X}$ and $\iota\colon D\hookrightarrow X$ are the inclusion morphisms.

By Proposition \ref{prop:pureness}, $F_{\mathcal{G}}(\mathcal{E})$ (resp.\ $F_{\mathcal{G}_{|\mathscr{D}}}(\mathcal{F}_{\mathscr{D}})$) is a torsion-free sheaf on $X$ (resp.\ $D$). Since $D$ is a smooth irreducible projective curve, $F_{\mathcal{G}_{|\mathscr{D}}}(\mathcal{F}_{\mathscr{D}})$ is locally free. Now $\mathrm{Supp}(\mathcal{E}\otimes\mathcal{O}_{\mathscr{X}}(-\mathscr{D}))$ is disjoint from $\mathscr{D}$, so that, by Corollary \ref{cor:supportpure}  the support of $F_{\mathcal{G}}(\mathcal{E}\otimes\mathcal{O}_{\mathscr{X}}(-\mathscr{D}))$ is disjoint from $D$ as well. Then $F_{\mathcal{G}}(\mathcal{E})_{| D} \simeq F_{\mathcal{G}_{|\mathscr{D}}}(\mathcal{F}_{\mathscr{D}})$.
\endproof

\begin{definition}\label{def:goodframing} 
An effective irreducible $\mathbb Q$-Cartier divisor $D$ in $X$ is a \emph{good framing divisor}  if there exists $a_D\in\mathbb{N}_{>0}$ such that $a_D D$ is a big and nef Cartier divisor on $X$ (i.e., $a_D D$ is a nef Cartier divisor, and $(a_D D)^2>0$).
\end{definition} 

\begin{theorem}\label{thm:boundedness-normal}
Let $X$ be a normal irreducible projective surface with rational singularities and $H$ an effective ample divisor on it. Let $D$ be a good framing divisor in $X$ which contains the singular locus of $X$, and $F_D$ a locally free sheaf on $D.$ Then for any numerical polynomial $P\in \mathbb{Q}[n]$ of degree two, the family $\mathscr{C}$ of torsion-free sheaves $E$ on $X$ such that $E_{|D} \simeq  F_{D}$ and $P(E)=P$ is bounded.
\end{theorem}
\proof
We shall adapt the arguments of \cite{phd:lehn1993}, Thm.~3.2.4. We  want to apply Kleiman's criterion (\cite{book:huybrechtslehn2010}, Thm.~1.7.8), so that we need to determine upper bounds for the quantities $h^0(X, E)$ and $h^0(H, E\vert_H)$, for any torsion-free sheaf $E$ in the family $\mathscr{C}$.

Let us fix a torsion-free sheaf $E$ on $X$ such that $E_{| D} \simeq  F_D$ and $P(E)=P.$ Consider the short exact sequence
\begin{equation}
0\to E(-(\nu+1)D)\to E(-\nu D)\to (E (-\nu D))_{| D}\to 0\ .
\end{equation}
By induction, we get $h^0(X, E)\leq h^0(X, E(-n D))+\sum_{\nu=0}^{n-1} h^0(D, F_D(-\nu D))$ for all $n\ge 1$. Let $n=ma_D+t$ with $0\le t\le a_D-1$, then by \cite{book:lazarsfeld2004}, Thm.~1.4.37, we have $h^0(X, E(-nD))=O(m^2).$ Since $\mathcal{O}_X(a_D D)$ is  big and nef, the line bundle $\mathcal{O}_X(a_D D)_{| D}$ is ample, hence there exists a positive integer $\nu_0$ such that for any $\nu\geq \nu_0$ one has $h^0(D, F_D (-\nu D)_{| D})=0$. Set  $K=\sum_{\nu=0}^{\nu_0-1} h^0(D, F_D (-\nu D)).$ This is independent of $E$ and 
\begin{equation}\label{eq:ine1}
h^0(X, E)\leq K\ .
\end{equation}

We want to estimate $h^0(H, E\vert_H)$. Since $h^0(H, E\vert_H)\leq h^0(X, E)+h^1(X, E(-H))$, we need only to estimate the quantities on the right-hand-side. First, note that $h^1(X, E (-H))= h^0(X, E (-H))+h^2(X, E (-H))-\chi(X, E (-H)).$ Since the Hilbert polynomial of $E$ is fixed, $\chi(X, E (-H))=P(-1).$ Moreover, the restriction of $E (-H)$ to $D$ is the fixed locally free sheaf $F_D\otimes \mathcal{O}_X(-H)_{| D}$, so we can adapt the previous arguments to obtain 
\begin{equation}
h^0(X, E (-H)))\leq L\ ,
\end{equation}
for some positive integer $L.$  Now we just need an estimate of $h^2(X, E (-H))$. Set $G=E (-H)$. By Serre duality (\cite{book:hartshorne1977}, Thm.~III-7.6), $H^2(X, G)\simeq \operatorname{Hom}(G,\omega_X)^\vee$, where $\omega_X$ is the dualizing sheaf of $X$.  Let $\pi\colon\hat{X}\to X$ be a resolution of singularities of $X$. Then we have the map
\begin{equation}
\operatorname{Hom}(G,\omega_X)\xrightarrow{\pi^\ast} \operatorname{Hom}(\pi^\ast G ,\pi^\ast \omega_X)\ .
\end{equation}
This map is injective. Indeed let $\varphi\colon G\to\omega_X$ be a morphism such that $\pi^\ast \varphi=0$. Since $\pi$ is an isomorphism over $X_{sm}$, the sheaf $\mathrm{Im}\,(\varphi)$ is supported on the singular locus $\mathrm{sing}(X)$. Since $\omega_X$ is a torsion free sheaf of rank one (cf.\ e.g. \cite{art:reid1980}, Appendix~1), $\varphi=0$.

By Kempf's Theorem (\cite{book:kkfmsd1973}, Ch.~I.3) we have $\pi_\ast \omega_{\hat{X}}\simeq \omega_X$, hence   a morphism $\pi^\ast\omega_X\to \omega_{\hat{X}}$, and maps 
\begin{equation}
\operatorname{Hom}(G,\omega_X)\to
\operatorname{Hom}(\pi^\ast G,\pi^\ast\omega_X) \to \operatorname{Hom}(\pi^\ast G,\omega_{\hat X})\ .
\end{equation}
The kernel of the composition $f\colon\operatorname{Hom}(G,\omega_X) \to \operatorname{Hom}(\pi^\ast G,\omega_{\hat X})$ lies in $\operatorname{Hom}(\pi^\ast G, T)$, where $T$ is the torsion of $\pi^\ast\omega_X$. Since the singularities of $X$ are in $D$, the group $\operatorname{Hom}(\pi^\ast G, T)$ injects into $\operatorname{Hom}(\pi^\ast G_{\vert \hat D}, T_{\vert \hat D})$, where $\hat D = \pi^{-1}(D)$. The dimension $M$ of the latter group does not depend on $E$ since $\pi^\ast G_{\vert \hat D}\simeq \pi^\ast (F_D\otimes \mathcal{O}_X(-H)_{\vert D})$. Thus $\dim\ker f\leq \dim \operatorname{Hom}(\pi^\ast G, T) \le M$ {for some positive integer $M$}.

Note that $\pi^\ast G$ is torsion-free since $G$ is locally free in a neighbourhood of $D$, and $D$ contains the singular locus of $X$. Consider the exact sequence
\begin{equation}
0 \to \pi^\ast G\to (\pi^\ast G)^{\vee\vee}\to Q\to 0\ ,
\end{equation}
where the support of $Q$ is zero-dimensional. 

By applying the functor $\operatorname{Hom}(\cdot,\omega_{\hat{X}})$ one gets $\operatorname{Hom}(\pi^\ast G ,\omega_{\hat{X}})\simeq\operatorname{Hom}((\pi^\ast G)^{\vee\vee},\omega_{\hat{X}})$.  The dual $(\pi^\ast G)^\vee$ is locally free, so that
\begin{equation}
\dim \operatorname{Hom}((\pi^\ast G)^{\vee\vee},\omega_{\hat{X}})=\dim \operatorname{Hom}(\mathcal{O}_{\hat{X}},(\pi^\ast G)^{\vee}\otimes\omega_{\hat{X}})=h^0((\pi^\ast G)^{\vee}\otimes\omega_{\hat{X}})\ .
\end{equation}
Moreover, $\hat{D}$ is a good framing divisor, since it is a pullback by a birational morphism, and the sheaf $\left(\pi^\ast(G)^{\vee}\otimes\omega_{\hat{X}}\right)_{|\hat{D}}$ is a fixed locally free sheaf on $\hat{D}$, so that we can use the same argument as before, and obtain
\begin{equation}
\dim \operatorname{Hom}((\pi^\ast G)^{\vee\vee},\omega_{\hat{X}}) \le N
\end{equation}
for some constant $N>0$. Then the dimension of the image of $f$ is bounded by $N$. Therefore, $h^2(X,E(-H))=\dim \operatorname{Hom}(E(-H),\omega_X)\le M+N$. Thus 
\begin{equation}%\label{eq:ine2}
h^0(H, E\vert_H)\leq h^0(X, E)+ h^1(X, E (-H))\leq K+L+M+N-P(-1)=:K'.
\end{equation}
Thus by Kleiman's criterion, the family $\mathscr{C}$ is bounded.
\endproof
\begin{theorem}\label{thm:boundedness}
Let $\mathscr X$ be a two-dimensional smooth projective irreducible stack with coarse moduli scheme $\pi\colon \mathscr{X}\to X$ a normal projective surface and $(\mathcal{G},\mathcal{O}_X(1))$ a polarization on it, where $\mathcal{G}$ is a direct sum of powers of a $\pi$-ample locally free sheaf. Fix a one-dimensional integral closed substack $\mathscr{D}\subset \mathscr{X}$ and a locally free sheaf $\mathcal{F}_{\mathscr{D}}$ on it. Assume that the coarse moduli space $\mathscr D\to D$ of $\mathscr{D}$ is a smooth curve containing the singular locus of $X$ and is a good framing divisor. For any numerical polynomial $P\in \mathbb{Q}[n]$ of degree two, the family $\mathscr{C}_{\mathscr{X}}$ of torsion-free sheaves $\mathcal{E}$ on $\mathscr{X}$ such that $\mathcal{E}_{|\mathscr{D}} \simeq  \mathcal{F}_{\mathscr{D}}$ and $P_{\mathcal{G}}(\mathcal{E})=P$ is bounded.
\end{theorem}
\proof
By Proposition \ref{prop:boundedness}, $\mathscr{C}_{\mathscr{X}}$ is bounded if and only if the family of torsion-free sheaves $F_{\mathcal{G}}(\mathcal{E})$ on $X$ such that $F_{\mathcal{G}}(\mathcal{E})_{| D} \simeq  F_{\mathcal{G}_{|\mathscr{D}}}(\mathcal{F}_{\mathscr{D}})$ and $P(F_{\mathcal{G}}(\mathcal{E}))=P$ is bounded. This is a subfamily of the family $\mathscr{C}_X$ of torsion-free sheaves $E$ on $X$ with Hilbert polynomial $P$ such that $E_{| D} \simeq  F_{\mathcal{G}_{|\mathscr{D}}}(\mathcal{F}_{\mathscr{D}}).$ This latter family is bounded by Theorem \ref{thm:boundedness-normal}.
\endproof

\subsection{Stability of $(\mathscr{D}, \mathcal{F}_{\mathscr{D}})$-framed sheaves}\label{sec:stability}

In this section we shall show that any $(\mathscr{D}, \mathcal{F}_{\mathscr{D}})$-framed sheaf on $\mathscr{X}$ is $\mu$-stable (according to Definition \ref{def:stab}) with respect to a suitable choice of an effective ample divisor on $X$ and of the parameter $\delta_1.$ From now on, we assume that $\mathscr{X}$ is an orbifold and $\mathscr{D}$ is smooth. As in the previous section, we assume that the coarse moduli space $\mathscr D\to D$ of $\mathscr{D}$ is a smooth curve and a good framing divisor.

Since $D\cap X_{sm}$ is an irreducible effective Cartier divisor, where $X_{sm}$ is the smooth locus of $X$, there exists a unique positive integer $a_{\mathscr{D}}$ such that $\pi^{-1}(D\cap X_{sm})=a_{\mathscr{D}}(\mathscr{D}\cap \pi^{-1}(X_{sm}))$. Then  $\pi^{-1}(a_D D\cap X_{sm})=a_{\mathscr{D}}a_D(\mathscr{D}\cap \pi^{-1}(X_{sm}))$. Since $X$ is normal, $\mathrm{codim}(X\setminus X_{sm})\geq 2$; moreover, since $\pi$ is a {\em codimension preserving morphism} (cf.\ \cite{art:fantechimannnironi2010}, Rem.~4.3), also $\mathrm{codim}(\pi^{-1}(X\setminus X_{sm}))$ is at least two and therefore
\begin{equation}\label{eq:equilinebundles}
\pi^\ast \mathcal{O}_X(a_D D)\simeq \mathcal{O}_{\mathscr{X}}(\mathscr{D})^{a_{\mathscr{D}}a_D}\ .
\end{equation}
This isomorphism will be useful later on.

%In the following we shall deal with the integral and rational Chow groups of $\mathscr{X}$: a well-written survey on intersection theory for Deligne-Mumford stacks is \cite{art:abramovichgrabervistoli2008}, Sect.~2.1.
%
%By \cite{art:vistoli1989}, Prop.~6.1, since $\mathscr{X}$ is smooth, its rational Chow groups are isomorphic to the rational Chow groups of $X$ via $\pi_\ast$. Moreover, 
%\begin{equation}\label{eq:pushforward}
%\pi_\ast([\mathscr{D}])=\frac{1}{k_{\mathscr{D}}}[D]\in A^1(X)_{\mathbb{Q}}\ ,
%\end{equation}
%where $k_{\mathscr{D}}$ is the order of the stabilizer of a generic geometric point of $\mathscr{D}$. The Edidin-Graham intersection product on the integral Chow groups of $\mathscr{X}$ is compatible with those of Vistoli and Gillet after extending to the rational Chow groups (cf.\ \cite{art:edidingraham1998}, Prop.~15). Hence
%\begin{equation}\label{eq:D2}
%\int_{\mathscr{X}}^{et} \mathrm{c}_1^{et}(\mathcal{O}_{\mathscr{X}}(\mathscr{D}))\cdot \mathrm{c}_1^{et}(\mathcal{O}_{\mathscr{X}}(\mathscr{D}))=\frac{(a_D D)^2}{a_D^2 k_{\mathscr D}^2}\ .
%\end{equation}
%Moreover, for any ample divisor $H$ on $X$, we have
%\begin{equation}\label{eq:DC}
%\int_{\mathscr{X}}^{et} \mathrm{c}_1^{et}(\mathcal{O}_{\mathscr{X}}(\mathscr{D}))\cdot\mathrm{c}_1^{et}(\pi^\ast(\mathcal{O}_X(H)))=\frac{(a_D D)\cdot H}{a_D k}\ .
%\end{equation}
\begin{definition}\label{def:goodframingsheaf}
Let $\mathscr{X}$ be an orbifold and $\mathscr{D}$ a smooth integral closed substack of $\mathscr{X}$ such that $D:=\pi(\mathscr{D})$ is a good framing divisor. A \emph{good framing sheaf} on $\mathscr{D}$ is a locally free sheaf $\mathcal{F}_{\mathscr{D}}$ for which there exists a real positive number $A_0$, with
\begin{equation}
0\leq A_0<\frac{1}{\mathrm{ork}(\mathcal{F}_{\mathscr{D}})}\int_{\mathscr{X}}^{et} \mathrm{c}_1^{et}(\mathcal{O}_{\mathscr{X}}(\mathscr{D}))\cdot \mathrm{c}_1^{et}(\mathcal{O}_{\mathscr{X}}(\mathscr{D}))=\frac{1}{\mathrm{ork}(\mathcal{F}_{\mathscr{D}})}\frac{(a_D D)^2}{a_D^2 k_{\mathscr D}^2}\ ,
\end{equation}
{where $k_{\mathscr{D}}$ is the order of the stabilizer of a generic geometric point of $\mathscr{D}$,} such that for any locally free subsheaf $\mathcal{F}'$ of $\mathcal{F}_{\mathscr{D}}$ we have 
\begin{equation}
\frac{1}{\mathrm{ork}(\mathcal{F}')}\int_{\mathscr{D}}^{et} \mathrm{c}_1^{et}(\mathcal{F}')\leq \frac{1}{\mathrm{ork}(\mathcal{F}_{\mathscr{D}})}\int_{\mathscr{D}}^{et} \mathrm{c}_1^{et}(\mathcal{F}_{\mathscr{D}})+A_0\ .
\end{equation}
\end{definition}
\begin{remark}
Note that if $\mathcal{F}_{\mathscr{D}}$ is a line bundle on $\mathscr{D}$, trivially it is a good framing sheaf. Moreover   a direct sum of line bundles $\mathcal L_i$ such that the value of
\begin{equation}
\int_{\mathscr{D}}^{et} \mathrm{c}_1^{et}(\mathcal L_i)
\end{equation}
is the same for all $i$ is a good framing sheaf.
\end{remark}

Let $H$ be an ample divisor on $X$; then $H_n=H+na_D D$ is ample for any positive integer $n$. Let $\mathcal{G}$ be a generating sheaf on $\mathscr{X}$. In the following we would like to compare the degree of a coherent sheaf $\mathcal{E}$ on $\mathscr{X}$ with respect to the polarizations $(\mathcal{G}, \mathcal{O}_X(H))$ and $(\mathcal{G}, \mathcal{O}_X(H_n))$. To avoid confusion, we shall write explicitly what polarization we use to compute the coefficients of the Hilbert polynomial.
 
\begin{remark}\label{rem:toenriemannroch}
Since $\mathscr{X}$ is smooth, to compute the degree of a coherent sheaf on $\mathscr{X}$ we can use the T\"oen-Riemann-Roch theorem \cite{phd:toen1999,art:toen1999}. Let $\mathcal{E}$ be a coherent sheaf on $\mathscr{X}$. Then the theorem states that the Euler characteristic $\chi(\mathcal{E})$ of $\mathcal{E}$ is equal to
\begin{equation}
\chi(\mathcal{E})=\int_{\mathscr{X}}^{rep} \mathrm{ch}^{rep}(\mathcal{E})\mathrm{td}^{rep}(\mathscr{X})\ .
\end{equation} 
The computation on the right-hand side is done not over $\mathscr{X}$ but over its inertia stack $\mathcal{I}(\mathscr{X})$: indeed, $\int_{\mathscr{X}}^{rep}$ denotes the pushforward $p_\ast\colon H^\bullet _{rep}(\mathscr{X}):=H^\bullet _{et}(\mathcal{I}(\mathscr{X}))\to H^\bullet _{et}(\mathrm{Spec}(k))\simeq \mathbb{Q}$ of the morphism $p\colon \mathcal{I}(\mathscr{X})\to \mathrm{Spec}(k)$. As explained in \cite{art:borne2007}, Sect.~3, there is a decomposition
\begin{equation}\label{eq:decomposition}
H^\bullet _{rep}(\mathscr{X})=H^\bullet _{et}(\mathcal{I}(\mathscr{X}))\cong H^\bullet _{et}(\mathscr{X})\oplus H^\bullet _{et}(\mathcal{I}(\mathscr{X})\setminus \mathscr{X})\ .
\end{equation}
For any class $\alpha\in H^\bullet _{rep}(\mathscr{X})$ we denote by $\alpha=\alpha_1+\alpha_{\neq 1}$ the corresponding decomposition. $\mathrm{ch}^{rep}(\mathcal{E})$ is the \'etale Chern character of the decomposition  of the K-theory class $\sigma^\ast[\mathcal{E}]$ with respect to the automorphism groups of the points of $\mathscr{X}$ (cf.\ \cite{art:borne2007}, Sect.~C.3.2), where $\sigma\colon \mathcal{I}(\mathscr{X})\to \mathscr{X}$ is the forgetful morphism. In particular, $[\mathrm{ch}^{rep}(\mathcal{E})]_1=\mathrm{ch}^{et}(\mathcal{E})$ by \cite{art:borne2007}, Lemma C.2. The definition of $\mathrm{td}^{rep}(\mathscr{X})$ is more involved, but in the following we will use only the fact that $[\mathrm{td}^{rep}(\mathscr{X})]_{1}$ is the usual Todd class $\mathrm{td}^{et}(\mathscr{X})$ in $H^\bullet _{et}(\mathscr{X})$.
\end{remark}

Now we introduce the following condition on $\mathcal{G}$:
\begin{condition}\label{item:integral} The number
\begin{equation}
\int^{et}_{\mathcal{I}(\mathscr{X})\setminus \mathscr{X}} \left[\left(\mathrm{ch}^{rep}(\mathcal{E})-\mathrm{ch}^{rep}(\mathcal{O}_{\mathscr{X}}^{\oplus \mathrm{ork}(\mathcal{E})})\right)\mathrm{ch}^{rep}(\mathcal{G}^\vee)\mathrm{c}_1^{rep}(\pi^\ast \mathcal L)\mathrm{td}^{rep}(\mathscr{X})\right]_{\ne 1}
\end{equation}
is zero for all coherent sheaves $\mathcal E$ on $\mathscr X$ and all ample line bundles $\mathcal L$ on $X$.
\end{condition}
 {\begin{remark}\label{rem:condition}
Recall that $\mathscr{X}$ is an orbifold, i.e., $\mathcal{I}(\mathscr{X})$ has exactly one two-dimensional component, which is $\mathscr{X}$ itself. We point out   that if $\mathcal{I}(\mathscr{X})\setminus \mathscr{X}$ has no one-dimensional components, the previous condition is trivially satisfied. If $\mathcal{I}(\mathscr{X})\setminus \mathscr{X}$ has one-dimensional components, the previous condition can be restated by saying that the zero degree part of
\begin{equation*}
\left(\mathrm{ch}^{rep}(\mathcal{E})-\mathrm{ch}^{rep}(\mathcal{O}_{\mathscr{X}}^{\oplus \mathrm{rk}(\mathcal{E})})\right)\mathrm{ch}^{rep}(\mathcal{G}^\vee)
\end{equation*}
is zero over the one-dimensional component of $\mathcal{I}(\mathscr{X})\setminus \mathscr{X}$ for any coherent sheaf $\mathcal E$ on $\mathscr X$.
\end{remark}}

\begin{lemma}\label{lem:degree}
Let $\mathscr X$ be a two-dimensional projective irreducible orbifold with coarse moduli scheme $\pi\colon \mathscr{X}\to X$ a normal projective surface and $\mathcal{G}$ a generating sheaf for it. Assume that condition \ref{item:integral} holds. Fix a one-dimensional smooth integral closed substack $\mathscr{D}\subset \mathscr{X}$, whose coarse moduli space $\mathscr D\to D$ is a smooth curve and a good framing divisor. Let $H$ be an ample divisor on $X$ and set $H_n=H+na_D D$ for any positive integer $n$. Then for any coherent sheaf $\mathcal{E}$ we have
\begin{equation}
\deg_{\mathcal{G},H_n}(\mathcal{E})=\deg_{\mathcal{G},H}(\mathcal{E})+n\,a_D\, a_{\mathscr{D}}\,\mathrm{ork}(\mathcal{G})\, \int_{\mathscr{D}}^{et} \mathrm{c}_1^{et}(\mathcal{E}_{|\mathscr{D}})\ .
\end{equation}
\end{lemma}
\proof
By the T\"oen-Riemann-Roch theorem the degree of $\mathcal{E}$ with respect to the polarization $(\mathcal{G}, \mathcal{O}_X(H_n))$ is
\begin{equation}
\deg_{\mathcal{G},H_n}(\mathcal{E})=\int^{rep}_{\mathscr{X}} \left(\mathrm{ch}^{rep}(\mathcal{E})-\mathrm{ch}^{rep}(\mathcal{O}_{\mathscr{X}}^{\oplus \mathrm{ork}(\mathcal{E})})\right)\mathrm{ch}^{rep}(\mathcal{G}^\vee)\mathrm{c}_1^{rep}(\pi^\ast \mathcal{O}_X(H_n))\mathrm{td}^{rep}(\mathscr{X})\ .
\end{equation}
Using the decomposition \eqref{eq:decomposition} and condition \ref{item:integral}, we obtain
\begin{equation}
\deg_{\mathcal{G},H_n}(\mathcal{E})=\int^{et}_{\mathscr{X}} \left[\left(\mathrm{ch}^{rep}(\mathcal{E})-\mathrm{ch}^{rep}(\mathcal{O}_{\mathscr{X}}^{\oplus \mathrm{ork}(\mathcal{E})})\right)\mathrm{ch}^{rep}(\mathcal{G}^\vee)\mathrm{c}_1^{rep}(\pi^\ast \mathcal{O}_X(H_n))\mathrm{td}^{rep}(\mathscr{X})\right]_{1}\ .
\end{equation}
By \cite{art:borne2007}, Lemma C.2 and the identity $\mathrm{c}_1^{et}(\pi^\ast \mathcal{O}_X(H_n))=\pi^\ast \mathrm{c}_1^{et}(\mathcal{O}_X(H_n))$, we have
\begin{equation}
\deg_{\mathcal{G},H_n}(\mathcal{E})=\int^{et}_{\mathscr{X}} \left(\mathrm{ch}^{et}(\mathcal{E})-\mathrm{ch}^{et}(\mathcal{O}_{\mathscr{X}}^{\oplus \mathrm{ork}(\mathcal{E})})\right)\mathrm{ch}^{et}(\mathcal{G}^\vee)\pi^\ast \mathrm{c}_1^{et}(\mathcal{O}_X(H_n))\mathrm{td}^{et}(\mathscr{X})\ .
\end{equation}
Since the zero degree part of $\mathrm{ch}^{et}(\mathcal{E})$ is $\mathrm{ork}(\mathcal{E})$, the degree becomes
\begin{equation}
\deg_{\mathcal{G},H_n}(\mathcal{E})=\mathrm{ork}(\mathcal{G})\int^{et}_{\mathscr{X}} \mathrm{c}_1^{et}(\mathcal{E})\pi^\ast \mathrm{c}_1^{et}(\mathcal{O}_X(H_n))\ .
\end{equation}
Moreover $\mathrm{c}_1^{et}(\mathcal{O}_X(H_n))=\mathrm{c}_1^{et}(\mathcal{O}_X(H))+n\mathrm{c}_1^{et}(\mathcal{O}_X(a_D D))$, so that we have
\begin{equation}
\deg_{\mathcal{G},H_n}(\mathcal{E})=\deg_{\mathcal{G},H}(\mathcal{E})+\mathrm{ork}(\mathcal{G})\int_{\mathscr{X}} \mathrm{c}_1^{et}(\mathcal{E})\pi^\ast \mathrm{c}_1^{et}(\mathcal{O}_X(n a_D D))\ .
\end{equation}
By Formula \eqref{eq:equilinebundles} we get
\begin{equation}
\int_{\mathscr{X}} \mathrm{c}_1^{et}(\mathcal{E})\pi^\ast \mathrm{c}_1^{et}(\mathcal{O}_X(a_D D))=a_D\,a_{\mathscr D}\,\int_{\mathscr{X}} \mathrm{c}_1^{et}(\mathcal{E})\mathrm{c}_1^{et}(\mathcal{O}_{\mathscr{X}}(\mathscr{D}))=a_D\,a_{\mathscr D}\,\int_{\mathscr{D}} \mathrm{c}_1^{et}(\mathcal{E}_{|\mathscr{D}})\ .
\end{equation}
Thus we obtain the assertion. 
\endproof
By using similar computations as before, we get also the following result.
\begin{lemma}\label{lem:degree2}
Under the same hypotheses of Lemma \ref{lem:degree} we have
\begin{equation}
\deg_{\mathcal{G},H_n}(\mathcal{E}\otimes \mathcal{O}_{\mathscr{X}}(\mathscr{D}))=\deg_{\mathcal{G},H_n}(\mathcal{E})+\mathrm{ork}(\mathcal{E})\deg_{\mathcal{G},H_n}(\mathcal{O}_{\mathscr{X}}(\mathscr{D}))\ .
\end{equation}
\end{lemma}

\begin{theorem}\label{thm:framed-stability}
Let $\mathscr X$ be a two-dimensional projective irreducible orbifold with coarse moduli scheme $\pi\colon \mathscr{X}\to X$ a normal projective surface and $\mathcal{G}$ a generating sheaf given as direct sum of powers of a $\pi$-ample locally free sheaf. Assume that condition \ref{item:integral} holds. Fix a one-dimensional smooth integral closed substack $\mathscr{D}\subset \mathscr{X}$, whose coarse moduli space $\mathscr D\to D$ is a smooth curve containing the singular locus of $X$ and a good framing divisor. Let $\mathcal{F}_{\mathscr{D}}$ be a good framing sheaf on $\mathscr{D}.$ For any numerical polynomial $P\in \mathbb{Q}[n]$ of degree two, there exist an effective ample divisor $C$ on $X$ and a positive rational number $\delta_1$ such that all the $(\mathscr{D}, \mathcal{F}_{\mathscr{D}})$-framed sheaves on $\mathscr{X}$ with Hilbert polynomial $P$ are $\mu$-stable with respect to $\delta_1$ and the polarization $(\mathcal{G}, \mathcal{O}_X(C)).$
\end{theorem}
\proof
%By arguing along the lines of the proof of the analogous theorem for framed sheaves on smooth projective surfaces (\cite{art:bruzzomarkushevich2011}, Thm.~3.1), and using many of the results so far proved in this section, we get the assertion. Indeed, let 
Let $H$ be an effective ample divisor on $X$ and let $n$ be a positive integer. Set $H_n=H+na_D D$. From now on, we shall use the polarizations $(\mathcal{G}, H)$ and $(\mathcal{G}, H_n)$ on $\mathscr{X}.$

Let us fix a numerical polynomial $P$ of degree two. The family of $(\mathscr{D}, \mathcal{F}_{\mathscr{D}})$-framed sheaves $\mathfrak{E}=(\mathcal{E},\phi_{\mathcal{E}})$ with Hilbert polynomial $P$ on $\mathscr{X}$ is bounded by Theorem \ref{thm:boundedness}. Then by the stacky version of Grothendieck Lemma (cf.\ \cite{art:nironi2008}, Lemma 4.13) and the Equation \eqref{eq:relationslope}, there exists a nonnegative constant $A_1$, depending only on $\mathcal{F}_{\mathscr D}$, $P$, $H$, such that for any $(\mathscr{D}, \mathcal{F}_{\mathscr{D}})$-framed sheaf $\mathfrak{E}=(\mathcal{E},\phi_{\mathcal{E}})$ with Hilbert polynomial $P$ on $\mathscr{X}$ and for any nonzero subsheaf $\mathcal{E}'\subset \mathcal{E}$
\begin{equation}
\mu_{\mathcal{G}, H}(\mathcal{E}')<\mu_{\mathcal{G}, H}(\mathcal{E})+A_1\ .
\end{equation}
Now we  {should} check that there exists $n$ such that the range of positive rational numbers $\delta_1$, for which all the $(\mathscr{D}, \mathcal{F}_{\mathscr{D}})$-framed sheaves with Hilbert polynomial $P$ on $\mathscr{X}$ are $\mu$-stable with respect to $\delta_1$ and the polarization $(\mathcal{G}, H_n)$, is nonempty. {By using the same arguments as in the proof of the analogous theorem for framed sheaves on smooth projective surfaces (\cite{art:bruzzomarkushevich2011}, Thm.~3.1) to study the $\mu$-stability condition (in particular  Lemma \ref{lem:locallyfree} is used in an essential way to compare the degree of a saturated subsheaf of $\mathcal{E}$ with the degree of its restriction to $\mathscr{D}$), this range of positive rational numbers $\delta_1$ is nonempty if} 

\begin{equation}
n>\max\left\{\frac{rA_1-\frac{\mathrm{ork}(\mathcal{G})(a_D D)\cdot H}{a_D k_{\mathscr{D}}}}{a_D\,a_{\mathscr{D}}\,\mathrm{ork}(\mathcal{G})\left(\frac{(a_D D)^2}{a_D^2 k_{\mathscr{D}}^2}-rA_0\right)},0\right\}\ .
\end{equation}
\endproof
\begin{remark}
When $\mathscr{X}=X$ is a smooth projective surface and $\mathcal{G}\simeq  \mathcal{O}_X$, this  proof reduces to the proof of \cite{art:bruzzomarkushevich2011}, Thm.~3.1. 
\end{remark}

By Theorems \ref{thm:deformations} and \ref{thm:framed-stability} we eventually have:
\begin{corollary} 
Under the same assumptions as in Theorem \ref{thm:framed-stability}, there exists a fine moduli space $\mathrm{M}_{\mathscr{X}/k}(P_0,\mathscr{D}, \mathcal{F}_{\mathscr{D}})$ parameterizing isomorphism classes of $(\mathscr{D}, \mathcal{F}_{\mathscr{D}})$-framed sheaves $(\mathcal E,\phi_{\mathcal E})$ on $\mathscr{X}$ with Hilbert polynomial $P$, which is a quasi-projective scheme. Its tangent space at a point $[(\mathcal E,\phi_{\mathcal E})]$ is $\operatorname{Ext}^1(\mathcal E,\mathcal E(-\mathscr{D}))$. If $\operatorname{Ext}^2(\mathcal E,\mathcal E(-\mathscr{D}))=0$ for all the points $[(\mathcal E,\phi_{\mathcal E})]$, the moduli space $\mathrm{M}_{\mathscr{X}/k}(P_0,\mathscr{D},$ $\mathcal{F}_{\mathscr{D}})$ is a smooth quasi-projective variety.
\end{corollary}
\begin{remark}
Since the moduli space $\mathrm{M}_{\mathscr{X}/k}(P_0,\mathscr{D}, \mathcal{F}_{\mathscr{D}})$ is fine, there exists a {\em universal} flat family $(\tilde{\mathcal{E}},L_{\tilde{\mathcal{E}}},\phi_{\tilde{\mathcal{E}}})$ of $(\mathscr{D}, \mathcal{F}_{\mathscr{D}})$-framed sheaves on $\mathscr{X}$ parameterized by $\mathrm{M}_{\mathscr{X}/k}(P_0,\mathscr{D}, \mathcal{F}_{\mathscr{D}})$. The fact the framing of a $(\mathscr{D}, \mathcal{F}_{\mathscr{D}})$-framed sheaf is an isomorphism after restricting to $\mathscr{D}$ implies that $\tilde{\phi}_{\tilde{\mathcal{E}}}\colon \tilde{\mathcal{E}} \to p_{\mathscr{X}}^\ast \mathcal{F}_{\mathscr{D}}$ is an isomorphism over $\mathrm{M}_{\mathscr{X}/k}(P_0,\mathscr{D}, \mathcal{F}_{\mathscr{D}})\times \mathscr{D}$. Moreover, this allows one to dispose of the homothety in the definition of morphisms of  $(\mathscr{D}, \mathcal{F}_{\mathscr{D}})$-framed sheaves, so that the line bundle $L_{\tilde{\mathcal{E}}}$ can be taken trivial.
\end{remark}

\bigskip\section{$(\mathscr{D}, \mathcal{F}_{\mathscr{D}})$-framed sheaves on two-dimensional projective toric orbifolds}\label{FramedsheavesToricorbifolds}

In this section we apply the theory of $(\mathscr{D},\mathcal{F}_{\mathscr{D}})$-framed sheaves developed in the previous section to the case of toric orbifolds. Let $\pi^{can}\colon\mathscr{X}^{can}\to X$ be the {\em canonical toric orbifold} of a normal projective toric surface $X$ and $\tilde{\mathscr{D}}\subset \mathscr{X}^{can}$ a smooth divisor whose coarse moduli scheme $D$ is a torus-invariant rational curve in $X$ containing the singular locus of $X$. By performing a {\em $k$-root construction} on $\mathscr{X}^{can}$ along $\tilde{\mathscr{D}}$ we obtain a two-dimensional projective toric orbifold $\mathscr{X}$, with coarse moduli scheme $X$, endowed with a smooth divisor $\mathscr{D}$ which is a $\mu_k$-gerbe over $\tilde{\mathscr{D}}$. We shall show that if $D$ is a good framing divisor and  {$\mathcal{O}_{\mathscr{X}^{can}}(\tilde{\mathscr{D}})$ is a $\pi^{can}$-ample line bundle}, Theorem \ref{thm:framed-stability} holds for any choice of a good framing sheaf $\mathcal{F}_{\mathscr{D}}$ on $\mathscr{D}$; hence for any numerical polynomial $P$ of degree two, there exists a fine moduli space parameterizing isomorphism classes of $(\mathscr{D},\mathcal{F}_{\mathscr{D}})$-framed sheaves on $\mathscr{X}$ with Hilbert polynomial $P$.

Our main reference for the theory of root stacks is \cite{art:cadman2007} (cf.\ also \cite{art:abramovichgrabervistoli2008}, Appendix~B); for the theory of toric stacks is \cite{art:fantechimannnironi2010} (see also \cite{art:borisovchensmith2004} for a combinatorial approach to the topic based on {\em stacky fans}). For the theory of toric varieties, see \cite{book:coxlittleschenck2011}.

In the following, we set $k=\mathbb{C}$.

\subsection{Root stack compactification of a smooth open toric surface}

Let $X$ be a normal projective toric surface, acted on by a torus $\mathbb C^\ast\times\mathbb C^\ast$, and let $\Sigma$ be its fan in $(N_2)_{\mathbb Q}$, where $N_2$ is the lattice of one-parameter subgroups of $\mathbb C^\ast\times\mathbb C^\ast$. Since $X$ is projective, the rays of $\Sigma$ generate  $(N_2)_{\mathbb Q}$. Let $n+2$ be the number of rays of $\Sigma$ for some positive integer $n$. By the orbit-cone correspondence there exist $n+2$ torus-invariant rational curves $D_0, \ldots, D_{n+1}$. We shall use also the letter $D$ to denote the curve $D_{n+1}$. 

The singular points of $X$ are necessarily torus-invariant, and, by the normality assumption, the singular locus $\mathrm{sing}(X)$ is zero-dimensional, i.e., $\mathrm{sing}(X)$ consists of a finite number of torus-fixed points. We assume that $\mathrm{sing}(X)$ is contained inside $D$. Then $\mathrm{sing}(X)$ consists at most of the two torus-fixed points of $D$, which we shall denote by $0,\infty$. Moreover, the complementary set $X_0:=X\setminus D$ is a smooth quasi-projective toric surface. Let us assume that the intersection point of $D_0$ and $D$ is $0$ and the intersection point of $D_n$ and $D$ is $\infty$. 

Let $\pi^{can}\colon\mathscr{X}^{can}\to X$ be the canonical toric orbifold of $X$ with Deligne-Mumford torus\footnote{A {\em Deligne-Mumford torus} is a product $T\times [pt/ G]$ where $T$ is a ordinary torus and $G$ a finite abelian group.} $\mathbb C^\ast\times\mathbb C^\ast$ (cf.~\cite{art:fantechimannnironi2010}, Sect.~4). Since $\pi^{can}$ is an isomorphism over $X_{sm}$, the ``orbifold" structure of $\mathscr{X}^{can}$ lies (at most) at the {\em stacky} points $\tilde{p}_0:=(\pi^{can})^{-1}(0)_{red}$ and $\tilde{p}_\infty:=(\pi^{can})^{-1}(\infty)_{red}$. So we have that $\tilde{\mathscr{D}}_i\simeq D_i$ for $i=1, \ldots, n-1$ and $\tilde{\mathscr{D}}_j$ is an orbifold for $j=0, n,n+1$. Since the coarse moduli scheme of $\tilde{\mathscr{D}}_j$ is $\mathbb{P}^1$, the stack $\tilde{\mathscr{D}}_j$ is a so-called {\em spherical orbicurve} (cf.\ \cite{art:behrendnoohi2006}, Sect.~5) for $j=0, n,n+1$. Since the number of orbifold points is at most two, by \cite{art:behrendnoohi2006}, Prop.~5.5, we have that
\begin{equation}
\tilde{\mathscr{D}}_0\simeq \mathscr{F}(a_0, 1)\ ,\qquad \tilde{\mathscr{D}}_n\simeq \mathscr{F}(a_\infty, 1)\ , \qquad \tilde{\mathscr{D}}_{n+1}\simeq \mathscr{F}(a_0, a_\infty)\ ,
\end{equation}
where we denote by $\mathscr{F}(p,q)$ the {\em football} with two orbifold points of order $p$ and $q$ respectively, where $p$ and $q$ are positive integers. A football is a one-dimensional complete orbifold with coarse moduli scheme $\mathbb{P}^1$ and at most two orbifold points. Note that $\mathscr{F}(1,1)\simeq \mathbb{P}^1$. 

A well-known consequence of the construction of the coarse moduli space is the existence for any geometric point $p$ of $\mathscr{X}^{can}$ with image $x$ in $X$ of an \'etale neighbourhood $U\to X$ of $x$ such that $U\times_X \mathscr{X}^{can}$ is a neighbourhood of $p$ and is a quotient stack of the form $[Y/\mathrm{Stab}(p)]$, where $Y$ is a scheme. In particular, there is an \'etale neighbourhood $U$ of $0$ in $X$ such that $U\times_X \mathscr{X}^{can}$ is an \'etale neighbourhood of $\tilde{p}_0$ and is a quotient stack of the form $[V/\mu_{a_0}]$, where $V$ is a smooth variety. Then $U=V/\mu_{a_0}$. So $a_0$ is the order of the singularity of $X$ at $0$. Similarly, $a_\infty$ is the order of the singularity of $X$ at $\infty$.

Since all toric footballs are fibred products of root stacks over $\mathbb{P}^1$ (cf.\ \cite{art:fantechimannnironi2010}, Example~7.31), {we can realized them as {\em root stacks}}
\begin{equation}
\tilde{\mathscr{D}}_0\simeq  \sqrt[a_0]{0/\mathbb{P}^1} \ ,\qquad \tilde{\mathscr{D}}_n \simeq \sqrt[a_\infty]{\infty/\mathbb{P}^1}\qquad\mbox{and}\qquad \tilde{\mathscr{D}}_{n+1}\simeq \sqrt[a_0]{0/\mathbb{P}^1}\times_{\mathbb{P}^1} \sqrt[a_\infty]{\infty/\mathbb{P}^1}\ .
\end{equation}
Here for root stacks along effective Cartier divisors  we use the notation in \cite{art:fantechimannnironi2010}, Sect.~1.3.b.

Denote by $\tilde{\mathscr{D}}$ the smooth effective Cartier divisor $\tilde{\mathscr{D}}_{n+1}$. From now on, assume that the line bundle $\mathcal{O}_{\mathscr{X}^{can}}(\tilde{\mathscr{D}})$ is $\pi^{can}$-ample\footnote{At the moment we do not know how strong this constraint is. On the other hand, if $X$ is smooth, $\mathscr{X}^{can}$ coincides with $X$ and the assumption is automatically satisfied. Moreover, the canonical toric stacks constructed in \cite{art:bruzzopedrinisalaszabo2013} satisfy this assumption.}.

\begin{remark}\label{rem:injectivecharacter}
As explained in \cite{art:fantechimannnironi2010}, Rem.~4.12-(2), $\mathscr{X}^{can}$ is isomorphic to the global quotient $[Z/G]$, where $G:=\operatorname{Hom}(A^1(X),\mathbb{C}^*)$ and $Z=\mathbb{A}^{n+2}\setminus V(J_\Sigma)$, where $J_\Sigma$ is the so-called {\em irrelevant ideal of $\Sigma$}. Since $\mathscr{X}^{can}$ is an orbifold, $G$ is a subgroup of $(\mathbb{C}^*)^{n+2}$ and its rank is $n$. Moreover, the action of $G$ on $Z$ is given via the composition of the inclusion map of $G$ into $(\mathbb{C}^*)^{n+2}$ and the standard action of $(\mathbb{C}^*)^{n+2}$ on $\mathbb{A}^{n+2}$. By \cite{art:fantechimannnironi2010}, Cor.~4.13, the divisor $\tilde{\mathscr{D}}$ is isomorphic to $[(Z\cap \{z_{n+1}=0\}/G]$ (the coordinates of $Z$ are $z_0, \ldots, z_{n+1}$). Since $\mathrm{codim}\,V(J_\Sigma)\geq 2$, the line bundle $\mathcal{O}_{\mathscr{X}^{can}}(\tilde{\mathscr{D}})$ corresponds to the trivial line bundle on $Z$ with a $G$-equivariant structure induced by a character of $G$. By \cite{art:fantechimannnironi2010}, Rem.~4.14-(2), the character corresponding to $\mathcal{O}_{\mathscr{X}^{can}}(\tilde{\mathscr{D}})$ is $\chi_{n+1}\colon G\hookrightarrow (\mathbb{C}^*)^{n+2}\xrightarrow{p_{n+1}} \mathbb{C}^*$ (the coordinates of $(\mathbb{C}^*)^{n+2}$ are $\lambda_0, \ldots, \lambda_{n+1}$). By the ampleness assumption on $\mathcal{O}_{\mathscr{X}^{can}}(\tilde{\mathscr{D}})$ the composition of the inclusion of $\mu_{a_0}$ into $G$ and $\chi_{n+1}$ is injective and the same holds for $\mu_{a_\infty}$. We shall use this fact later on.
\end{remark}

Let $k$ be a positive integer and denote by $\mathscr{X}$ the root stack $\sqrt[k]{\tilde{\mathscr{D}}/\mathscr{X}^{can}}$. It is a two-dimensional toric orbifold with Deligne-Mumford torus $\mathbb{C}^\ast\times \mathbb{C}^\ast$ and with coarse moduli scheme $X$. By \cite{art:fantechimannnironi2010}, Prop.~5.1, the structure morphism $\pi\colon \mathscr{X}\to X$ factorizes as
\begin{equation}\label{eq:factorization}
\begin{aligned}
  \begin{tikzpicture}[xscale=3.5,yscale=-1.3,->, bend angle=40]
  \node (A0_0) at (0, 0) {$\mathscr{X}$};
\node (A1_0) at (1, 0) {$\mathscr{X}^{can}$};
\node (A2_0) at (2, 0) {$X$};
\node (A0_3) at (2.2, 0) {$,$};
    \path (A0_0) edge [bend right]node [auto] {$\scriptstyle{\pi}$} (A2_0);
        \path (A0_0) edge [->]node [below] {$\scriptstyle{\psi}$} (A1_0);
        \path (A1_0) edge [->]node [below] {$\scriptstyle{\pi^{can}}$} (A2_0);
  \end{tikzpicture}
\end{aligned}
\end{equation}
Moreover, by \cite{art:fantechimannnironi2010}, Lemma 7.1, $\mathscr{X}$ is isomorphic to the global quotient $[\tilde{Z}/\tilde{G}]$, where $\tilde{Z}$ and $\tilde{G}$ are defined by the following fibre products
\begin{equation}
  \begin{tikzpicture}[xscale=1.5,yscale=-1.2]
    \node (A0_0) at (0, 0) {$\tilde{Z}$};
    \node (A0_2) at (2, 0) {$\mathbb{A}^1$};
    \node (A1_1) at (1, 1) {$\square$};
    \node (A2_0) at (0, 2) {$Z$};
    \node (A2_2) at (2, 2) {$\mathbb{A}^1$};

    \node (A0_3) at (4, 0) {$\tilde{G}$};
    \node (A0_5) at (6, 0) {$\mathbb{C}^\ast$};
    \node (A1_4) at (5, 1) {$\square$};
    \node (A2_3) at (4, 2) {$G$};
    \node (A2_5) at (6, 2) {$\mathbb{C}^\ast$};
 
     \node (A1_3) at (7, 1) {$.$};
        
    \path (A0_0) edge [->]node [auto] {$\scriptstyle{}$} (A2_0);
    \path (A0_0) edge [->]node [auto] {$\scriptstyle{}$} (A0_2);
    \path (A0_2) edge [->]node [auto] {$\scriptstyle{(-)^{k}}$} (A2_2);
    \path (A2_0) edge [->]node [auto] {$\scriptstyle{z_{n+1}}$} (A2_2);

    \path (A0_3) edge [->]node [auto] {$\scriptstyle{\varphi}$} (A2_3);
    \path (A0_3) edge [->]node [auto] {$\scriptstyle{\tilde{\chi}_{n+1}}$} (A0_5);
    \path (A0_5) edge [->]node [auto] {$\scriptstyle{(-)^{k}}$} (A2_5);
    \path (A2_3) edge [->]node [auto] {$\scriptstyle{\chi_{n+1}}$} (A2_5);
  \end{tikzpicture}   
\end{equation}
The action of $\tilde{G}$ on $\tilde{Z}$ is given by
\begin{equation}
(g, \lambda)\cdot (z, x):= (gz, \lambda x)\ ,
\end{equation}
for any $(g, \lambda)\in \tilde{G}$ and $(z, x)\in\tilde{Z}$.

The effective Cartier divisor $\mathscr{D}:=\pi^{-1}(D)_{red}$ is an \'etale $\mu_k$-gerbe over $\tilde{\mathscr{D}}$. As a global quotient 
$\mathscr{D}$ is the stack $[Z\cap \{z_{n+1}=0\}/\tilde{G}]$, where the $\tilde{G}$-action is given via $\varphi$,  and $\ker\varphi=\{(1, \lambda)\,\vert\,\lambda^k=1\}\simeq \mu_k$. Moreover, the line bundle $\mathcal{O}_{\mathscr{X}}(\mathscr{D})$ corresponds to the morphism $\mathscr{X}\to [\mathbb{A}^1/\mathbb{G}_m]$ and then to the character $\tilde{\chi}_{n+1}$.

Now we check if the hypotheses of Theorem \ref{thm:framed-stability} hold for the pair $(\mathscr{X}, \mathscr{D})$. The first thing we shall prove is that the line bundle $\mathcal{O}_{\mathscr{X}}(\mathscr{D})$ is $\pi$-ample. 
Since $\tilde{\mathscr{D}}$ is the rigidification of $\mathscr{D}$ with respect to $\mu_k$ (cf.\ \cite{art:fantechimannnironi2010}, Sect.~6.3), the stabilizer group of a geometric point $p$ of $\mathscr{D}$  is an extension
\begin{equation}
1\to \ker\varphi \to \mathrm{Stab}(p)\to \mathrm{Stab}(\tilde{p})\to 1\ ,
\end{equation}
where $\tilde{p}:=\psi(p)\in \tilde{\mathscr{D}}$. In particular, if $\tilde{p}$ is not   $\tilde{p}_0$ or $\tilde{p}_\infty$, the stabilizer group of $p$ is $\ker\varphi$. Since the character $(\tilde{\chi}_{n+1})_{\vert \ker\varphi}$ is $(1,\lambda)\mapsto \lambda$, the representation of the stabilizer group at the fibre of $\mathcal{O}_{\mathscr{X}}(\mathscr{D})$ at $p$ is faithful. If $\tilde{p}=\tilde{p}_0$, denote by $p_0$ the corresponding geometric point in $\mathscr{D}$. The kernel of the character $(\tilde{\chi}_{n+1})_{\vert \mathrm{Stab}(p_0)}$ is the set $\{(g,1)\,\vert\, g\in \mathrm{Stab}(\tilde{p}_0)\mbox{ and }\chi_{n+1}(g)=1\}$. By Remark \ref{rem:injectivecharacter}, $(\chi_{n+1})_{\vert \mathrm{Stab}(\tilde{p}_0)}$ is injective, and $(\tilde{\chi}_{n+1})_{\vert \mathrm{Stab}(p_0)}$ is injective as well. Hence the representation of $\mathrm{Stab}(p_0)$ on the fibre of $\mathcal{O}_{\mathscr{X}}(\mathscr{D})$ at the point $p_0$ is faithful. One can argue similarly for the geometric point $p_\infty\in \mathscr{D}$ such that $\psi(p_\infty)=\tilde{p}_\infty$. Thus $\mathcal{O}_{\mathscr{X}}(\mathscr{D})$ is $\pi$-ample. Therefore,   
\begin{equation}
\mathcal{G}:=\oplus_{i=1}^r \mathcal{O}_{\mathscr{X}}(\mathscr{D})^{\otimes i}
\end{equation}
is a generating sheaf for $\mathscr{X}$ for any positive integer $r\geq N_{\mathscr{X}}$, where $N_{\mathscr{X}}=\max\{k\,a_0,k\,a_\infty\}$, by Proposition \ref{prop:veryampleness}. We fix a positive integer $a$ such that $r:=k\,a\geq N_{\mathscr{X}}$.

Now we check that Condition \ref{item:integral} holds. We shall use some arguments of \cite{art:borne2007}, Sect.~4.2.4. The inertia stack $\mathcal{I}(\mathscr{X})$ of $\mathscr{X}$ has only one two-dimensional component, i.e., the stack $\mathscr{X}$ associated with the trivial stabilizer. The one-dimensional components of $\mathcal{I}(\mathscr{X})$ are $\bigsqcup_{j=1}^{k-1} \mathscr{D}^j$, where $\sigma_{\mathscr X}(\mathscr D^j)=\mathscr D$ (here $\sigma_{\mathscr{X}}\colon \mathcal{I}(\mathscr{X})\to \mathscr{X}$ is the forgetful morphism), hence $\mathcal{I}(\mathscr{X})\setminus \mathscr{X}$ has a nontrivial one-dimensional component. On the other hand, the one-dimensional component of the inertia stack $\mathcal{I}(\mathscr{D})$ of $\mathscr{D}$ is $\bigsqcup_{j=0}^{k-1}\mathcal{I}(\mathscr{D})^j$, where $\sigma_{\mathscr{D}}(\mathcal{I}(\mathscr{D})^j)=\mathscr{D}$ for any $j=1, \ldots, k-1$ (here $\sigma_{\mathscr{D}}\colon \mathcal{I}(\mathscr{D})\to \mathscr{D}$ is the forgetful morphism). 

After fixing a primitive $k$-root of unity $\omega$, the substack $\mathcal{I}(\mathscr{D})^j$ is associated with the automorphism induced by the multiplication by $\omega^j$ for $j=0, \ldots, k-1$. Thus --- roughly speaking --- $\mathcal{I}(\mathscr{D})^j$ consists of pairs of the form $(p, \omega^j)$, where $p$ is a point of $\mathscr{D}$. The same argument holds for $\mathscr D^j$ with $j=1, \ldots, k-1$. 
Let us denote by $\imath\colon \mathscr{D}\to \mathscr{X}$ the inclusion morphism and by $\mathcal{I}(\imath)\colon \mathcal{I}(\mathscr{D})\setminus \mathscr{D}\to \mathcal{I}(\mathscr{X})\setminus \mathscr{X}$ the corresponding inclusion morphism at the level of inertia stacks. Then $\mathcal{I}(\imath)_{\vert \mathcal{I}(\mathscr{D})^j}\colon \mathcal{I}(\mathscr{D})^j\to \mathscr D^j$ is an isomorphism for any $j=1, \ldots, k-1$.  
\begin{equation}
x:=\left(\mathrm{ch}^{rep}(\mathcal{E})-\mathrm{ch}^{rep}(\mathcal{O}_{\mathscr{X}}^{\oplus \mathrm{ork}(\mathcal{E})})\right)\mathrm{ch}^{rep}(\mathcal{G}^\vee)\mathrm{c}_1^{rep}(\pi^\ast \mathcal L)\mathrm{td}^{rep}(\mathscr{X})\ .
\end{equation}
Since the integral of $x_{\ne 1}$ is zero over the zero-dimensional components of $\mathcal{I}(\mathscr{X})\setminus \mathscr{X}$, we have  
\begin{equation}
\int^{et}_{\mathcal{I}(\mathscr{X})\setminus \mathscr{X}} x_{\ne 1}=\int^{et}_{\mathcal{I}(\mathscr{D})\setminus \mathscr{D}} \mathcal{I}(\imath)^\ast(x_{\ne 1})=\int^{et}_{\mathcal{I}(\mathscr{D})\setminus \mathscr{D}} [\imath^\ast x]_{\neq 1}\ .
\end{equation}
Now, note that $\displaystyle\int^{et}_{\mathscr{D}} [\imath^\ast x]_{1}=0$.  Indeed, $[\mathcal{I}(\imath)^\ast\sigma_{\mathscr X}^\ast\mathrm{c_1}^{rep}(\mathcal{O}_X(1))]_1=\imath^\ast\pi^\ast\mathrm{c_1}^{et}(\mathcal{O}_X(1))$ and the degree zero part of the difference
\begin{equation}
\left[\mathrm{ch}^{rep}(\imath^\ast\mathcal{E})-\mathrm{ch}^{rep}(\mathcal{O}_{\mathscr{D}}^{\oplus \mathrm{ork}(\mathcal{E})})\right]_1=\mathrm{ch}^{et}(\mathcal{F}_{\mathscr{D}})-\mathrm{ch}^{et}(\mathcal{O}_{\mathscr{D}}^{\oplus \mathrm{ork}(\mathcal{E})})
\end{equation}
is zero. Thus we get
\begin{equation}
\int^{et}_{\mathcal{I}(\mathscr{D})\setminus \mathscr{D}} [\imath^\ast x]_{\neq 1}=\int^{rep}_{\mathscr{D}} \imath^\ast x\ .
\end{equation}
It remains to prove that the last integral is zero. Let us fix $j\in\{1, \ldots, k-1\}$. By \cite{art:borne2007}, Lemma 4.14, we have
\begin{equation}
\mathrm{ch}^{rep}(\mathcal{G}^\vee_{\vert \mathcal{I}(\mathscr{D})^j})=\mathrm{ch}^{rep}(\oplus_{i=1}^r \mathcal{O}_{\mathscr{X}}(\mathscr{D})^{\otimes -i}_{\vert \mathcal{I}(\mathscr{D})^j})=\sum_{i=1}^r \omega^{-ij}\mathrm{ch}^{et}(\mathcal{O}_{\mathscr{X}}(\mathscr{D})^{\otimes -i}_{\vert \mathcal{I}(\mathscr{D})^j})\ ,
\end{equation}
So the zero degree part of it over $\mathcal{I}(\mathscr{D})^j$ is $\sum_{i=1}^r \omega^{-ij}$. Recall that
\begin{equation*}
\frac{1}{k}\sum_{i=0}^{k-1} \omega^{is}=\left\{
\begin{array}{ll}
0 & s\not\equiv 0\,\mathrm{mod}\,k\ ,\\
1 & s\equiv 0\,\mathrm{mod}\,k\ .
\end{array}\right.
\end{equation*}
Thus 
\begin{equation}
\sum_{i=1}^r \omega^{-ij}=\sum_{l=1}^a\sum_{i=k(l-1)}^{kl-1} \omega^{-ij}+\omega^{-rj}-1=\sum_{l=1}^a\sum_{i=0}^{k-1} \omega^{-ij+k(l-1)j}+1-1=a\sum_{i=0}^{k-1} \omega^{-ij}=0\ .
\end{equation}
Since the zero degree part of $\mathrm{ch}^{rep}(\imath^\ast\mathcal{E})-\mathrm{ch}^{rep}(\mathcal{O}_{\mathscr{D}}^{\oplus \mathrm{ork}(\mathcal{E})})$ is zero over $\mathcal{I}(\mathscr{D})^0$ and the zero degree part of $\mathrm{ch}^{rep}(\mathcal{G}^\vee_{\vert \mathscr{D}})$ is zero over $\mathcal{I}(\mathscr{D})^j$ for $j=1, \ldots, k-1$, the zero degree part of 
\begin{equation}
\Big(\mathrm{ch}^{rep}(\imath^\ast\mathcal{E})-\mathrm{ch}^{rep}(\mathcal{O}_{\mathscr{D}}^{\oplus \mathrm{ork}(\mathcal{E})})\Big)\mathrm{ch}^{rep}(\mathcal{G}^\vee_{\vert \mathscr{D}})
\end{equation}
is zero over $\mathcal{I}(\mathscr{D})$ and this implies that $\displaystyle\int^{rep}_{\mathscr{D}} \imath^\ast x=0$.

Theorem \ref{thm:framed-stability} implies the following result.
\begin{theorem}\label{thm:moduli-toricroot}
Let $X$ be a normal projective toric surface and $D$ a torus-invariant rational curve which contains the singular locus $\mathrm{sing}(X)$ of $X$ and is a good framing divisor. Let $\pi^{can}\colon\mathscr{X}^{can}\to X$ be the canonical toric orbifold of $X$ and $\tilde{\mathscr{D}}$ the smooth effective Cartier divisor $(\pi^{can})^{-1}(D)_{red}$. Assume that the line bundle $\mathcal{O}_{\mathscr{X}^{can}}(\tilde{\mathscr{D}})$ is $\pi^{can}$-ample. Let $\mathscr{X}:=\sqrt[k]{\tilde{\mathscr{D}}/\mathscr{X}^{can}}$, for some positive integer $k$, and $\mathscr{D}\subset\mathscr{X}$ the effective Cartier divisor corresponding to the morphism $\mathscr{X}\to [\mathbb{A}^1/\mathbb{G}_m]$. Then for any good framing sheaf $\mathcal{F}_{\mathscr{D}}$ on $\mathscr{D}$ and any numerical polynomial $P\in \mathbb{Q}[n]$ of degree two, there exists a fine moduli space parameterizing isomorphism classes of $(\mathscr{D}, \mathcal{F}_{\mathscr{D}})$-framed sheaves on $\mathscr{X}$ with Hilbert polynomial $P$, which is a quasi-projective scheme over $\mathbb{C}$.
\end{theorem}

\appendix

\bigskip\section{A semicontinuity theorem for the Hom group of framed sheaves}

In this section we provide a semicontinuity result for framed sheaves on stacks. As usual we need to start from a semicontinuity result in commutative algebra.

\begin{lemma}\label{lem:semicontinuity}
Let $A$ be a Noetherian ring. Let $L^\bullet$ be a complex of finite generated free $A$-modules bounded above. The function
\begin{equation}
p\mapsto \dim_{A_p/m_p} h^i(L^\bullet\otimes_A A_p/m_p)
\end{equation}
is upper semicontinuous for any $i$ and for any $p\in \mathrm{Spec}(A)$, where $m_p$ is the unique maximal ideal of the local ring $A_p$.
\end{lemma}
\proof
The proof follows from the same arguments as in the proof of the usual semicontinuity theorem, see e.g.~\cite{book:hartshorne1977}, Thm.~12.8.
\endproof

\begin{proposition}\label{prop:semicontinuity}
Let $\mathscr{X}$ be a projective stack and $\mathcal{F}$ a framing sheaf on it. Let $\mathfrak{E}=(\mathcal{E},L_{\mathcal{E}},\phi_{\mathcal{E}})$ and $\mathfrak{H}=(\mathcal{H},L_{\mathcal{H}},\phi_{\mathcal{H}})$ be flat families of framed sheaves on $\mathscr{X}$ parameterized by a scheme $S$ of finite type over $k$. Assume that $(\tilde{\phi}_\mathcal{E})_s$ is zero either for all $s\in S$ or for none. Then the function
\begin{equation}
s\mapsto \dim_{k(s)}\operatorname{Hom}((\mathcal{E}_s, (\tilde{\phi}_{\mathcal{E}})_s), (\mathcal{H}_s, (\tilde{\phi}_{\mathcal{H}})_s))
\end{equation}  
is upper semicontinuous. % in $s\in S$.
\end{proposition}
\proof
We combine the arguments of the proofs of \cite{art:huybrechtslehn1995-II}, Lemma 3.4 and \cite{art:nironi2008}, Lemma 6.18: the first is a semicontinuity theorem for Hom groups of flat families of framed sheaves on smooth projective varieties, while the second is the semicontinuity theorem for Hom groups of flat families of coherent sheaves on families of projective stacks.

The problem is local in $S$ so we can shrink $S$ to an affine scheme $\mathrm{Spec}(A)$, where $A$ is a $k$-algebra of finite type, such that $L_{\mathcal{E}}$ and $L_{\mathcal{H}}$ are trivial line bundles. Since $\mathscr{X}\times S\to S$ is a family of projective stacks, there exists a locally free resolution 
\begin{multline}
p_{\mathscr{X}\times S, \mathscr{X}}^\ast(\mathcal{G})^{\oplus N_1}\otimes p_{\mathscr{X}\times S, \mathscr{X}}^\ast\pi^\ast \mathcal{O}_X(-m_1)\to \\  p_{\mathscr{X}\times S, \mathscr{X}}^\ast(\mathcal{G})^{\oplus N_2}\otimes p_{\mathscr{X}\times S, \mathscr{X}}^\ast\pi^\ast \mathcal{O}_X(-m_2)\to \mathcal{E}\to 0\ ,
\end{multline}
where $m_1, m_2$ are large enough integers (cf.\ \cite{book:hartshorne1977}, Cor.~II-5.18). Then for any $A$-module $M$ we have the following exact sequences
\begin{multline}
0\to \operatorname{Hom}(\mathcal{E}, \mathcal{H}\otimes_A M)\to H^0(X\times S, F_{p_{\mathscr{X}\times S, 
\mathscr{X}}^\ast \mathcal{G}}(\mathcal{H})^{\oplus N_2}(m_2))\otimes_A M\\
\to H^0(X\times S, F_{p_{\mathscr{X}\times S, \mathscr{X}}^\ast \mathcal{G}}(\mathcal{H})^{\oplus N_1}(m_1))\otimes_A M\ ,
\end{multline}
\begin{multline}
0\to \operatorname{Hom}(\mathcal{E}, p_{\mathscr{X}\times S, \mathscr{X}}^\ast \mathcal{F}\otimes_A M)\to  H^0(X\times S, F_{p_{\mathscr{X}\times S, \mathscr{X}}^\ast \mathcal{G}}(p_{\mathscr{X}\times S, \mathscr{X}}^\ast \mathcal{F})^{\oplus N_2}(m_2))\otimes_A M\\ 
\to   H^0(X\times S, F_{p_{\mathscr{X}\times S, \mathscr{X}}^\ast \mathcal{G}}(p_{\mathscr{X}\times S, \mathscr{X}}^\ast \mathcal{F})^{\oplus N_1}(m_1))\otimes_A M\ .
\end{multline}
As explained in the proof of \cite{art:nironi2008}, Lemma 6.18, since $\mathcal{H}$ and $p_{\mathscr{X}\times S, \mathscr{X}}^\ast \mathcal{F}$ are $A$-flat, for $m_1, m_2$ sufficiently large the $A$-modules
\begin{align}
K^i_{\mathcal{H}}&:=H^0(X\times S, F_{p_{\mathscr{X}\times S, \mathscr{X}}^\ast \mathcal{G}}(\mathcal{H})^{\oplus N_i}(m_i))\ ,\\ K^i_{\mathcal{F}}&:=H^0(X\times S, F_{p_{\mathscr{X}\times S, \mathscr{X}}^\ast \mathcal{G}}(p_{\mathscr{X}\times S, \mathscr{X}}^\ast \mathcal{F})^{\oplus N_i}(m_i))
\end{align}
are free for $i=1, 2$. Define the complexes (concentrated in degrees zero and one)
\begin{equation}
K^\bullet_{\mathcal{H}}\colon 0\to K^2_{\mathcal{H}}\to K^1_{\mathcal{H}}\to 0\qquad\mbox{and}\qquad K^\bullet_{\mathcal{F}}\colon 0\to K^2_{\mathcal{F}}\to K^1_{\mathcal{F}}\to 0\ .
\end{equation}
Then for any $A$-module $M$ we have
\begin{equation}
\operatorname{Hom}(\mathcal{E}, \mathcal{H}\otimes_A M)\simeq H^0(K^\bullet_\mathcal{H}\otimes_A M)\;\;\mbox{and}\;\; \operatorname{Hom}(\mathcal{E}, p_{\mathscr{X}\times S, \mathscr{X}}^\ast \mathcal{F}\otimes_A M) \simeq  H^0(K^\bullet_\mathcal{F}\otimes_A M)\ .
\end{equation}
An element  $\phi\in K^2_\mathcal{F}$ which represents the framing $\phi_\mathcal{E}$  defines a homomorphism of complexes $\phi\colon A^\bullet \to K^\bullet_\mathcal{F}$, where $A^\bullet$ is the complex $A^0=A$ and $A^i=0$ for $i\neq 0$. Consider the homomorphism of complexes $\psi\colon K^\bullet_\mathcal{H}\oplus A^\bullet\to K^\bullet_\mathcal{F}$ induced by $(\phi_\mathcal{H}, -\phi)$ and let $\mathrm{Cone}(\psi)^\bullet$ be the cone of $\psi$. Then we have the short exact sequence
\begin{equation}
0\to K^\bullet_\mathcal{F}\to \mathrm{Cone}(\psi)^\bullet \to (K^\bullet_\mathcal{H}\oplus A^\bullet)[1]\to 0
\end{equation}
of bounded complexes of free $A$-modules. Thus we get the long exact sequence in cohomology
\begin{align}
0=\operatorname{Ext}^{-1}(\mathcal{E}, p_{\mathscr{X}\times S, \mathscr{X}}^\ast \mathcal{F}\otimes_A M)&\to h^{-1}(\mathrm{Cone}(\psi)^\bullet\otimes_A M)\to \operatorname{Hom}(\mathcal{E}, \mathcal{H}\otimes_A M)\oplus M\\ &\to \operatorname{Hom}(\mathcal{E}, p_{\mathscr{X}\times S, \mathscr{X}}^\ast \mathcal{F}\otimes_A M)\to \cdots
\end{align}
In particular for any $s\in S$ and $M=k(s)$ there is a cartesian diagram
\begin{equation}
  \begin{tikzpicture}[xscale=2.5,yscale=-1.2]
    \node (A0_0) at (0, 0) {$h^{-1}(\mathrm{Cone}(\psi)^\bullet\otimes k(s))$};
    \node (A0_2) at (2, 0) {$k(s)$};
    \node (A1_1) at (1, 1) {$\square$};
    \node (A2_0) at (0, 2) {$\operatorname{Hom}(\mathcal{E}, \mathcal{H}\otimes_A k(s))$};
    \node (A2_2) at (2, 2) {$\operatorname{Hom}(\mathcal{E}, p_{\mathscr{X}\times S, \mathscr{X}}^\ast \mathcal{F}\otimes_A k(s))$};
    \path (A0_0) edge [->]node [auto] {$\scriptstyle{}$} (A2_0);
    \path (A0_0) edge [->]node [auto] {$\scriptstyle{}$} (A0_2);
    \path (A0_2) edge [->]node [auto] {$\scriptstyle{\cdot \phi_{\mathcal{E}}}$} (A2_2);
    \path (A2_0) edge [->]node [auto] {$\scriptstyle{\phi_{\mathcal{H}}\circ}$} (A2_2);
  \end{tikzpicture} 
\end{equation}
Therefore $\dim_{k(s)}\operatorname{Hom}((\mathcal{E}_s, (\tilde{\phi}_{\mathcal{E}})_s), (\mathcal{H}_s, (\tilde{\phi}_{\mathcal{H}})_s))=\dim_{k(s)} h^{-1}(\mathrm{Cone}(\psi)^\bullet\otimes k(s))-1+\epsilon((\tilde{\phi}_\mathcal{E})_s)$ (cf.\ Remark \ref{rem:morphisms}). By Lemma \ref{lem:semicontinuity}, the function $s\mapsto \dim_{k(s)} h^{-1}(\mathrm{Cone}(\psi)^\bullet\otimes k(s))$ is upper semicontinuous. Since $(\tilde{\phi}_\mathcal{E})_s$ is zero either for all $s\in S$ or for none, the function $$s\mapsto \dim_{k(s)}\operatorname{Hom}((\mathcal{E}_s, (\tilde{\phi}_{\mathcal{E}})_s), (\mathcal{H}_s, (\tilde{\phi}_{\mathcal{H}})_s))$$ is upper semicontinuous as well.
\endproof

\bigskip\section{Serre duality for smooth projective stacks}\label{sec:serreduality}

In this section we prove Serre duality theorems for coherent sheaves on smooth projective stacks. Since these results are only sketched  in Nironi's papers \cite{art:nironi2008, art:nironi2008-II},  we give here a more complete treatment. First, we recall two results from \cite{art:nironi2008-II}, Thm.~1.16 and Cor.~2.10. Then we prove Serre duality for Deligne-Mumford stacks.
\begin{proposition}
Let $f\colon \mathscr{X}\to \mathscr{Y}$ be a proper morphism of separated Deligne-Mumford stacks of finite type over $k$. The functor $\mathbf{R}f_\ast\colon \mathrm{D}^+(\mathscr{X})\to \mathrm{D}^+(\mathscr{Y})$ has a right adjoint $f^!\colon \mathrm{D}^+(\mathscr{Y})\to \mathrm{D}^+(\mathscr{X})$. Moreover, for $\mathcal{E}^\bullet \in \mathrm{D}_c^+(\mathscr{X})$ and $\mathcal{F}^\bullet\in \mathrm{D}^+(\mathscr{Y})$ the natural morphism
\begin{equation}\label{eq:duality}
\mathbf{R}f_\ast\mathbf{R}\mathcal{H}om_{\mathscr{X}}(\mathcal{E}^\bullet, f^!\mathcal{F}^\bullet)\to \mathbf{R}\mathcal{H}om_{\mathscr{Y}}(\mathbf{R}f_\ast\mathcal{E}^\bullet, \mathbf{R}f_\ast f^!\mathcal{F}^\bullet)\xrightarrow{\mathrm{tr}_{f}} \mathbf{R}\mathcal{H}om_{\mathscr{Y}}(\mathbf{R}f_\ast\mathcal{E}^\bullet, \mathcal{F}^\bullet)
\end{equation}
is an isomorphism.
\end{proposition}
\begin{theorem}{(Serre duality - I)}\label{thm:serreduality-I}
Let $p\colon \mathscr{X}\to \mathrm{Spec}(k)$ be a proper Cohen-Macaulay Deligne-Mumford stack of pure dimension $d$. For any coherent sheaf $\mathcal{E}$ on $\mathscr{X}$ one has
\begin{equation}
H^i(\mathscr{X},\mathcal{E})^\vee \simeq  \operatorname{Ext}^{d-i}(\mathcal{E}, \omega_{\mathscr{X}})\ ,
\end{equation}
where $\omega_{\mathscr{X}}$ is the dualizing sheaf of $\mathscr{X}$.
\end{theorem}
\proof
By \cite{art:nironi2008-II}, Cor.~2.30, $p^!\mathcal{O}_{\mathrm{Spec}(k)}$ is isomorphic to the complex $\omega_{\mathscr{X}}[d]$, where $\omega_{\mathscr{X}}$ is the dualizing sheaf of $\mathscr{X}$. Let $\mathcal{E}$ be a coherent sheaf on $\mathscr{X}$. By applying the Formula \eqref{eq:duality} to the coherent sheaves $\mathcal{E}$ and $\mathcal{O}_{\mathrm{Spec}(k)}$ (regarded in the derived category as complexes concentrated in degree zero), we obtain
\begin{equation}
\mathbf{R}p_\ast\mathbf{R}\mathcal{H}om_{\mathscr{X}}(\mathcal{E}, \omega_{\mathscr{X}}[d])\xrightarrow{\sim} \mathbf{R}\mathcal{H}om_{\mathrm{Spec}(k)}(\mathbf{R}p_\ast\mathcal{E}, \mathcal{O}_{\mathrm{Spec}(k)})\simeq  \mathbf{R}\Gamma(\mathscr{X},\mathcal{E})^\vee\ .
\end{equation}
By taking cohomology, we get for any $i\geq 0$
\begin{equation}
\operatorname{Hom}_{\mathrm{D}(\mathscr{X})}(\mathcal{E}, \omega_{\mathscr{X}}[d-i]) \simeq  H^i(\mathscr{X},\mathcal{E})^\vee\ .
\end{equation}
Since the category of quasi-coherent sheaves on $\mathscr{X}$ has enough injectives (\cite{art:nironi2008-II}, Prop.~1.13), we get $\operatorname{Ext}^{d-i}(\mathcal{E}, \omega_{\mathscr{X}})\simeq \operatorname{Hom}_{\mathrm{D}(\mathscr{X})}(\mathcal{E}, \omega_{\mathscr{X}}[d-i])$, and therefore we obtain the desired result.
\endproof
Now we would like to prove a Serre duality theorem for Ext groups. We readapt the proof of the analogous theorem in the case of coherent sheaves on proper Gorenstein varieties (cf.\ \cite{book:bartoccibruzzohernandez2009}, Appendix~C). From now on, we   assume that $\mathscr{X}$ is a smooth projective stack of dimension $d$, so that it is of the form $[Z/G]$ with $Z$ a smooth quasi-projective variety (cf.\ Remark \ref{rem:equiv}). Recall that any $G$-equivariant coherent sheaf on $Z$ admits a finite resolution consisting of $G$-equivariant locally free sheaves of finite rank (\cite{book:chrissginzburg2010}, Prop.~5.1.28). Then we get the following result.
\begin{lemma}
A coherent sheaf on $\mathscr{X}$ admits a finite resolution by locally free sheaves of finite rank.
\end{lemma}

Before proving Serre duality theorem for Ext group we need some technical results about the relation between the derived dual $(\cdot)^\ast$ of a coherent sheaf and the tensor product $\stackrel{\mathbf L}{\otimes}$ in the derived category of $\mathscr{X}$. The techniques we shall use are similar to those in the proofs of \cite{book:bartoccibruzzohernandez2009}, Prop.~A.86, Prop.~A.87 and Cor.~A.88.
\begin{lemma}\label{eq:propA86}
Let $\mathcal{E}$, $\mathcal{F}$ and $\mathcal{G}$ be coherent sheaves on $\mathscr{X}$. There is a functorial isomorphism
\begin{equation}
\mathbf{R}\mathcal{H}om_{\mathscr{X}}^\bullet(\mathcal{E},\mathcal{F})
\stackrel{\mathbf L}{\otimes}\mathcal{G}\simeq  \mathbf{R}\mathcal{H}om_{\mathscr{X}}^\bullet(\mathcal{E},\mathcal{F}\stackrel{\mathbf L}{\otimes}\mathcal{G})
\end{equation}
in the derived category.
\end{lemma}
\proof
Let $\mathcal{E}^\bullet\to \mathcal{E}$ and $\mathcal{G}^\bullet\to \mathcal{G}$ be  finite resolutions of $\mathcal{E}$ and $\mathcal{G}$, respectively, consisting of locally free sheaves of finite rank. There is a quasi-isomorphism of complexes 
\begin{equation}\label{eq:homtensor}
\mathcal{H}om_{\mathscr{X}}^\bullet(\mathcal{E}^\bullet, \mathcal{F})\otimes \mathcal{G}^\bullet \simeq  \mathcal{H}om_{\mathscr{X}}^\bullet(\mathcal{E}^\bullet, \mathcal{F}\otimes \mathcal{G}^\bullet)\ .
\end{equation}
Let $\mathcal{F}\to \mathcal{F}^\bullet$ be an injective resolution of $\mathcal{F}$. Then $\mathcal{J}^\bullet=\mathcal{F}^\bullet\otimes\mathcal{G}^\bullet$ is injective and quasi-isomorphic to $\mathcal F\otimes \mathcal{G}^\bullet$. There is an induced quasi-isomorphism 
\begin{equation}
\mathcal{H}om_{\mathscr{X}}^\bullet(\mathcal{E}^\bullet, \mathcal{F}^\bullet)\otimes \mathcal{G}^\bullet \to   \mathcal{H}om_{\mathscr{X}}^\bullet(\mathcal{E}^\bullet, \mathcal{J}^\bullet)\ ,
\end{equation}
which yields in derived category the required isomorphism.
\endproof

\begin{lemma}
Let $\mathcal{E}$, $\mathcal{F}$ be coherent sheaves on $\mathscr{X}$ and $\mathcal{M}^\bullet$ a finite complex of locally free sheaves of finite rank. Then
\begin{equation}\label{eq:tensor}
\mathbf{R}\mathcal{H}om_{\mathscr{X}}^\bullet(\mathcal{E}\stackrel{\mathbf L}{\otimes}\mathcal{M}^\bullet,\mathcal{F})\simeq  \mathbf{R}\mathcal{H}om_{\mathscr{X}}^\bullet(\mathcal{E},\mathbf{R}\mathcal{H}om_{\mathscr{X}}^\bullet(\mathcal{M}^\bullet,\mathcal{F}))\ .
\end{equation}
\end{lemma}
\proof
Let $\mathcal I^\bullet$ be an injective resolution of $\mathcal F$. There is an isomorphism of complexes
\begin{equation}\label{eq:localhom}\mathcal{H}om^\bullet_{\mathcal X}(\mathcal E\otimes\mathcal M^\bullet,\mathcal I^\bullet) \simeq 
\mathcal{H}om^\bullet_{\mathcal X}(\mathcal E, \mathcal{H}om^\bullet_{\mathcal X}(\mathcal M^\bullet,\mathcal I^\bullet))\,.
\end{equation}
The left-hand side produces in derived category the object 
$\mathbf{R}\mathcal{H}om_{\mathscr{X}}^\bullet(\mathcal{E}\stackrel{\mathbf L}{\otimes}\mathcal{M}^\bullet,\mathcal{F})$.
To deal with the right-hand side, we note that since $\mathcal M^\bullet$ is flat and $\mathcal I^\bullet$ is injective,
the complex $\mathcal{H}om^\bullet_{\mathcal X}(\mathcal M^\bullet,\mathcal I^\bullet)$ is injective (and is quasi-isomorphic
to $\mathbf{R}\mathcal{H}om_{\mathscr{X}}^\bullet(\mathcal{M}^\bullet,\mathcal{F})$). Therefore the right-hand side of eq.~\eqref{eq:localhom} in derived category yields
$\mathbf{R}\mathcal{H}om_{\mathscr{X}}^\bullet(\mathcal{E},\mathbf{R}\mathcal{H}om_{\mathscr{X}}^\bullet(\mathcal{M}^\bullet,\mathcal{F}))$.
\endproof
\begin{proposition}\label{prop:tensor}
Let $\mathcal{E}$, $\mathcal{F}$ and $\mathcal{G}$ coherent sheaves on $\mathscr{X}$. Then in the derived category of $\mathscr{X}$ there are functorial isomorphisms
\begin{align}
 \operatorname{Hom}_{\mathrm{D}(\mathscr{X})}(\mathcal{E}\stackrel{\mathbf L}{\otimes}\mathcal{G}^{*\bullet},\mathcal{F})&\simeq \operatorname{Hom}_{\mathrm{D}(\mathscr{X})}(\mathcal{E},\mathcal{F}\stackrel{\mathbf L}{\otimes}\mathcal{G})\ ,\\
\operatorname{Hom}_{\mathrm{D}(\mathscr{X})}(\mathcal{E}\stackrel{\mathbf L}{\otimes}\mathcal{G},\mathcal{F})&\simeq \operatorname{Hom}_{\mathrm{D}(\mathscr{X})}(\mathcal{E},\mathcal{F}\stackrel{\mathbf L}{\otimes}\mathcal{G}^{*\bullet})\ ,
\end{align}
where $\mathcal{G}^{*\bullet}$ denotes the derived dual $\mathbf{R}\mathcal{H}om_{\mathscr{X}}^\bullet(\mathcal{G},\mathcal{O}_{\mathscr{X}})$ of $\mathcal{G}$.
\end{proposition}
\proof Since $\mathcal{G}$ admits a finite resolution consisting of finite rank locally free sheaves, its derived dual $\mathcal{G}^{*\bullet}$  is isomorphic, in the derived category of $\mathscr{X}$, to a finite complex consisting of finite rank locally free sheaves. By Lemma \ref{eq:propA86}, we get
\begin{equation}
\mathcal{G}^{*\bullet}\stackrel{\mathbf L}{\otimes}\mathcal{F}\simeq  \mathbf{R}\mathcal{H}om_{\mathscr{X}}(\mathcal{G},\mathcal{F})\ .
\end{equation}
By eq. \eqref{eq:tensor}, we have
\begin{align}
\mathbf{R}\mathcal{H}om_{\mathscr{X}}^\bullet(\mathcal{E}\stackrel{\mathbf L}{\otimes}\mathcal{G}^{*\bullet},\mathcal{F})
 &\simeq  \mathbf{R}\mathcal{H}om_{\mathscr{X}}^\bullet(\mathcal{E},\mathbf{R}\mathcal{H}om_{\mathscr{X}}^\bullet(\mathcal{G}^{*\bullet},\mathcal{F}))\\
 &\simeq  \mathbf{R}\mathcal{H}om_{\mathscr{X}}^\bullet(\mathcal{E},\mathcal{F}\stackrel{\mathbf L}{\otimes}\mathcal{G})\ .
\end{align}
By taking cohomology, we obtain
\begin{equation}
 \operatorname{Hom}_{\mathrm{D}(\mathscr{X})}(\mathcal{E}\stackrel{\mathbf L}{\otimes}\mathcal{G}^{*\bullet},\mathcal{F})\simeq \operatorname{Hom}_{\mathrm{D}(\mathscr{X})}(\mathcal{E},\mathcal{F}\stackrel{\mathbf L}{\otimes}\mathcal{G})\ .
\end{equation}
Similarly, we get the second formula of the statement.
\endproof
\begin{theorem}{(Serre duality - II)}
Let $p\colon \mathscr{X}\to \mathrm{Spec}(k)$ be a smooth projective stack of pure dimension $d$. Let $\mathcal{E}$ and $\mathcal{F}$ be coherent sheaves on $\mathscr{X}$. Then
\begin{equation}
\operatorname{Ext}^i(\mathcal{E}, \mathcal{G})\simeq  \operatorname{Ext}^{d-i}(\mathcal{G},\mathcal{E}\otimes \omega_{\mathscr{X}})^\vee\ ,
\end{equation} 
where $\omega_{\mathscr{X}}$ is the canonical line bundle of $\mathscr{X}.$
\end{theorem}
\proof
By \cite{art:nironi2008-II}, Thm.~2.22, the dualizing sheaf of $\mathscr{X}$ is the canonical line bundle $\omega_{\mathscr{X}}$. By applying the Formula \eqref{eq:duality} to the complexes $\mathcal{E}^{*\bullet}\stackrel{\mathbf L}{\otimes}\mathcal{G}$ and $\mathcal{O}_{\mathrm{Spec}(k)}$ we get
\begin{equation}
\mathbf Rp_\ast\mathbf{R}\mathcal{H}om_{\mathscr{X}}(\mathcal{E}^{*\bullet}\stackrel{\mathbf L}{\otimes}\mathcal{G}, \omega_{\mathscr{X}}[d])\simeq  \mathbf{R}\Gamma(\mathscr{X},\mathcal{E}^{*\bullet}\stackrel{\mathbf L}{\otimes}\mathcal{G})^\vee\ .
\end{equation}
By taking cohomology and by Proposition \ref{prop:tensor}, we obtain
\begin{equation}
\operatorname{Hom}_{\mathrm{D}(\mathscr{X})}(\mathcal{E}^{*\bullet}\stackrel{\mathbf L}{\otimes}\mathcal{G}, \omega_{\mathscr{X}}[d])\simeq  \operatorname{Hom}_{\mathrm{D}(\mathscr{X})}(\mathcal{G},\mathcal{E}\stackrel{\mathbf L}{\otimes} \omega_{\mathscr{X}}[d]) 
\end{equation}
in the left-hand-side, and 
\begin{equation} 
H^0(\mathbf{R}\Gamma(\mathscr{X},\mathcal{E}^{*\bullet}\stackrel{\mathbf L}{\otimes}\mathcal{G}))^\vee \simeq \operatorname{Hom}_{\mathrm{D}(\mathscr{X})}(\mathcal{E},\mathcal{G})^\vee
\end{equation}
in the right-hand side.
Therefore
\begin{multline}
\operatorname{Ext}^i(\mathcal{E},\mathcal{G})\simeq  \operatorname{Hom}_{\mathrm{D}(\mathscr{X})}(\mathcal{E},\mathcal{G}[i])\simeq \operatorname{Hom}_{\mathrm{D}(\mathscr{X})}(\mathcal{G},\mathcal{E}\stackrel{\mathbf L}{\otimes} \omega_{\mathscr{X}}[d-i])^\vee \\
\simeq  \operatorname{Ext}^{d-i}(\mathcal{G},\mathcal{E}\otimes \omega_{\mathscr{X}})^\vee\ .
\end{multline}
\endproof

\bigskip
\section{The dual of a coherent sheaf on a smooth projective stack}

In this section we introduce in the case of  stacks the notion of dual of a coherent sheaf as given in \cite{book:huybrechtslehn2010}, and some of its properties. This definition works well also with sheaves of nonmaximal dimension. In \cite{art:nironi2008} Nironi states without proofs the results we describe here.

Let $\mathscr{X}$ be a smooth projective stack of dimension $d$ and $\pi\colon\mathscr{X}\to X$ its coarse moduli scheme. Denote by $\omega_{\mathcal{X}}$ the canonical line bundle of $\mathscr{X}$ and fix a polarization $(\mathcal{G},\mathcal{O}_X(1))$ on $\mathscr{X}$. Since $\mathscr{X}$ is a smooth projective stack, for any \'etale presentation $(U,u)$ of $\mathscr{X}$, the scheme  $U$ is  a smooth and separated  of finite type over $k$, i.e., it is a smooth variety.

Let $\mathcal{E}$ be a coherent sheaf of dimension $n$ on  $\mathscr{X}$. The codimension of $\mathcal{E}$ is by definition $c=d-n$. We say that $\mathcal{E}$ satisfies the {\em generalized  Serre condition} $S_{k,c}$ for $k\geq 0$ if for all \'etale presentations $(U,u)$ of $\mathscr{X}$ and for all points $x\in\mathrm{supp}(\mathcal{E}_{U,u})$ the conditions  
\begin{equation}
S_{k,c}\colon \mathrm{depth}((\mathcal{E}_{U,u})_x)\geq \min\{k, \dim(\mathcal{O}_{U,x})-c\}
\end{equation}
hold. The condition $S_{0,c}$ is always satisfied. $S_{1,c}$ is equivalent to the purity of $\mathcal{E}_{U,u}$ for all \'etale presentations $(U,u)$ of $\mathscr{X}$ (cf.\ \cite{book:huybrechtslehn2010}, Sect.~1), hence it is equivalent to the purity of $\mathcal{E}$ by Remark \ref{rem:pureness}.
\begin{proposition}\label{prop:ext}
Let $\mathcal{E}$ be a coherent sheaf of dimension $n$ and codimension $c=d-n$.
\begin{itemize}\setlength{\itemsep}{2pt}
\item[(i)] The coherent sheaves $\mathcal{E}xt_{\mathscr{X}}^q(\mathcal{E},\omega_{\mathscr{X}})$ are supported on $\mathrm{supp}(\mathcal{E})$ and $\mathcal{E}xt_{\mathscr{X}}^q(\mathcal{E},\omega_{\mathscr{X}})=0$ for all $q<c$.
\item[(ii)] $\mathrm{codim}(\mathcal{E}xt_{\mathscr{X}}^c(\mathcal{E},\omega_{\mathscr{X}}))\geq c$. Moreover, the sheaf $\mathcal{E}$ satisfies the condition $S_{k,c}$ if and only if $\mathrm{codim}(\mathcal{E}xt_{\mathscr{X}}^q(\mathcal{E},\omega_{\mathscr{X}}))\geq q+k$ for all $q>c$.
\end{itemize}
\end{proposition}
\proof
The first statement in (i) is trivial. For the second, by \cite{art:nironi2008-II}, Lemma 2.28, one can choose an integer $m$ large enough to ensure that
\begin{equation}
H^0(\mathscr{X}, \mathcal{E}xt_{\mathscr{X}}^q(\mathcal{E}\otimes \mathcal{G}\otimes \pi^\ast(\mathcal{O}_X(-m)),\omega_{\mathscr{X}}))\simeq  \operatorname{Ext}^q(\mathcal{E}\otimes \mathcal{G}\otimes \pi^\ast(\mathcal{O}_X(-m)),\omega_{\mathscr{X}})\ .
\end{equation}
By Serre's duality (Theorem \ref{thm:serreduality-I}) we have
\begin{multline}
\operatorname{Ext}^q(\mathcal{E}\otimes \mathcal{G}\otimes \pi^\ast(\mathcal{O}_X(-m)),\omega_{\mathscr{X}})\simeq  H^{d-q}(\mathscr{X}, \mathcal{E}\otimes \mathcal{G}\otimes \pi^\ast(\mathcal{O}_X(-m))) \\ 
\simeq  H^{d-q}(X, \pi_\ast(\mathcal{E}\otimes \mathcal{G})\otimes \mathcal{O}_X(-m))\ .
\end{multline}
The coherent sheaf $\pi_\ast(\mathcal{E}\otimes \mathcal{G})\otimes \mathcal{O}_X(-m)$ ha dimension less or equal than $n$. Thus for $q<c$, the cohomology group $H^0(\mathscr{X}, \mathcal{E}xt_{\mathscr{X}}^q(\mathcal{E}\otimes \mathcal{G}\otimes \pi^\ast(\mathcal{O}_X(-m)),\omega_{\mathscr{X}}))$ vanishes. Since 
\begin{multline}
H^0(\mathscr{X}, \mathcal{E}xt_{\mathscr{X}}^q(\mathcal{E}\otimes \mathcal{G}\otimes \pi^\ast(\mathcal{O}_X(-m)),\omega_{\mathscr{X}})) \simeq \\  H^0(\mathscr{X}, \mathcal{E}xt_{\mathscr{X}}^q(\mathcal{E},\omega_{\mathscr{X}})\otimes \mathcal{G}^\vee\otimes \pi^\ast(\mathcal{O}_X(m))) 
 \simeq H^0(X, F_{\mathcal{G}}(\mathcal{E}xt_{\mathscr{X}}^q(\mathcal{E},\omega_{\mathscr{X}}))\otimes\mathcal{O}_X(m))\ , 
\end{multline}
also $H^0(X, F_{\mathcal{G}}(\mathcal{E}xt_{\mathscr{X}}^q(\mathcal{E},\omega_{\mathscr{X}}))\otimes\mathcal{O}_X(m))=0$ for $q<c$ and, by choosing $m$ sufficiently large,  $F_{\mathcal{G}}(\mathcal{E}xt_{\mathscr{X}}^q(\mathcal{E},\omega_{\mathscr{X}}))=0$ for $q< c$. By Lemma \ref{lem:support}, $\mathcal{E}xt_{\mathscr{X}}^q(\mathcal{E},\omega_{\mathscr{X}})=0$ for $q<c$ as well. 

The statements in (ii) can be proved  by resorting to \'etale presentations of $\mathscr X$  and using the same arguments as in \cite{book:huybrechtslehn2010}, Prop.~1.1.6; indeed, for any \'etale presentation $(U,u)$ of $\mathscr{X}$, the representative of $\mathcal{E}xt_{\mathscr{X}}^q(\mathcal{E},\omega_{\mathscr{X}})$ on $U$ is $\mathcal{E}xt_{U}^q(\mathcal{E}_{U,u},\omega_{U})$, where $\omega_{U}$ is the canonical line bundle of $U$.
\endproof
\begin{definition}
Let $\mathcal{E}$ be a coherent sheaf of dimension $n$ and $c=d-n$ its codimension. The dual sheaf of $\mathcal{E}$ is defined as $\mathcal{E}^D:=\mathcal{E}xt_{\mathscr{X}}^c(\mathcal{E},\omega_{\mathscr{X}})$.
\end{definition}
If $c=0$, then $\mathcal{E}^D\simeq \mathcal{E}^\vee \otimes\omega_{\mathscr{X}}$, where $\mathcal{E}^\vee$ is $\mathcal{H}om(\mathcal{E}, \mathcal{O}_{\mathscr{X}})$.
\begin{remark}
If the stack $\mathscr{X}$ is not smooth we could think of studying the dual sheaf by choosing a closed embedding $i\colon \mathscr{X}\to \mathscr{P}$ in a smooth projective stack $\mathscr{P}$ and using $i_\ast(\mathcal{E})^D = \mathcal{E}xt^e_{\mathscr{P}}(i_\ast(\mathcal{E}),\omega_{\mathscr{P}})$ where $e$ is the codimension of $i_\ast(\mathcal{E})$ in $\mathscr{P}$. This is reasonable because of \cite{art:nironi2008}, Lemma 6.8: let $\mathscr{X}$ be a smooth projective stack and $i\colon \mathscr{X}\to \mathscr{P}$ a closed embedding in a smooth projective stack $\mathscr{P}$. Let $c$ be the codimension of $\mathcal{E}$ in $\mathscr{X}$ and $e$ the codimension of $i_\ast(\mathcal{E})$ in $\mathscr{P}$. Then $i_\ast(\mathcal{E}xt^c_{\mathscr{X}}(\mathcal{E},\omega_{\mathscr{X}}))\simeq  \mathcal{E}xt^e_{\mathscr{P}}(i_\ast(\mathcal{E}),\omega_{\mathscr{P}})$. 
\end{remark}
By the same arguments as in \cite{book:huybrechtslehn2010}, Lemma 1.1.8, one has:
\begin{lemma}
There is a spectral sequence
\begin{equation}
E_2^{p,q}=\mathcal{E}xt_{\mathscr{X}}^p(\mathcal{E}xt_{\mathscr{X}}^{-q}(\mathcal{E},\omega_{\mathscr{X}}),\omega_{\mathscr{X}})\to \mathcal{E}\ .
\end{equation}
In particular, there is a natural morphism $\theta_\mathcal{E}\colon \mathcal{E}\to E_2^{c,-c}=\mathcal{E}^{DD}$.
\end{lemma}
\begin{definition}
A coherent sheaf $\mathcal{E}$ of codimension $c$ is called {\em reflexive} if $\theta_\mathcal{E}$ is an isomorphism.
\end{definition}
By using Proposition \ref{prop:ext} and the same arguments of the proof of \cite{book:huybrechtslehn2010}, Prop.~1.1.10, one can prove the following:
\begin{proposition}\label{prop:dual}
Let $\mathscr{X}$ be a smooth projective stack of dimension $d$ and $\mathcal{E}$ a coherent sheaf of codimension $c$. The following conditions are equivalent:
\begin{itemize}\setlength{\itemsep}{2pt}
\item[1)] $\mathcal{E}$ is pure,
\item[2)] $\mathrm{codim}(\mathcal{E}xt_{\mathscr{X}}^{q}(\mathcal{E},\omega_{\mathscr{X}}))\geq q+1$ for $q>c$,
\item[3)] $\mathcal{E}$ satisfies $S_{1,c}$,
\item[4)] $\theta_\mathcal{E}$ is injective.
\end{itemize}
Similarly, the following conditions are equivalent:
\begin{itemize}\setlength{\itemsep}{2pt}
\item[1')] $\mathcal{E}$ is reflexive,
\item[2')] $\mathcal{E}$ is the dual of a coherent sheaf of codimension $c$,
\item[3')] $\mathrm{codim}(\mathcal{E}xt_{\mathscr{X}}^{q}(\mathcal{E},\omega_{\mathscr{X}}))\geq q+2$ for $q>c$,
\item[4')] $\mathcal{E}$ satisfies $S_{2,c}$.
\end{itemize}
\end{proposition}
\begin{example}\label{ex:saturatedloc}
Let $\mathscr{X}$ be a two-dimensional smooth projective stack and $\mathcal{E}$ a reflexive sheaf of dimension two on $\mathscr{X}$. By Proposition \ref{prop:dual} for any $q>0$ the coherent sheaf $\mathcal{E}xt_{\mathscr{X}}^{q}(\mathcal{E},\mathcal{O}_{\mathscr{X}})$ is zero. Thus for any $q>0$ the coherent sheaf $\mathcal{E}xt_{\mathscr{X}}^{q}(\mathcal{E},\mathcal{O}_{\mathscr{X}})_{U,u}=\mathcal{E}xt_{U}^{q}(\mathcal{E}_{U,u},\mathcal{O}_{U})$ vanishes for any \'etale presentation $(U,u)$ of $\mathscr{X}$, so that  the homological dimension of $\mathcal{E}_{U,u}$ is zero, and $\mathcal{E}_{U,u}$ is a locally free sheaf. Therefore $\mathcal{E}$ is locally free. If $\mathcal{E}$ is not reflexive, but only torsion-free, by Proposition \ref{prop:dual} the coherent sheaf $\mathcal{E}xt_{\mathscr{X}}^{1}(\mathcal{E},\mathcal{O}_{\mathscr{X}})$ has dimension zero and $\mathcal{E}xt_{\mathscr{X}}^{2}(\mathcal{E},\mathcal{O}_{\mathscr{X}})=0$.  Now, if $\varphi\colon \mathcal{H}\to \mathcal{E}$ is a surjective morphism, and $\mathcal{H}$ locally free, it turns out that that $\mathcal{E}xt_{\mathscr{X}}^{1}(\ker(\varphi),\mathcal{O}_{\mathscr{X}})=\mathcal{E}xt_{\mathscr{X}}^{2}(\ker(\varphi),\mathcal{O}_{\mathscr{X}})=0$, so that $\ker(\varphi)$ is locally free. As a consequence, the saturated subsheaf of any locally free sheaf on $\mathscr{X}$ is locally free. 
\end{example}

\bigskip\section{An example: framed sheaves on stacky Hirzebruch surfaces}\label{Mattia}
\centerline{\em by Mattia Pedrini}
\par\medskip
In this Appendix theory developed in this paper is used to give a mathematically rigorous foundation to the ideas and the results in \cite{art:bruzzopoghossiantanzini2011}. In that paper the authors argued that their results about the structure and the Poincaré polynomial of the moduli spaces of framed sheaves on Hirzebruch surfaces $\mathbb{F}_p$ made sense in a more general setting, namely, they noticed that their framed sheaves $(\mathcal{E}, \phi_{\mathcal{E}})$ were allowed to have ``fractionary" first Chern class, i.e., $\operatorname{c}_1(\mathcal{E}) = \frac{m}{p}C$, where $C$ is the exceptional curve in $\mathbb{F}_p$. They also argued that this came from some hidden ``stacky" structure.

The ground field in this Appendix will be the complex field $\mathbb C$. We follow the procedure described in Section \ref{FramedsheavesToricorbifolds} to construct a stacky compactification $\mathscr{X}_p$ of the total space $\operatorname{Tot}(\mathcal{O}_{\mathbb{P}^1}(-p))$, which will be called the {\em $p$-th stacky Hirzebruch surface}. This turns out to be a toric orbifold with Deligne-Mumford torus $\mathbb{C}^\ast\times \mathbb{C}^\ast$ and with coarse moduli scheme $\mathbb{F}_p$. The stack $\mathscr{X}_p$ is obtained by adding to $\mathrm{Tot}(\mathcal{O}_{\mathbb{P}^1}(-p))$ a trivial gerbe $\mathscr{D}_\infty=\mathbb{P}^1\times [pt/\mu_p]$; this provides a rigorous construction of the stack $\mathcal{X}_p$ in \cite{art:bruzzopoghossiantanzini2011}. We use the results in \cite{art:fantechimannnironi2010, art:cadman2007, art:borisovchensmith2004, art:jiangtseng2010} to study the geometry of this orbifold, obtaining a full proof of Lemma 4.1 in \cite{art:bruzzopoghossiantanzini2011}. In particular, we introduce a {\em tautological class} $\omega$ which plays the role of $-\frac{1}{p}C$ in \cite{art:bruzzopoghossiantanzini2011}, and so gives a formal justification to the ``fractionary" first Chern classes of the framed sheaves.

In fact, after choosing a good framing sheaf $\mathcal{F}_\infty^{\,\vec{w}}$ on $\mathscr{D}_\infty$, we use Theorem \ref{thm:moduli-toricroot} to construct moduli spaces of $(\mathscr{D}_\infty,\mathcal{F}_\infty^{\,\vec{w}}\,)$-framed sheaves on $\mathscr{X}_p$, with fixed rank and discriminant and fixed first Chern class which is an integer multiple of $\omega$. The study of these moduli spaces allows us to give rigorous proofs of Propositions 3.1 and 3.2 of \cite{art:bruzzopoghossiantanzini2011} in the stacky case, namely, that the moduli spaces are smooth quasi-projective varieties, and the classification of their torus-fixed points.

Our notation somehow differs  from that of \cite{art:bruzzopoghossiantanzini2011}, for example, we denote by $E$, instead of $C$, the exceptional curve in $\mathbb{F}_p$, and by $D_\infty$ instead of $C_\infty$ the divisor at infinity.

\subsection{The root toric stack over $\mathbb{F}_p$}

Let $N\cong\mathbb{Z}^2$ be the lattice of $1$-parameter subgroups of the $2$-dimensional torus $T_t:=(\mathbb{C}^\ast)^2$. Fix a $\mathbb{Z}$-basis $\{e_1,e_2\}$ of $N$. Let $p\in\mathbb{N}$ and define the vectors $v_0:=e_2$, $v_1:=e_1$, $v_2:=pe_1-e_2$ and $v_\infty:=-e_1$. Define the fan $\Sigma_p\subset N_{\mathbb{Q}} := N \otimes_{\mathbb{Z}} \mathbb{Q}$ given by
\begin{align}
\Sigma_p(0) & :=  \left\{\{ 0 \} \right\}\ ,\\
\Sigma_p(1) & :=  \left\{ \rho_i:=\operatorname{Cone}(v_i)\vert i =0,1,2,\infty \right\}\ ,\\
\Sigma_p(2) & :=  \left\{ \sigma_i:=\operatorname{Cone}(v_{i-1}, v_i)\vert i=1,2 \right\} \cup \left\{ \sigma_{\infty,2}:=\operatorname{Cone}(v_2,v_\infty), \sigma_{\infty,0}:=\operatorname{Cone}(v_\infty, v_0) \right\}\ .
\end{align}
The fan $\Sigma_p$ defines the {\em $p$-th Hirzebruch surface} $\mathbb{F}_p$, which is a smooth projective toric surface, and is a projective compactification of the total space of the line bundle $\mathcal{O}_{\mathbb{P}^1}(-p)$ over $\mathbb{P}^1$. By the orbit-cone correspondence (\cite{book:coxlittleschenck2011}, Thm.~3.2.6) the rays $\rho_0, \rho_2$ correspond to the {\em fibres} $F_1,F_2$ with respect to the natural projection $\mathbb{F}_p\to \mathbb{P}^1$, $\rho_1$ corresponds to the {\em exceptional curve} $E$ and $\rho_\infty$ corresponds to the {\em divisor at infinity} $D_\infty$ such that $\mathbb{F}_p = \operatorname{Tot} (\mathcal{O}_{\mathbb{P}^1}(-p)) \cup D_\infty$ (cf. also \cite{book:beauville1996}, Chp.~IV). The Picard group $\operatorname{Pic}(\mathbb{F}_p)$ is freely generated over $\mathbb{Z}$ by $F_1$ and $E$, with the conditions $F_1=F_2$ and $D_\infty=E+pF_1$. The intersection product is determined by $F_1^2=0$, $E\cdot F_1 = 1$, $E^2 = -p$.

We define the stack $\mathscr X_p:=\sqrt[p]{D_\infty / \mathbb{F}_p}$ obtained performing a $p$-th root construction along the divisor $D_\infty$ and call it {\em $p$-th stacky Hirzebruch surface}. By \cite{art:fantechimannnironi2010}, Thm.~5.2, $\mathscr{X}_p$ is a two-dimensional toric orbifold with Deligne-Mumford torus $T_t$ and with coarse moduli scheme $\pi_p\colon \mathscr{X}_p \to \mathbb{F}_p$. Its {\em stacky fan} is $\mathbf{\Sigma}_p := (N,\Sigma_p, \beta)$ where the map $\beta\colon\mathbb{Z}^4 \to N$ is given by $\{v_0,v_1,v_2,pv_\infty\}$.

By computing the Gale dual\footnote{We refer to \cite{art:borisovchensmith2004}, Sect.~2, for an introduction to Gale duality.} of the map $\beta\colon \mathbb{Z}^4 \to N$ we have the quotient stack description of $\mathscr{X}_p$ as $[Z_{\Sigma_p} / G_{\mathbf{\Sigma}_p}]$, where $Z_{\Sigma_p}\subset\mathbb{C}^4$ is the union over the cones $\sigma\in\Sigma_p(2)$ of the open subsets
\begin{equation}
Z_\sigma := \{ x\in\mathbb{C}^4\vert x_i \neq 0 \mbox{ if }\rho_i \notin \sigma \} \subset \mathbb{C}^4\ .
\end{equation}
The group $G_{\mathbf{\Sigma}_p}$ is
\begin{equation}
G_{\mathbf{\Sigma}_p}:= \operatorname{Hom}_\mathbb{Z} (\operatorname{DG}(\beta),\mathbb{C}^\ast)\ ,
\end{equation}
where $\operatorname{DG}(\beta)$ is the cokernel of the map $\beta^\ast \colon \mathbb{Z}^2 \to \mathbb{Z}^4$, dual to the map $\beta$. Thus $\operatorname{DG}(\beta)\simeq \mathbb{Z}^2$. By applying $\operatorname{Hom}_\mathbb{Z}(\cdot,\mathbb{C}^\ast)$ to the quotient map $\beta^\vee \colon \mathbb{Z}^4 \to \operatorname{DG}(\beta)\simeq \mathbb{Z}^2$, we obtain an injective group morphism $i_{G_{\mathbf{\Sigma}_p}}\colon G_{\mathbf{\Sigma}_p} \to (\mathbb{C}^\ast)^4$. By composing $i_{G_{\mathbf{\Sigma}_p}}$ with the standard action of $(\mathbb{C}^\ast)^4$ on $\mathbb{C}^4$, we obtain the action of $G_{\mathbf{\Sigma}_p}\simeq(\mathbb{C}^\ast)^2$ on $Z_{\Sigma_p}\subset\mathbb{C}^4$:
\begin{equation}
(t_1,t_2)\cdot (z_1,z_2,z_3,z_4) = (t_1^p\, t_2^{-p}\, z_1, t_2\, z_2, t_1 z_3, t_2 z_4)
\end{equation}
for any $(t_1,t_2)\in G_{\mathbf{\Sigma}_p}$ and $(z_1,z_2,z_3,z_4)\in Z_{\Sigma_p}$.

Define the $1$-dimensional integral closed substacks 
\begin{equation}
\mathscr{F}_i:=\pi_p^{-1}(F_i)_{red}\ ,\quad\mathscr{E}:=\pi_p^{-1}(E)_{red}\ ,\quad\mathscr{D}_\infty:=\pi_p^{-1}(D_\infty)_{red}\ ,
\end{equation}
for $i=1,2$. These are Cartier and $T_t$-invariant. In addition, they are the irreducible component of the simple normal crossing divisor $\mathscr{X}_p\setminus T_t$ (cf. \cite{art:fantechimannnironi2010}, Cor.~5.4). By \cite{art:fantechimannnironi2010}, Remark 7.19, $\operatorname{Pic}(\mathscr{X}_p)$ is isomorphic to $\operatorname{DG}(\beta) \simeq \mathbb{Z}^2$ and, by the computation of the Gale dual of $\beta$, it is generated by $\mathscr{F}_1, \mathscr{D}_\infty$, with the relations $\mathscr{F}_1=\mathscr{F}_2$, $\mathscr{E} = p\mathscr{D}_\infty - p\mathscr{F}_1$. In addition, we have (cf. \cite{art:cadman2007}, Sect.~3.1)
\begin{align}
\pi_p^{\ast} \mathcal{O}_{\mathbb{F}_p}(F_i) & \simeq  \mathcal{O}_{\mathscr{X}_p}(\mathscr{F}_i)\ ,\\
\pi_p^{\ast} \mathcal{O}_{\mathbb{F}_p}(E) & \simeq  \mathcal{O}_{\mathscr{X}_p}(\mathscr{E})\ ,\\
\pi_p^{\ast} \mathcal{O}_{\mathbb{F}_p}(D_\infty) & \simeq  \mathcal{O}_{\mathscr{X}_p}(p \mathscr{D}_\infty)\ .
\end{align}

\subsubsection*{Characterization of $\mathscr{D}_\infty$} 
The divisor $\mathscr{D}_\infty$ is by construction (cf. \cite{art:cadman2007}, Example~2.4.3) isomorphic to the root stack $\sqrt[p]{\mathcal{O}_{\mathbb{F}_p}(D_\infty)_{\vert D_\infty} / D_\infty}$. By \cite{art:fantechimannnironi2010}, Thm.~6.25, $\mathscr{D}_\infty$ is a smooth toric Deligne-Mumford stack with Deligne-Mumford torus $T_t \times [pt/\mu_p]$ and with coarse moduli scheme $r_p\colon \mathscr{D}_\infty\to D_\infty \simeq \mathbb{P}^1$ and is an {\em essentially trivial\footnote{We refer to \cite{art:fantechimannnironi2010}, Sect.~6.1, for an introduction to (essentially trivial) banded gerbes.}} gerbe with banding group $\mu_p$ over its rigidification, which is $D_\infty$. Since $\mathcal{O}_{\mathbb{F}_p}(D_\infty)_{\vert D_\infty}\simeq \mathcal{O}_{\mathbb{P}^1}(p)$, $\mathscr{D}_\infty$ is actually a {\em trivial} gerbe, indeed by \cite{art:fantechimannnironi2010}, Rmk.~6.4, $\sqrt[p]{\mathcal{O}_{\mathbb{F}_p}(D_\infty)_{\vert D_\infty} / D_\infty}\simeq \sqrt[p]{\mathcal{O}_{\mathbb{P}^1}(p) / \mathbb{P}^1}\simeq \sqrt[p]{\mathcal{O}_{\mathbb{P}^1}/ \mathbb{P}^1}$ and the latter is the trivial $\mu_p$-gerbe $\mathbb{P}^1\times [pt/\mu_p]$ over $D_\infty\simeq \mathbb{P}^1$.

The stacky fan of $\mathscr{D}_\infty$ is $\mathbf{\Sigma}_p/\rho_\infty := (N(\rho_\infty), \Sigma_p/\rho_\infty, \beta(\rho_\infty))$, where $N(\rho_\infty):=N/p\mathbb{Z} v_\infty\simeq \mathbb{Z} \oplus \mathbb{Z}_p$, the quotient fan $\Sigma_p/\rho_\infty \subset N(\rho_\infty)_\mathbb{Q} \simeq \mathbb{Q}$ is given by the images of $v_0, v_2$ in $N(\rho_\infty)_\mathbb{Q}$, so it is
\begin{align}
\Sigma_p/\rho_\infty(0) & =  \left\{\{0\}\right\}\ ,\\
\Sigma_p/\rho_\infty(1) & =  \left\{ \operatorname{Cone}(1), \operatorname{Cone}(-1)\right\}\ .
\end{align}
The map $\beta(\rho_\infty)$ is given by the images of $v_0, v_2$ in $N(\rho_\infty)\simeq \mathbb{Z} \oplus \mathbb{Z}_p$, and so it is represented by the matrix
\begin{equation*}
\left(
\begin{array}{cc}
1 & -1 \\
0 & 0
\end{array}
\right)\ .
\end{equation*}

The computation of the Gale dual of the map $\beta(\rho_\infty)$ gives the quotient representation $\mathscr{D}_\infty \simeq [\mathbb{C}^2\setminus \{ 0\} / \mathbb{C}^\ast\times \mu_p]$, with the action given by $(t,\omega)\cdot(z_1,z_2)=(t\, z_1,t\, z_2)$. Moreover, one has $\operatorname{DG}(\beta(\rho_\infty)) \simeq \mathbb{Z} \oplus \mathbb{Z}_p$, so that the Picard group of $\mathscr{D}_\infty$ is generated over $\mathbb{Z}$ by the line bundles $\mathcal{L}_1, \mathcal{L}_2$ corresponding to the characters $\chi^{(1,0)}\colon (t,\omega) \mapsto t$, $\chi^{(0,1)}\colon (t,\omega) \mapsto \omega$, with $\mathcal{L}_2^{\otimes p}\simeq \mathcal{O}_{\mathscr{D}_\infty}$.

Define the $0$-dimensional integral closed substacks $P_0:=r_p^{-1}(0)_{red}$, $P_\infty:=r_p^{-1}(\infty)_{red}$. By the above computation of the Gale dual we have
\begin{equation}
\mathcal{O}_{\mathscr{D}_\infty}(P_0) \simeq \mathcal{L}_1  \simeq \mathcal{O}_{\mathscr{D}_\infty}(P_\infty)\ .
\end{equation}

The following result is proved by using the explicit descriptions of the Gale duals $\operatorname{DG}(\beta)$ and $\operatorname{DG}(\beta(\rho_\infty))$, and by following the technique described in \cite{art:bruzzopedrinisalaszabo2013}, Remark 2.20, to compute the restriction map $\operatorname{Pic}(\mathscr X_p)\to \operatorname{Pic}(\mathscr D_\infty)$.
\begin{proposition}\label{prop:restrictions}
The restriction of the line bundle $\mathcal{O}_{\mathscr{X}_p}(\mathscr{E})$ to $\mathscr{D}_\infty$ is trivial. Moreover the following relations hold:
\begin{equation}
\mathcal{O}_{\mathscr{X}_p}(\mathscr{F}_i)_{\vert \mathscr{D}_\infty} \simeq \mathcal{L}_1\qquad \mathcal{O}_{\mathscr{X}_p}(\mathscr{D}_\infty)_{\vert \mathscr{D}_\infty} \simeq \mathcal{L}_1\otimes \mathcal{L}_2\ .
\end{equation}
\end{proposition}
So far one gets the following chain of isomorphisms
\begin{equation}
r_p^\ast\mathcal{O}_{D_\infty}(1) \simeq r_p^\ast \mathcal{O}_{\mathbb{F}_p}(F_1)_{\vert D_\infty}\simeq \mathcal{O}_{\mathscr{X}_p}(\mathscr{F}_i)_{\vert \mathscr{D}_\infty}\simeq \mathcal{L}_1\simeq \mathcal{O}_{\mathscr{D}_\infty}(P_\infty)\ .
\end{equation}

\begin{lemma}\label{lem:degree-infty}
For any line bundle $\mathcal{L} \simeq \mathcal{L}_1^{\otimes a} \otimes \mathcal{L}_2^{\otimes b}$ on $\mathscr{D}_\infty$ with $a,b\in\mathbb{Z}$, we have
\begin{equation}
\int^{et}_{\mathscr{D}_\infty} \operatorname{c}_1(\mathcal{L}) = \frac{a}{p}\ .
\end{equation}
\end{lemma}
\proof 
First note that
\begin{equation}
\mathcal{L}^{\otimes p} \simeq \mathcal{L}_1^{\otimes ap} \simeq \mathcal{O}_{\mathscr{D}_\infty}(ap\, P_\infty)\ .
\end{equation}
By \cite{art:vistoli1989}, Proposition 6.1, the morphism $r_p$ induces an isomorphism of the rational Chow groups and by \cite{art:vistoli1989}, Example~6.7, we have ${r_p}_\ast[P_\infty]=\frac{1}{p}[\infty]$. Therefore
\begin{equation}
\int^{et}_{\mathscr{D}_\infty} \operatorname{c}_1(\mathcal{L}^{\otimes p}) = \int^{et}_{\mathscr{D}_\infty} \operatorname{c}_1( \mathcal{O}_{\mathscr{D}_\infty}(ap\,P_\infty) ) = \int^{et}_{D_\infty} {r_p}_\ast \operatorname{c}_1 (\mathcal{O}_{\mathbb{P}^1}(ap\,P_\infty)) = \int^{et}_{D_\infty}\frac{1}{p}\operatorname{c}_1(\mathcal{O}_{\mathbb{P}^1}(ap)) = a\ .
\end{equation}
The result follows.
\endproof

\subsubsection*{The tautological line bundle}
We can consider the $p$-th root of $\mathscr{E}$ in the following sense. Define  the {\em tautological class} $\omega:= \mathscr{F}_1 - \mathscr{D}_\infty$ in $\operatorname{Pic}(\mathscr{X}_p)$.  This class has the following properties:
\begin{equation}
p\omega = -\mathscr{E}\qquad\mbox{and}\qquad \omega^2=-\frac{1}{p}\ .
\end{equation}
Define the {\em tautological line bundle} $\mathcal{R} := \mathcal{O}_{\mathscr{X}_p}(\omega)$. By construction, $\mathcal{R}$, together with $\mathcal{O}_{\mathscr{X}_p}(\mathscr{D}_\infty)$, generates the Picard group $\operatorname{Pic}(\mathscr{X}_p)$. Moreover, $\mathcal{R}^{\otimes p} \simeq \mathcal{O}_{\mathscr{X}_p}(-\mathscr{E})$ and by Proposition \ref{prop:restrictions} we also have
\begin{equation*}
\mathcal{R}_{\vert \mathscr{D}_\infty} \simeq \mathcal{L}_2\ .
\end{equation*}

\subsection{Moduli spaces of $(\mathscr{D}_\infty,\mathcal{F}_{\infty}^{s,\vec{w}})$-framed sheaves on $\mathscr{X}_p$}

For $i=0, \ldots, p-1$ define the line bundles $\mathcal{O}_{\mathscr{D}_\infty}(i):= \mathcal{L}_2^{\otimes\, i}$. In addition, fix $\vec{w}:=(w_0, \ldots, w_{p-1})\in \mathbb{N}^{p}$ and define the {\em framing sheaf}
\begin{equation*}
\mathcal{F}_{\infty}^{\,\vec{w}}:=\oplus_{i=0}^{p-1}\mathcal{O}_{\mathscr{D}_\infty}(i)^{\oplus w_i}\ .
\end{equation*}
It is a locally free sheaf on $\mathscr{D}_\infty$ of orbifold rank $r:=\sum_{i=0}^{p-1} w_i$.

Let $(\mathcal{E}, \phi_{\mathcal{E}})$ be a $(\mathscr{D}_\infty,\mathscr{F}_\infty^{\,\vec{w}})$-framed sheaf on $\mathscr{X}_p$. First note that its orbifold rank $\operatorname{ork}(\mathcal{E})$ is $r$. Indeed, by Remark \ref{rem:orbirank}, $\operatorname{ork}(\mathcal{E})$ is the zero degree part of $\mathrm{ch}^{et}(\mathcal E)$, and since the K-theory groups $K(\mathscr{X}_p)$ and $K(\mathscr{D}_\infty)$ are both generated by line bundles (see \cite{art:borisovhorja2006}, Thm.~4.6), the latter is preserved under the restriction to $\mathscr{D}_\infty$; so by the isomorphism $\phi_{\mathcal{E}}$ we get $\operatorname{ork}(\mathcal{E})=\operatorname{ork}(\mathcal{F}_{\infty}^{\,\vec{w}})=r$.

Let $\det(\mathcal{E})\simeq \mathcal{R}^{\otimes u}\otimes \mathcal{O}_{\mathscr{X}_p}(u_\infty\,\mathscr{D}_\infty)$ be the determinant of $\mathcal{E}$ for some $u,u_\infty\in\mathbb{Z}$. Since $\det(\mathscr{F}_\infty^{\,\vec{w}})\simeq\det(\mathcal{E}_{\vert \mathscr{D}_\infty})$, we get
\begin{equation}
\otimes_{i=0}^{k-1}\mathcal{O}_{\mathscr{D}_\infty}(i)^{\otimes w_i}\simeq  \mathcal{R}^{\otimes u}_{\vert \mathscr{D}_\infty}\otimes \mathcal{O}_{\mathscr{X}_p}(u_\infty\,\mathscr{D}_\infty)_{\vert \mathscr{D}_\infty}\ .
\end{equation}
By Proposition \ref{prop:restrictions}, we have ${\mathcal{R}}_{\vert\mathscr{D}_\infty}\simeq\mathcal{O}_{\mathscr{D}_\infty}(1)$ and $\mathcal{O}_{\mathscr{X}_p}(\mathscr{D}_\infty)_{\vert\mathscr{D}_\infty}\simeq\mathcal{L}_1\otimes\mathcal{O}_{\mathscr{D}_\infty}(1)$, hence we get that $u_\infty=0$, and 
\begin{equation}\label{eq:condition-determinant}
u = \sum_{i=0}^{p-1}iw_i\bmod p \ .
\end{equation}

The divisor $D_\infty$ is a good framing divisor, indeed it is a nef Cartier divisor with $D_\infty^2=p$, and the locally free sheaf $\mathcal{F}_{\infty}^{\,\vec{w}}$ is a good framing sheaf because it is the direct sum of line bundles $\mathcal{O}_{\mathscr{D}_\infty}(i)$ such that $\int^{et}_{\mathscr{D}_\infty}\operatorname{c}_1(\mathcal{O}_{\mathscr{D}_\infty}(i))=0$ for $i=0, \ldots, p-1$. So as explained in Section \ref{FramedsheavesToricorbifolds} we can construct moduli spaces of $(\mathscr{D}_\infty, \mathcal{F}_{\infty}^{\,\vec{w}})$-framed sheaves on $\mathscr{X}_p$. In particular by Theorem \ref{thm:moduli-toricroot} there exists a fine moduli space $\mathrm{M}_{r,u,\Delta}(\mathscr{X}_p,\mathscr{D}_\infty,\mathcal{F}_\infty^{\,\vec{w}}\,)$, parameterizing isomorphism classes of $(\mathscr{D}_\infty,\mathcal{F}_\infty^{\,\vec{w}}\,)$-framed sheaves on $\mathscr{X}_p$ of orbifold rank $r$, determinant $\mathcal{R}^{\otimes u}$  and discriminant $\Delta$, such that $u$ satisfies the equation \eqref{eq:condition-determinant}; this moduli space is a quasi-projective scheme over $\mathbb{C}$.
Here the discriminant of a sheaf $\mathcal{E}$ on $\mathscr{X}_p$ is
\begin{equation}
\Delta(\mathcal{E})=\int^{et}_{\mathscr{X}_p}\left(\operatorname{c}_2(\mathcal{E})-\frac{\operatorname{ork}(\mathcal{E})-1}{2\operatorname{ork}(\mathcal{E})}\operatorname{c}_1^2(\mathcal{E})\right)\ .
\end{equation}

\begin{remark}
The Picard group of $\mathscr{X}_p$ is isomorphic to the second singular cohomology group of $\mathscr{X}_p$ with integral coefficients via the \emph{first Chern class} map (see \cite{art:iritani2009}, Sect.~3.1.2). So, fixing the determinant line bundle of a sheaf $\mathcal{E}$ is equivalent to fixing its first Chern class.
\end{remark}

\subsubsection*{Smoothness of the moduli spaces}

\begin{lemma}
For any $s\in \mathbb{Z}$ and $i=0, \ldots, p-1$ the pushforward ${r_p}_\ast\Big(\mathcal{L}_1^{\otimes \, s}\otimes\mathcal{O}_{\mathscr{D}_\infty}(i)\Big)$ is
\begin{equation*}
{r_p}_\ast\Big(\mathcal{L}_1^{\otimes \, s}\otimes\mathcal{O}_{\mathscr{D}_\infty}(i)\Big) = \left\{
\begin{array}{ll}
\mathcal{O}_{\mathbb{P}^1}(s)\qquad & \mbox{ if } i= 0\ ,\\
0 & \mbox{ otherwise}\ .
\end{array}
\right.
\end{equation*} 
\end{lemma}
\proof
Let $s\in \mathbb{Z}$ and $i=0, \ldots, p-1$. Recall that $r_p\colon \mathscr{D}_\infty \to D_\infty\cong \mathbb{P}^1$ is a trivial gerbe, so its banding group $\mu_p$ fits into the exact sequence
\begin{equation}
1\to \mu_p \xrightarrow{i} \mathbb{C}^\ast\times \mu_p\xrightarrow{q} \mathbb{C}^\ast\to 1\ ,
\end{equation}
where $i\colon \omega\mapsto (1,\omega)$ and $q\colon (t,\omega)\mapsto t$. Moreover, any coherent sheaf on $\mathscr{D}_\infty$ decomposes as a direct sum of eigensheaves with respect to the characters of $\mu_p$. The pushforward with respect to $r_p$ is nonzero only on $\mu_p$-invariant part of a coherent sheaf on $\mathscr{D}_\infty$. So, by the previous exact sequence,   the pushforward ${r_p}_\ast\Big(\mathcal{L}_1^{\otimes \, s}\otimes\mathcal{O}_{\mathscr{D}_\infty}(i)\Big)$ is nonzero if and only if $i=0$. For $i=0$, the projection formula, which holds for the rigidification morphism $r_p$ (cf.\ \cite{question:vistoli2013}), easily implies
\begin{equation}
{r_p}_\ast\Big(\mathcal{L}_1^{\otimes \, s}\otimes\mathcal{O}_{\mathscr{D}_\infty}(0)\Big) = {r_p}_\ast \mathcal{O}_{\mathscr{D}_\infty}(s\, P_\infty) = {r_p}_\ast r_p^\ast \mathcal{O}_{\mathbb{P}^1}(s) = \mathcal{O}_{\mathbb{P}^1}(s)\ .
\end{equation}
\endproof
\begin{corollary}\label{cor:vanishing}
For all negative integers $s$ we have 
\begin{equation}
H^0(\mathscr{D}_\infty,\mathcal{L}_1^{\otimes \, s}\otimes\mathcal{O}_{\mathscr{D}_\infty}(i))=H^0\Big(\mathbb{P}^1,{r_p}_\ast\Big(\mathcal{L}_1^{\otimes \, s}\otimes\mathcal{O}_{\mathscr{D}_\infty}(i)\Big)\Big)=0\ .
\end{equation}
\end{corollary}

Thanks to Corollary \ref{cor:vanishing}, we can argue as in the proof of Proposition~2.1 of \cite{art:gasparimliu2010} to get the following result. The proof involves Serre duality for stacks (see Appendix \ref{sec:serreduality}).
\begin{proposition}\label{prop:vanishing}
The Ext group $\operatorname{Ext}^i(\mathcal{E}',\mathcal{E}\otimes\mathcal{O}_{\mathscr{X}_p}(-\mathscr{D}_\infty))$ vanishes for $i=0,2$ and for any pairs of $(\mathscr{D}_\infty,\mathcal{F}_\infty^{\,\vec{w}}\,)$-framed sheaves $(\mathcal{E},\phi_{\mathcal{E}})$ and $(\mathcal{E}',\phi_{\mathcal{E}'})$ on $\mathscr{X}_p$. In addition, $H^i(\mathscr{X}_p, \mathcal{E}\otimes \mathcal{O}_{\mathscr{X}_p}(-\mathscr{D}_\infty))=0$ for $i=0,2$.
\end{proposition}

Eventually, by using Theorem \ref{thm:deformations} and Proposition \ref{prop:vanishing}, we can prove the smoothness of the moduli spaces
and compute its dimension. 
\begin{theorem}\label{thm:moduli} 
The moduli space $\mathrm{M}_{r,u,\Delta}(\mathscr{X}_p,\mathscr{D}_\infty,\mathcal{F}_\infty^{\,\vec{w}}\,)$ is a smooth quasi-projective variety.
\end{theorem}
\begin{corollary}
The dimension of $\mathrm{M}_{r,u,\Delta}(\mathscr{X}_p,\mathscr{D}_\infty,\mathcal{F}_\infty^{\,\vec{w}}\,)$ is 
\begin{equation}
\operatorname{dim}_{\mathbb{C}} \mathrm{M}_{r,u,\Delta}(\mathscr{X}_p,\mathscr{D}_\infty,\mathcal{F}_\infty^{\,\vec{w}}\,) = 2 r \Delta - \sum_{j=0}^{p-1} \vec{w}\cdot\vec{w}(j) \frac{j(j+p-2)}{2p}\ ,
\end{equation}
where by $\vec{w}(j)$ we denote the translated vectors $(w_j, \ldots, w_{p-1}, w_0, \ldots, w_{j-1})$.
\end{corollary}
\proof
By Theorem \ref{thm:deformations}, $\operatorname{dim}_{\mathbb{C}} \mathrm{M}_{r,u,\Delta}(\mathscr{X}_p,\mathscr{D}_\infty,\mathcal{F}_\infty^{\,\vec{w}}\,) = \operatorname{dim}_{\mathbb{C}} \operatorname{Ext}^1(\mathcal{E},\mathcal{E}\otimes \mathcal{O}_{\mathscr{X}_p}(-\mathscr{D}_\infty))$, where $\mathcal E$ is the underlying torsion free sheaf of a point $[(\mathcal{E},\phi_{\mathcal{E}})]$ of $\mathrm{M}_{r,u,\Delta}(\mathscr{X}_p,\mathscr{D}_\infty,\mathcal{F}_\infty^{\,\vec{w}}\,)$. By Proposition \ref{prop:vanishing}, $\operatorname{Ext}^i(\mathcal{E},\mathcal{E}\otimes \mathcal{O}_{\mathscr{X}_p}(-\mathscr{D}_\infty))=0$ for $i=1,2$. So
\begin{equation}
\operatorname{dim}_{\mathbb{C}} \mathrm{M}_{r,u,\Delta}(\mathscr{X}_p,\mathscr{D}_\infty,\mathcal{F}_\infty^{\,\vec{w}}\,) = -\chi(\mathcal{E},\mathcal{E}\otimes \mathcal{O}_{\mathscr{X}_p}(-\mathscr{D}_\infty))\ ,
\end{equation}
where $\chi(\mathcal F, \mathcal G)$ is the Euler characteristic of a pair $(\mathcal{F}, \mathcal{G})$ of coherent sheaves on $\mathscr X_p$. In the following, we apply a version of T\"oen-Riemann-Roch theorem for Euler characteristic of pairs as stated in \cite{art:bruzzopedrinisalaszabo2013}, App.~A.4. As pointed out there, since $(\mathcal{E},\phi_{\mathcal{E}})$ is a $(\mathscr{D}_\infty,\mathcal{F}_\infty^{\,\vec{w}}\,)$-framed sheaf, to compute the dimension of the moduli space is enough to assume that $\mathcal{E}$ is locally free and compute $\chi(\mathscr{X}_p,\mathcal{E}\otimes\mathcal{E}^\vee\otimes \mathcal{O}_{\mathscr{X}_p}(-\mathscr{D}_\infty))$.

To apply T\"oen-Riemann-Roch theorem, we need to compute the inertia stack of $\mathscr X_p$. It is easy to see that it is $\mathscr{I}(\mathscr{X}_p)=\mathscr{X}_p \sqcup \left( \bigsqcup_{i=1}^{p-1} \mathscr{D}_\infty^i \right)$. Thus we get
\begin{align}
\chi(\mathscr{X}_p,\mathcal{E}\otimes\mathcal{E}^\vee\otimes \mathcal{O}_{\mathscr{X}_p}(-\mathscr{D}_\infty)) &= \int_{\mathscr{I}(\mathscr{X}_p)}^{et} \frac{\operatorname{ch}^{rep}(\sigma^\ast(\mathcal{E}\otimes\mathcal{E}^\vee\otimes \mathcal{O}_{\mathscr{X}_p}(-\mathscr{D}_\infty)))}{\operatorname{ch}^{rep}(\lambda_{-1}\mathcal{N}^\vee)}\cdot\operatorname{Td}^{et}(\mathcal{T}_{\mathscr{I}(\mathscr{X}_p)})\\
&= A+ B\ ,
\end{align}
where $\mathcal{N}$ is the normal bundle to the local immersion $\sigma\colon\mathscr{I}(\mathscr{X}_p)\to\mathscr{X}_p$, and
\begin{align}
A & := \int_{\mathscr{X}_p}^{et} \operatorname{ch}^{et}(\mathcal{E}\otimes\mathcal{E}^\vee\otimes \mathcal{O}_{\mathscr{X}_p}(-\mathscr{D}_\infty)) \operatorname{Td}^{et}(\mathscr{X}_p)\ ,\\
B & :=  \sum_{i=1}^{p-1} \int_{\mathscr{D}_\infty^i}^{et} \frac{\operatorname{ch}^{rep}\left((\mathcal{E}\otimes\mathcal{E}^\vee\otimes \mathcal{O}_{\mathscr{X}_p}(-\mathscr{D}_\infty))_{\vert \mathscr{D}_\infty^i} \right)}{\operatorname{ch}^{rep} \left( (\lambda_{-1}(\mathcal{N}_{\mathscr{D}_\infty^i/\mathscr{X}_p}^\vee) \right)}\operatorname{Td}^{et}(\mathscr{D}_\infty)\ .
\end{align}
First we compute $A$:
\begin{equation}
A = -2r\Delta + r^2\left( \int_{\mathscr{X}_p}^{et}\operatorname{Td}_2(\mathscr{X}_p) + \int_{\mathscr{X}_p}^{et}\left( \frac{1}{2}[\mathscr{D}_\infty]^2 - [\mathscr{D}_\infty]\operatorname{Td}_1(\mathscr{X}_p) \right) \right)\ .
\end{equation}
To compute the first integral in the right-hand side of the equation we need the following two facts. First the canonical bundle of $\mathscr X_p$ is $\omega_{\mathscr X_p}\simeq\mathcal{O}_{\mathscr X_p} (-(\mathscr{E}+ \mathscr{F}_1 +\mathscr{F}_2+ \mathscr{D}_\infty))$. This can be seen as a generalization of the analogous result for toric varieties \cite{book:coxlittleschenck2011}, Thm.~8.2.3, see \cite{art:kawamata2006}. 
%First, $\operatorname{c}_1(\mathcal{T}_{\mathscr{X}_p}) = [\mathscr{E}] + [\mathscr{F}_1] + [\mathscr{F}_2] + [\mathscr{D}_\infty]$, and the intersection product on $\mathscr{X}_p$ is determined by $\mathscr{F}_1^2=0$, $\mathscr{D}_\infty^2 = \mathscr{D}_\infty \cdot \mathscr{F}_1 = \frac{1}{p}$. 
Secondly, by the decomposition of the inertia stack,
\begin{equation}
\int_{\mathscr{X}_p}^{et} \operatorname{c}_2(\mathcal{T}_{\mathscr{X}_p}) = \int_{\mathscr{I}(\mathscr{X}_p)}^{et} \operatorname{c}(\mathcal{T}_{\mathscr{I}(\mathscr{X}_p)}) - (p-1)\int_{\mathscr{D}_\infty}^{et}\operatorname{c}_1(\mathcal{T}_{\mathscr{D}_\infty})\ .
\end{equation}

Moreover, $\omega_{\mathscr{D}_\infty}\simeq \mathcal{O}_{\mathscr D_\infty}(-(P_0+P_\infty))$, thus by Proposition \ref{prop:restrictions} $\int_{\mathscr{D}_\infty}^{et}\operatorname{c}_1(\mathcal{T}_{\mathscr{D}_\infty})=\frac{2}{p}$. On the other hand by \cite{art:tseng2011}, Thm. 3.4, the term $\int_{\mathscr{I}(\mathscr{X}_p)}^{et} \operatorname{c}(\mathcal{T}_{\mathscr{I}(\mathscr{X}_p)})$ equals the Euler characteristic of $\mathbb{F}_p$, which in turn is the number of rays in its fan, that is, $4$. Summing up
\begin{equation}
\int_{\mathscr{X}_p}^{et}\operatorname{Td}_2(\mathscr{X}_p) = \frac{1}{2}\int_{\mathscr{X}_p}^{et}\operatorname{c}_1(\mathcal{T}_{\mathscr{X}_p})^2 + \frac{1}{12}\int_{\mathscr{X}_p}^{et}\operatorname{c}_2(\mathcal{T}_{\mathscr{X}_p}) = -\frac{p^2-6p-7}{12p}\ .
\end{equation}
To compute the second integral we use the adjunction formula (\cite{art:nironi2008-II}, Thm.~3.8) and Proposition \ref{prop:restrictions}
\begin{equation}
\int_{\mathscr{X}_p}^{et}\left( \frac{1}{2}[\mathscr{D}_\infty]^2 - [\mathscr{D}_\infty]\operatorname{Td}_1(\mathscr{X}_p) \right) = \frac{1}{2}\int_{\mathscr{D}_\infty}^{et}\operatorname{c}_1(\omega_{\mathscr{D}_\infty}) = -\frac{1}{p}
\end{equation}
obtaining
\begin{equation}
A = -2r\Delta -r^2\frac{p^2-6p+5}{12p}\ .
\end{equation}

To compute $B$, first note that $(\mathcal{E}\otimes\mathcal{E}^\vee\otimes \mathcal{O}_{\mathscr{X}_p}(-\mathscr{D}_\infty))_{\vert \mathscr{D}_\infty^i}\simeq \bigoplus_{j=0}^{p-1} \mathcal{L}_1^{\otimes\, -1}\otimes\mathcal{O}_{\mathscr{D}_\infty}(j-1)^{\oplus \vec{w}\cdot\vec{w}(j)}$. To evaluate $\operatorname{ch}^{rep}((\mathcal{E}\otimes\mathcal{E}^\vee\otimes \mathcal{O}_{\mathscr{X}_p}(-\mathscr{D}_\infty))_{\vert \mathscr{D}_\infty^i})$ we need to know its decomposition with respect to the action of $\mu_p$, which is the generic stabilizer of $\mathscr{D}_\infty$, given by multiplication by $\omega^i$ for a primitive root of unity $\omega$. Here the action of the generic stabilizer is given by the maps $\phi^i_p\colon \omega \in \mu_p \to (1,\omega^{i})\in\mathbb{C}^\ast \times\mu_p$. So we have
\begin{multline}
\operatorname{ch}^{rep}\left( (\mathcal{E}\otimes\mathcal{E}^\vee\otimes \mathcal{O}_{\mathscr{X}_p}(-\mathscr{D}_\infty))_{\vert \mathscr{D}_\infty^i} \right) = \\ \sum_{j=0}^{p-1}\vec{w}\cdot\vec{w}(j)\operatorname{ch}^{rep}\Big(\mathcal{L}_1^{\otimes\, -1}\otimes\mathcal{O}_{\mathscr{D}_\infty}(j-1)^{\oplus \vec{w}\cdot\vec{w}(j)}\Big) 
= \\  \sum_{j=0}^{p-1} \vec{w}\cdot\vec{w}(j) \, \omega^{i(j-1)}\,\Big(1+\operatorname{c}_1\Big(\mathcal{L}_1^{\otimes\, -1}\otimes\mathcal{O}_{\mathscr{D}_\infty}(j-1)^{\oplus \vec{w}\cdot\vec{w}(j)}\Big)\Big)\ .
\end{multline}
The normal bundle $\mathcal{N}_{\mathscr{D}_\infty/\mathscr{X}_p}$ is isomorphic to $\mathcal{O}_{\mathscr{X}_p}(\mathscr{D}_\infty)_{\vert \mathscr{D}_\infty} \simeq \mathcal{L}_1\otimes \mathcal{L}_2$, so for the denominator in $B$ we have
\begin{equation}
\operatorname{ch}^{rep}( \lambda_{-1}(\mathcal{N}_{\mathscr{D}_\infty^i/\mathscr{X}_p}^\vee) ) = 1 - \omega^{-i} + \omega^{-i}\operatorname{c}_1(\mathcal{L}_1\otimes \mathcal{L}_2)\ .
\end{equation}
Using again Proposition \ref{prop:restrictions}, we reduce the computation of $B$ to
\begin{equation}
B = -\frac{1}{p}\sum_{j=0}^{p-1} \vec{w}\cdot\vec{w}(j) \sum_{i=1}^{p-1} \frac{\omega^{i(j-2)}}{(1-\omega^{-i})^2}\ .
\end{equation}
By \cite{art:bruzzopedrinisalaszabo2013}, Lemma B.1 with $s=p+j-2$ for $j=0, \ldots, p-1$, the sum over $i$ is
\begin{equation}
\sum_{i=1}^{p-1} \frac{\omega^{-i(j+2)}}{(1-\omega^{-i})^2} = -\frac{p^2-6p+5}{12} - \frac{j(j+p-2)}{2}\ .
\end{equation}
This yields the claimed formula.
\endproof

\begin{remark} 
If the framing sheaf is trivial one has $w_0=r$ and $w_j=0$ for $j>0$, so that the dimension of the moduli space is $2r\,\Delta$
(cf.\ \cite{art:bruzzopoghossiantanzini2011}).
\end{remark}

\subsubsection*{The rank-one case}

Note that $\operatorname{Tot}(\mathcal{O}_{\mathbb{P}^1}(-p))$ is exactly $\mathscr{X}_p \setminus \mathscr{D}_\infty \simeq \mathbb{F}_p \setminus D_\infty$. Set  $T(-p):=\operatorname{Tot}(\mathcal{O}_{\mathbb{P}^1}(-p))$ and consider the Hilbert scheme $\hilb{n}{T(-p)}$ of $n$ points of $T(-p)$. We can define a morphism
\begin{equation}
i_{u,n}\colon \hilb{n}{T(-p)} \to \mathrm{M}_{1,u,n}(\mathscr{X}_p,\mathscr{D}_\infty,\mathcal{O}_{\mathscr{D}_\infty}(i))
\end{equation}
in the following way. Take a point $Z$ of $\hilb{n}{T(-p)}$ and consider its ideal sheaf $\mathcal{I}_Z$. The pushforward $i_\ast(\mathcal{I}_Z)$ with respect to the inclusion $i\colon T(-p) \to \mathscr{X}_p$ is a rank 1 torsion-free sheaf on $\mathscr{X}_p$. Its determinant $\det(i_\ast(\mathcal{I}_Z))$ is trivial and $\int^{et}_{\mathscr{X}_p}\operatorname{c}_2(i_\ast(\mathcal{I}_Z))=n$. Moreover $i_\ast(\mathcal{I}_Z)$ is locally free in a neighbourhood of $\mathscr{D}_\infty$ because $Z$ is disjoint from $\mathscr{D}_\infty$. Now take $u\in\mathbb{Z}$ and $i\in \{0, \ldots, p-1\}$ such that
\begin{equation}
i = u \bmod p \ .
\end{equation}
Then $\mathcal{E}:=i_\ast(\mathcal{I}_Z)\otimes \mathcal{R}^{\otimes u}$ is a rank one torsion-free sheaf on $\mathscr{X}_p$, locally free in a neighbourhood of $\mathscr{D}_\infty$, with a canonical framing induced by the isomorphism $\mathcal{R}^{\otimes u}_{\vert\mathscr{D}_\infty}\simeq \mathcal{O}_{\mathscr{D}_\infty}(i)$. So $(\mathcal{E}, \phi_{\mathcal{E}})$ is a $(\mathscr{D}_\infty,\mathcal{F}_\infty^{\,\vec{w}}\,)$-framed sheaf on $\mathscr{X}_p$, for $\vec{w}=(\delta_{ji})_{j=0, \ldots, p-1}$. Moreover, $\det(\mathcal{E})\simeq \mathcal{R}^{\otimes u}$ and 
\begin{equation}
\int^{et}_{\mathscr{X}_p} \operatorname{c}_2(\mathcal{E})=n \ .
\end{equation}
Thus $[(\mathcal{E},\phi_{\mathcal{E}})]$ gives a point in $\mathrm{M}_{1,u,n}(\mathscr{X}_p,\mathscr{D}_\infty,\mathcal{O}_{\mathscr{D}_\infty}(i))$.

It is easy to see that $i_{u,n}$ is injective, that the construction extends to families of zero-dimensional subschemes of $T(-p)$ of length $n$, and that $i_{u,n}$ admits an inverse. It follows that  $i_{u,n}$ is an isomorphism of fine moduli spaces.

\begin{remark}\label{rem:rankone}
We can sum up the previous construction in the following way. Fix $i\in \{0, 1, \ldots,$ $ p-1\}$, $u\in \mathbb{Z}$, such that $ u= i \bmod p$. For any $[(\mathcal{E}, \phi_{\mathcal{E}})]\in\mathrm{M}_{1,u,n}(\mathscr{X}_p,\mathscr{D}_\infty,\mathcal{O}_{\mathscr{D}_\infty}(i))$, there exists a zero-dimensional subscheme of $T(-p)$, with ideal sheaf $\mathcal{I}$, such that the torsion-free sheaf $\mathcal{E}$ isomorphic to $i_\ast(\mathcal{I})\otimes \mathcal{R}^{\otimes u}$. Moreover, the framing $\phi_{\mathcal{E}}$ is canonically induced by the isomorphism $\mathcal{R}^{\otimes u}_{\vert\mathscr{D}_\infty}\simeq \mathcal{O}_{\mathscr{D}_\infty}(i)$.
\end{remark}

\subsubsection*{Torus action and fixed points}
For the sake of completeness we analyze the torus action on the moduli space and characterize its fixed points.
Let $T_\rho$ be the maximal torus of $GL(r,\mathbb{C})$ consisting of diagonal matrices and set $T:=T_t\times T_\rho$. For any element $(t_1, t_2)\in T_t$ let $F_{(t_1, t_2)}$ be the automorphism of $\mathscr{X}_p$ induced by the torus action of $T_t$ on $\mathscr{X}_p$. Define an action of the torus $T$ on $\mathrm{M}_{r,u,\Delta}(\mathscr{X}_p,\mathscr{D}_\infty,\mathcal{F}_\infty^{\,\vec w})$ by
\begin{equation}
(t_1,t_2,\vec\rho\,)\cdot [(\mathcal{E},\phi_{\mathcal{E}})]
:=\big[\big( \big(F_{(t_1,t_2)}^{-1}\big)^\ast(\mathcal{E}), \phi_\mathcal{E}'\big)\big] \ ,
\end{equation}
where $\vec\rho=(\rho_1,\dots,\rho_r)\in T_\rho$ and $\phi_\mathcal{E}'$ is the composition of isomorphisms
\begin{equation}
\phi_\mathcal{E}' \colon  \big(F_{(t_1,t_2)}^{-1}\big)^\ast\mathcal{E}_{\vert \mathscr{D}_\infty} \xrightarrow{\big(F_{(t_1,t_2)}^{-1}\big)^\ast(\phi_{\mathcal{E}})} 
\big(F_{(t_1,t_2)}^{-1}\big)^\ast \mathcal{F}_\infty^{\,\vec w} \to \mathcal{F}_\infty^{\,\vec w} \xrightarrow{\vec\rho \ \cdot} \mathcal{F}_\infty^{\,\vec w} \ .
\end{equation}
Here the middle arrow is given by the $T_t$-equivariant structure induced on $\mathcal{F}_\infty^{\,\vec{w}}$ by the restriction of the torus action of $\mathscr{X}_p$ to $\mathscr D_\infty$.

\begin{proposition} 
A $T$-fixed point $[(\mathcal{E},\phi_\mathcal{E})]\in \mathrm{M}_{r,u,\Delta}(\mathscr{X}_p,\mathscr{D}_\infty,\mathcal{F}_\infty^{\,\vec w})^T$  decomposes as a direct sum of rank-one framed sheaves
\begin{equation}
(\mathcal{E}, \phi_{\mathcal{E}})=\bigoplus_{\alpha=1}^r (\mathcal{E}_\alpha, \phi_{\alpha})\ ,
\end{equation}
where for $i=0, \ldots, p-1$ and $\sum_{j=0}^{i-1}w_j<\alpha\leq\sum_{j=0}^i w_j$ one has:
\begin{itemize}
\item $\mathcal{E}_\alpha$ is a tensor product $i_\ast(\mathcal{I}_\alpha)\otimes\mathcal{R}^{\otimes u_\alpha}$, where $i\colon T(-p)\to \mathscr{X}_p$ is the inclusion morphism, $\mathcal{I}_\alpha$ is an ideal sheaf of zero-dimensional subscheme $Z_\alpha$ of  $T(-p)$ supported at the $T_t$-fixed points $P_1,P_2$ and $u_\alpha\in\mathbb{Z}$ is such that 
\begin{equation}\label{eq:decompositionalpha}
u_\alpha =  i\bmod p \ ;
\end{equation}
\item the framing ${\phi_\alpha}_{\vert \mathscr{D}_\infty}\colon \mathcal{E}_\alpha \xrightarrow{\sim}\mathcal{O}_{\mathscr{D}_\infty}(i)$ is induced by the isomorphism $\mathcal{R}^{\otimes\, u_\alpha}_{\vert \mathscr{D}_\infty}\simeq \mathcal{O}_{\mathscr{D}_\infty}(i)$.
\end{itemize}
\end{proposition}
\proof
One can argue along the lines of the proof of Proposition~3.2 in \cite{art:bruzzopoghossiantanzini2011} to obtain a decomposition $\mathcal{E}=\bigoplus_{\alpha=1}^r \mathcal{E}_\alpha$, where each $\mathcal{E}_\alpha$ is a $T$-invariant rank-one torsion-free sheaf on $\mathscr{X}_p$. The restriction ${\phi_\mathcal{E}}_{\vert\mathcal{E}_\alpha}$ gives a canonical framing to a direct summand of $\mathcal{F}_\infty^{\,\vec w}$. Reordering the indices $\alpha$, for $i=0,\ldots,p-1$ and for each $\alpha$ such that $\sum_{j=0}^{i-1}w_j<\alpha\leq\sum_{j=0}^i w_j$ we have an induced framing on $\mathcal{E}_\alpha$
\begin{equation}
\phi_\alpha:={\phi_\mathcal{E}}_{\vert\mathcal{E}_\alpha}\colon\mathcal{E}_\alpha\xrightarrow{\sim}\mathcal{O}_{\mathscr{D}_\infty}(i)\ .
\end{equation}
Thus $(\mathcal{E}_\alpha,\phi_\alpha)$ is a $(\mathscr{D}_\infty, \mathcal{O}_{\mathscr{D}_\infty}(i))$-framed sheaf of rank one on $\mathscr{X}_p$. As explained in Remark \ref{rem:rankone}, the torsion-free sheaf $\mathcal{E}_\alpha$ is a tensor product of the ideal sheaf $\mathcal{I}_\alpha$ of a zero-dimensional subscheme $Z_\alpha$ of length $n_\alpha$ whose support is contained in $T(-p)$, by a line bundle $\mathcal{R}^{\otimes u_\alpha}$ for a $u_\alpha\in \mathbb{Z}$ which satisfies Formula \eqref{eq:decompositionalpha}. Since the torsion-free sheaf $\mathcal{E}$ is fixed by the $T_t$-action, $Z_\alpha$ is fixed as well, hence it is supported at the $T_t$-fixed points $P_1,P_2$ of $T(-p)$.
\endproof

\begin{remark}
Let $[(\mathcal{E}, \phi_{\mathcal{E}})]=[\oplus_{\alpha=1}^r (\mathcal{E}_\alpha, \phi_{\alpha})]$ be a $T$-fixed point in $\mathrm{M}_{r,u,\Delta}(\mathscr{X}_p,\mathscr{D}_\infty,\mathcal{F}_\infty^{\,\vec w})$. Then 
\begin{equation}
\mathcal{R}^{\otimes u}\simeq \det(\mathcal{E})\simeq\otimes_{\alpha=1}^r \det(\mathcal{E}_\alpha)\simeq \otimes_{\alpha=1}^r \mathcal{R}^{\otimes u_\alpha}\ ,
\end{equation}
hence $\sum_{\alpha=1}^r u_\alpha=u$.
One can write $Z_\alpha=Z_\alpha^1\cup Z_\alpha^2$, where   $Z_\alpha^i$ is supported at the $T_t$-fixed point $P_i$ for $i=1,2$. Each $Z_\alpha^i$ corresponds to a Young tableau $Y_\alpha^i$ (cf.\ \cite{art:ellingsrudstromme1987}), so that  $Z_\alpha$ corresponds to a pair of Young tableaux $Y_\alpha=(Y_\alpha^1,Y_\alpha^2)$ such that $\vert Y_\alpha^1\vert + \vert Y_\alpha^2\vert= n_\alpha$. (A {\em Young tableau} is a finite set $Y\subset\mathbb{N}_{>0} \times\mathbb{N}_{>0}$, that we think of as sitting ``in the first quadrant,'' i.e., its elements are the coordinates of the right-top vertices of cells arranged in left-justified columns, with the columns lengths weakly decreasing (each column has the same or a shorter length than its predecessor.)

Thus we can denote the point $[(\mathcal{E}, \phi_{\mathcal{E}})]$ by a pair $(\vec{Y},\vec{u})$, where 
\begin{itemize}
\item $\vec{Y}=((Y_1^1,Y_1^2), \ldots, (Y_r^1,Y_r^2))$ and for any $\alpha=1, \ldots, r$ the pair $(Y_\alpha^1,Y_\alpha^2)$ is such that $\vert Y_\alpha^1\vert + \vert Y_\alpha^2\vert= n_\alpha$,
\item $\vec{u}=(u_1, \ldots, u_r)$ and for any $\alpha=1, \ldots, r$, $u_\alpha$ is such that the relation \eqref{eq:decompositionalpha} holds and $\sum_{\alpha=1}^r u_\alpha=u$.
\end{itemize}
\end{remark}

\end{document}